\documentclass[a4paper,12pt]{article}
\usepackage{ifstyle}
\usepackage[dvips]{graphicx}
\usepackage{epsfig} 
\usepackage[usenames]{color} 
\numberwithin{equation}{section}  

\usepackage[english]{babel}               
\usepackage[T1]{fontenc}
\usepackage[latin1]{inputenc}
\usepackage{amssymb}             
\usepackage{array}                 

\def\modintr#1{\textcolor{blue}{#1}}
\def\modcomm#1{\textcolor{magenta}{#1}}
\def\modsuppr#1{\textcolor{red}{#1}}

\def\modintr#1{#1}\def\modcomm#1{}\def\modsuppr#1{}

\def\Xint#1{\mathchoice
{\XXint\displaystyle\textstyle{#1}}%
{\XXint\textstyle\scriptstyle{#1}}%
{\XXint\scriptstyle\scriptscriptstyle{#1}}%
{\XXint\scriptscriptstyle\scriptscriptstyle{#1}}%
\!\int}
\def\XXint#1#2#3{{\setbox0=\hbox{$#1{#2#3}{\int}$}
\vcenter{\hbox{$#2#3$}}\kern-.5\wd0}}

\def\dashint{\Xint-}
\newtheorem{theorem}{Theorem}[section]
\newtheorem{claim}[theorem]{Claim}
\newtheorem{definition}[theorem]{Definition}
\newtheorem{corollary}[theorem]{Corollary}
\newtheorem{lemma}[theorem]{Lemma}

\newtheorem{conjecture}[theorem]{Conjecture}

\def\opn#1#2{\def#1{\operatorname{#2} } } 
\opn\id{id}
\opn\Rm{Rm}
\opn\Ric{Ric} 
\opn\SB{SB}
\opn\Rc{Rc}
\opn\Scal{Sc}
\opn\Tr{Tr}
\opn\Trac{Tr} 
\opn\dist{dist}
\opn\diam{Diam}
\opn\det{det} 
\opn\div{div}
\opn\Ker{Ker} 
\opn\Cap{Cap}
\opn\BT{BT}
\opn\ras{ras}
\opn\exp{exp}
\opn\Exp{Exp}
\opn\exph{exph}
\opn\Osc{Osc} 
\opn\Herm{Herm}
\opn\End{End} 
\opn\Psh{Psh} 
\opn\Sh{Sh} 
\opn\Supp{Supp} 
\opn\Hess{Hess} 
\opn\UB{UB}
\opn\Vol{Vol} 
\opn\Div{Div}
\opn\Exc{Exc}

\newcommand{\N}{\mathbb N}
  
\newcommand{\Q}{\mathbb Q}
\newcommand{\R}{\mathbb R}
\newcommand{\C}{\mathbb C}

\newcommand{\I}{\mathbb I}

\newcommand{\proj}{\mathbb P}
\newcommand{\contract}{\mathrel{\kern-1.5pt\vrule width6.0pt height0.4pt depth0pt
                \vrule width0.4pt height4.0pt depth0pt}}
\newcommand{\retract}{\mathrel{\kern-1.5pt\vrule width0.4pt height4.0pt depth0pt
                          \vrule width6.0pt height0.4pt depth0pt}}

\newdimen\boxrulethickness \boxrulethickness=.07em
\newdimen\Openboxwidth \Openboxwidth=.72em
\newcommand{\Openbox}{%
 \leavevmode
 \hbox{%
   \hfil\vrule width\boxrulethickness
   \vbox to\Openboxwidth{%
     \advance\Openboxwidth -2\boxrulethickness
     \hrule height \boxrulethickness width\Openboxwidth\vfil
     \hrule height\boxrulethickness}%
   \vrule width\boxrulethickness\hfil
 }}
 			  
\begin{document}
{\def\thefootnote{\relax}} 
\begin{center} 
%\Large{\bf{K\"ahler-Einstein metrics over singular K\"ahler spaces}}
%\Large{\bf{Degenerate K\"ahler-Einstein metrics over varieties of general type}}
\Large{\bf{Degenerate complex Monge-Amp\`ere equations over compact K\"ahler manifolds}}
\\
\vspace{0.4cm}
\large{Jean-Pierre Demailly and Nefton Pali}
\end{center} 
\begin{abstract}
We prove the existence and uniqueness of the solutions of some very general type of degenerate complex Monge-Amp\`ere equations, and investigate their regularity. This type of equations are precisely what is needed in order to construct K\"ahler-Einstein metrics over irreducible singular K\"ahler spaces  with ample or trivial canonical sheaf and singular K\"ahler-Einstein metrics over varieties of general type. 
\end{abstract}
\tableofcontents 
\section{Introduction}
In a celebrated paper \cite{Yau} published in 1978, Yau solved the Calabi conjecture. As is well known, the problem of prescribing the Ricci curvature can be formulated in terms of non-degenerate complex Monge-Amp\`ere equations.
\begin{theorem}{\bf (Yau).}
Let $X$ be a compact K\"ahler manifold of complex dimension $n$ and let $\chi$ be a  K\"ahler class. Then for any smooth density $v>0$ on $X$ such that $\int_Xv=\int_X\chi^n$, there exists a unique $($smooth$)$ K\"ahler metric $\omega\in \chi$ $($i.e.\ $\omega=\omega_0+i\partial\bar\partial\varphi$ with $\omega_0\in\chi$ $)$ such that $\omega^n=(\omega_0+i\partial\bar\partial\varphi)^n=v$.
\end{theorem}
Another breakthrough concerning the study of complex Monge-Amp\`ere equations was achieved by Bedford-Taylor \cite{Be-Te}. They initiated a new method for the study of very degenerate complex Monge-Amp\`ere equations. In fact, by combining  these results, Ko{\l}odziej \cite{Kol1} proved the existence of solutions for equations of type 
$$
(\omega+i\partial\bar\partial\varphi )^n=v\,,
$$
where $\omega$ a K\"ahler metric and $v\geq 0$ a density in $L^p$ or in some general Orlicz spaces. However, in various geometric applications, it is necessary to consider the case where $\omega$ is merely semipositive. This more difficult situation has been examined first by Tsuji \cite{Ts}, and his technique has been reconsidered in the recent works \cite{Ca-La}, \cite{Ti-Zha}, \cite{E-G-Z1} and \cite{Pau}.

In this paper we push further the techniques developed so far and we obtain some very general and sharp results on the existence, uniqueness and regularity of the solutions of degenerate complex Monge-Amp\`ere equations. In order to define the relevant concept of uniqueness of the solutions, we introduce a suitable subset of the space of closed $(1,1)$-currents, namely the domain of definition $\BT$ of the complex Monge-Amp\`ere operator ``in the sense of Bedford-Taylor'': a current $\Theta$ is in $\BT$ if the the successive exterior powers can be computed as 
$$
\Theta^{k+1}\;=\;i\partial\bar\partial(\varphi \Theta^k)\,,
$$ 
where $\varphi$ is a potential of $\Theta$ and $\varphi \Theta^k$ is locally of finite mass. Then for every pseudoeffective $(1,1)$-cohomology class~$\chi$, we prove a monotone convergence result for exterior powers of currents in the subset 
$$
\BT_\chi\;:=\;\BT\cap\;\chi\,.
$$

The uniqueness of the solutions of the degenerate complex Monge-Amp\`ere equations in a reasonable class of unbounded potentials has been a big issue and the object of intensive studies, see e.g.\ \cite{Ts}, \cite{Ti-Zha}, \cite{Blo1}, \cite{E-G-Z1}. In this direction, we introduce the subset 
$$
\BT^{\log}_{\chi}\;\subset \;\BT_\chi\,,
$$ 
of (closed positive) currents $T\in \BT_{\chi}$ which have a Monge-Amp\`ere product $T^n$ possessing an $L^1$-density such that 
$$
\int_X-\log(T^n/\Omega)\,\Omega<+\infty\,,
$$ 
for some smooth volume form $\Omega>0$. For example this is the case when the current $T^n$ possesses an $L^1$-density with complex analytic singularities (see Theorem \ref{HgRegMA}). We observe that the Ricci operator is well defined in the class $\BT^{\log}_{\chi}$.

In the last section we prove existence and fine regularity properties of the solutions of complex Monge-Amp\`ere equations with respect to a given degenerate metric $\omega\geq 0$, when the right hand side possesses an $L\log^{n+\varepsilon}L$-density or a density carrying complex analytic singularities (see Theorems \ref{ExistSol} and \ref{HgRegMA}). As a consequence of this results, we derive the following generalization of Yau's theorem.
\begin{theorem}\label{yau-tm}
Let $X$ be a compact K\"ahler manifold of complex dimension $n$ and let $\chi$ be a $(1,1)$-cohomology class admitting a smooth closed semipositive $(1,1)$-form $\omega$ such that $\int_X\omega^n>0$.
\\
\\
{\bf (A)} For any $L\log^{n+\varepsilon}L$-density $v\geq 0$, $\varepsilon>0$ such that $\int_Xv=\int_X\chi^n$, there exists a unique closed positive current $T\in \BT_{\chi}$ such that \hbox{$T^n=v$}. Moreover, this current possesses bounded local potentials over $X$ and continuous local potentials outside a complex analytic set $\Sigma_{\chi}\subset X$. This set depends only on the class $\chi$ and can be taken to be empty if the class $\chi$ is K\"ahler.
\\
\\
{\bf (B)}
In the special case of a density $v\geq 0$ possessing complex analytic singularities the current $T$ is also smooth outside the complex analytic subset  $\Sigma_{\chi}\cup Z(v)\subset X$, where $Z(v)$ is the set of zeros and poles of $v$. 
\end{theorem}
The type of complex Monge-Amp\`ere equation solved in Theorem \ref{HgRegMA} is precisely what is needed in order to construct K\"ahler-Einstein metrics over irreducible singular K\"ahler spaces  with ample or trivial canonical sheaf. It can be also used to construct singular K\"ahler-Einstein metrics over varieties of general type and to solve generalized equations of the form 
$$
\Ric(\omega)\;=\;-\,\lambda\,\omega+\rho\,, \qquad\lambda\geq 0\,.
$$ 

The relevant $L^{\infty}$-estimate needed in the proof of Theorem \ref{HgRegMA} (in the case related with K\"ahler-Einstein metrics) is obtained combining the $L^{\infty}$-estimate in Statement (A) of Theorem \ref{Kolo} with an important iteration method invented by Yau \cite{Yau} (see the Lemma \ref{YauIth}). The main issue here is that one can not use directly the maximum principle since the reference metric is degenerate.

The proof of our Laplacian estimate in Theorem \ref{HgRegMA}, which is obtained as a combination of the ideas of in \cite{Yau}, \cite{Ts}, \cite{Blo2}, provides in particular a drastic simplification of Yau's most general argument for complex Monge-Amp\`ere equations with degenerate right hand side. Moreover, it can be applied immediately to certain singular situations considered in \cite{Pau} and it reduces the Laplacian estimate in \cite{Pau} to \modintr{a quite simple consequence}\modsuppr{a triviality. Namely if $b_p\geq 0$, $p=1,...,n$ are real numbers then $-b_n+b_1+...+b_n\geq 0$.}\modcomm{\ (J'ai supprim\'e cette phrase car je pense qu'elle n'\'eclaire pas le lecteur en l'absence du contexte...)} \modintr{(however, one should point out that the argument in \cite{Pau} contains a gap due to the fact that the $L^p$-norm of the exponential $\exp(\psi_{1,\varepsilon}-\psi_{2,\varepsilon})$ of $\varepsilon$-regularized quasi-plurisubharmonic functions need not be uniformly bounded in $\varepsilon$ under the assumption that $\exp(\psi_1-\psi_2)$ is $L^p$, as our Lemma \ref{IntegLp} clearly shows if we do not choose carefully the constant $A$ there).}
Theorem \ref{HgRegMA} gives also some metric results for the geometry of varieties of general type. In this direction, we obtain the following results.
\begin{theorem}\label{GenTyp} 
Let $X$ be a smooth complex projective variety of general type. If the canonical bundle is nef, then there exists a unique closed positive current $\omega_{_E}\in \BT^{\log}_{\, 2\pi c_1(K_X)}$ solution of the Einstein equation 
\begin{eqnarray}\label{KE-equat}
\Ric(\omega_{_E})\;=\;-\;\omega_{_E}\,.
\end{eqnarray}
This current possesses bounded local potentials over $X$ and defines a smooth K\"ahler metric outside a complex analytic subset $\Sigma$, which is empty if and only if the canonical bundle is ample.
\end{theorem}
The existence part has been studied in \cite{Ts}, \cite{Ca-La} and \cite{Ti-Zha} by a K\"ahler-Ricci flow method. The importance of the uniqueness statement in Theorem \ref{GenTyp} is the following. If a current 
$$
\omega_{_E}\in \BT^{\log}_{\, 2\pi c_1(K_X)}
$$ 
satisfies the Einstein equation \eqref{KE-equat} then it has bounded local potentials. In the non nef case we obtain the following statement.
\begin{theorem}\label{G-GenTyp}
Let $X$ be a smooth variety of general type and let $\SB\subset\Sigma$ be respectively the stable and augmented stable base locus of the canonical bundle $K_X$. Then there exists a  closed positive current $\omega_{_{E}}\in 2\pi c_1(K_X)$ over~$X$, with locally bounded potentials over $X\smallsetminus SB$, solution of the Einstein equation \eqref{KE-equat}
over $X\smallsetminus \SB$, which restricts to a smooth $($non-degenerate$)$ K\"ahler-Einstein metric over $X\smallsetminus \Sigma$. If $\omega_{_{E}}$ has minimal singularities, then $\omega_{_{E}}$ is unique in the class of currents with minimal singularities in $2\pi c_1(K_X)$.
\end{theorem}
%It is easy to see that in this setting the K\"ahler-Einstein current $\omega_{_{E}}$ is not unique. However the uniqueness follows directly from the comparison principle in the case the potential of $\omega_{_{E}}$ is with minimal singularities.

Quite recently Tian and Ko{\l}odziej \cite{Ti-Ko} proved a very particular case of our  $L^{\infty}$-estimate. Their method, which is completely different, is based on an idea developed in \cite{De-Pa}. Our $L^{\infty}$-estimate allows us to completely solve the following conjecture of Tian stated in \cite{Ti-Ko}.
\begin{conjecture}\label{ti-conj}
Let $(X,\omega_X)$ be a polarized compact connected K\"ahler manifold of complex dimension $n$, let $(Y,\omega_Y)$ be a compact irreducible K\"ahler space of complex dimension $m\leq n$, let $\pi:X\rightarrow Y$ be a surjective holomorphic map and let  $0\leq f\in L\log^{n+\varepsilon}L(X,\omega^n_X)$, for some $\varepsilon>0$ such that $1=\int_Xf\omega^n_X$. Set $K_t:=\{\pi^*\omega_Y+t\omega_X\}^n>0$ for $t\in (0,1)$. Then the solutions of the complex Monge-Amp\`ere equations 
$$
(\pi^*\omega_Y+t\omega_X+i\partial\bar\partial \psi_t)^n=K_t\,f\,\omega^n_X\,,
$$
satisfy the uniform $L^{\infty}$-estimate $\Osc(\psi_t):=\sup_X\psi_t-\inf_X\psi_t\leq C<+\infty$ for all $t\in (0,1)$.
\end{conjecture}

The present manuscript expands and completes a paper accepted for publication in the
International Journal of Mathematics, which had to be shortened
in view of the length of the manuscript and of the demands of referees -
in particular it gives more details about the relation with the existing
litterature (see Appendix~C).

%%%%%%%%%%%%%%%%%%%%%%%%%%%%%%%%%%%%%%%%%%%%%%%%%%%%%%%%%%%%%%%%%%%%%%%%%%%%%%%%%%%%%%%%%%%%%%%%%%%%%%%%%%%%%%%%%%%%%%%%%%%%%%%%%%%%%%%%%%%%%%%%%%%%%%%%%%%%%%%%%%%%%%%%%%%%%%%%%%%%%%%%%%%%%%%%%%%%%%%%%%%%%%%%%%%%%%%%%%%%%%%%%%%%%%%%%%%%%%%%%%%%%%%%%%%%%%%%%%%%%%%%%%%%%%%%%%%%%%%%%%%%%%%%%%%%%%%%%%%%%%%%%%%%%%%%%%%%%%%%%%%%%%%%%%%%%%%%%%%%%%%%%%%%%%%%%%%%%%%%%%%%%%%%%%%%%%%%%%%%%%%%%%%%%%%%%%%%%%%%%%%%%%%%%%%%%%%%%%%%%%%%%%%%%%%%%%%%%%%%%%%%%%%%%%%%%%%%%%%%%%%%%%%%%%%%%%%%%%%%%%%%%%%%%%%%%%%%%%%%%%%%%%%%%%%%%%%%%%%%%%%%%%%%%%%%%%%%%%%%%%%%%%%
\section{General $L^{\infty}$-estimates for the solutions}
Let $X$ be a compact connected complex manifold of complex dimension~$n$ and let $\gamma$ be a closed real $(1,1)$-current with continuous local potentials or a closed positive $(1,1)$-current with bounded local potentials. Then to any distribution $\Psi$ on $X$ such that $\gamma+i\partial\bar\partial \Psi\geq 0$ we can associate a unique locally integrable and bounded from above function $\psi:X\rightarrow [-\infty,+\infty)$ such that the corresponding distribution coincides with $\Psi$ and such that for any continuous or plurisubharmonic local potential $h$ of $\gamma$ the function $h+\psi$ is plurisubharmonic.  The set of functions $\psi$ obtained in this way will be denoted by ${\cal P}_{\gamma}$. We set
${\cal P}^0_{\gamma}:=\{\psi\in {\cal P}_{\gamma}\,\mid\,\sup_X\psi=0 \}$. 
\begin{definition}
Let $X$ be a compact complex manifold of complex dimension~$n$. A closed positive $(1,1)$-current with bounded local potentials such that $\{\gamma\}^n:=\int_X\gamma^n>0$, will be called big.
\end{definition}
If $X$ is compact K\"ahler, one knows by \cite{De-Pa} that the class $\{\gamma\}$ is big if and only if it contains a K\"ahler current $T=\gamma+i\partial\bar\partial\psi\ge\varepsilon\omega$ (the inequality is in the sense of currents), for some K\"ahler metric $\omega$ on~$X$ and~$\varepsilon>0$. 
\\
\\
{\bf Basic facts about Orlicz spaces.}
Let $P:\R_{\geq 0}\rightarrow\R_{\geq 0}$, $P(0)=0$, be a convex increasing function and $\Omega>0$ be a smooth volume form over a manifold $M$ and let $X\subset M$ be a Borel set of $\Omega$-finite volume. According to \cite{Ra-Re} we introduce the vector space
\begin{eqnarray*}
L^{P}(X):=\left\{f:X\rightarrow\R\cup\{\pm \infty\}\,|\,\exists \lambda>0\,:\,\int_XP(|f|/\lambda)\,\Omega<+\infty\right\}\,, 
\end{eqnarray*}
(with the usual identification of functions equal a.e.), equipped with the norm
\begin{eqnarray*}
\Vert f \Vert_{L^{P}(X)}:=\inf\left\{\lambda>0\,|\, \int_XP(|f|/\lambda)\,\Omega\leq 1\right\} \,.
\end{eqnarray*}
The space $L^{P}(X)$ equipped with this norm is called the Orlicz space associated with the convex function~$P$. Moreover this norm is order preserving, i.e 
$$
\Vert f \Vert_{L^{P}(X)}\leq \Vert g \Vert_{L^{P}(X)}\,,
$$ 
if $|f|\leq |g|$ a.e.
If $P(t)=|t|^p$, $p\geq 1$, then $L^{P}(X)$ is the usual Lebesgue space. More refined examples of Orlicz spaces are given by the functions 
$$
P_{\beta}:=t\log^{\beta}(e+t)\,,
$$ 
and 
$$
Q_{\beta}:=e^{t^{1/\beta}}-1\,,
$$
with $\beta\geq 1\,$. In these cases, we set 
$$
L\log^{\beta}L(X):=L^{P_{\beta}}(X)\,,
$$ 
and 
$$
\Exp^{1/\beta}L(X):=L^{Q_{\beta}}(X)\,.
$$
An important class of Orlicz spaces is given by considering functions $P$ satisfying the ``doubling property'': $P(2t)\leq 2^CP(t)$ for some constant $C\geq 1$. This is the case of the functions $|t|^p$ and $P_{\beta}(t)$, but not the case of $Q_\beta(t)$. For functions satisfying the doubling condition one has (see proposition 6 page 77 in \cite{Ra-Re})
\begin{eqnarray*}
L^{P}(X)=\left\{f:X\rightarrow\R\cup\{\pm \infty\}\,|\,\int_XP(|f|)\,\Omega<+\infty\right\}\,,
\end{eqnarray*}
and
\begin{eqnarray*}
\int_XP(\Vert f \Vert_{L^{P}(X)}^{-1}|f|)\,\Omega=1
\end{eqnarray*}
for all $f\in L^{P}(X)\smallsetminus \{0\}$. So in the particular case of the function $P_{\beta}$, one obtains the inequality
\begin{eqnarray}\label{LlogL-Int}
\|f\|_{L\log^{\beta}L(X)}\leq \int_X|f|\log^{\beta}\left(e+\|f\|_{L^1(X)}^{-1}|f|\right)\,\Omega\,,
\end{eqnarray}
since $\|f\|_{L^1(X)}\leq \|f\|_{L\log^{\beta}L(X)}$. It is quite hard to get precise estimates of the norm $\Exp^{1/\beta}L(X)$, however it is easy to see that
\begin{eqnarray}\label{Exp1-Vol}
\|1\|_{\Exp^{1/\beta}L(X)} \,=\,\frac{1}{\log^{\beta}(1+1/\Vol_{\Omega}(X))}\,.
\end{eqnarray}
The relation between the Orlicz spaces $L\log^{\beta}L(X)$ and $\Exp^{1/\beta}L(X)$ is expressed by the H\"{o}lder inequality (see \cite{Iw-Ma})
\begin{eqnarray}\label{Hold-Ineq}
\left|\int_X fg\,\Omega\right|\leq 2C_{\beta}\, \|f\|_{L\log^{\beta}L(X)}\,\|g\|_{\Exp^{1/\beta}L(X)} \,,
\end{eqnarray}
which follows from the inequality $xy\leq C_{\beta}(P_{\beta}(x)+Q_{\beta}(y))$ for all $x,y\geq 0$. (Observe moreover that $C_1=1$.)
\\
\\
We define the oscillation operator $\Osc:=\sup-\inf$. With the notations so far introduced we state the following result.
\begin{theorem}\label{Kolo}
Let $X$ be a compact connected K\"ahler manifold of complex dimension~$n$, let $\Omega>0$ be a smooth volume form, let $\gamma$ be a big closed positive $(1,1)$-current with continuous local potentials. Let also
$\psi\in{\cal P}_{\gamma}\cap L^{\infty}(X)$ be a solution of the degenerate complex Monge-Amp\`ere equation 
$$
(\gamma+i\partial\bar\partial \psi)^n=f\,\Omega\,,
$$
with $f\in L\log^{n+\varepsilon_0}L(X)$ for some $\varepsilon_0>0$.
Then the following conclusions hold.
\\
{\bf (A)}
There exists a uniform constant $C_1=C_1(\varepsilon_0,\gamma,\Omega)>0$ such that for all $\varepsilon\in (0,\varepsilon_0]$ we have an estimate
$$
\Osc(\psi)\leq (C_1/\varepsilon)^{n^2/\varepsilon}\,
I_{\gamma,\varepsilon}(f)^{\frac{n}{\varepsilon}}+\;1\,,
$$
where
$$
I_{\gamma,\varepsilon}(f):= \{\gamma\}^{-n}\int\limits_Xf\log^{n+\varepsilon}\left(e+\{\gamma\}^{-n}f\right)\Omega\,.
$$
{\bf (B)} Assume that the solution $\psi$ is normalized by the condition $\sup_X \psi=0$ and consider also a solution $\varphi\in {\cal P}_{\gamma}\cap L^{\infty}(X)$, $\sup_X \varphi=0$ of the degenerate complex Monge-Amp\`ere equation 
$$
(\gamma+i\partial\bar\partial \varphi)^n=g\,\Omega\,,
$$
with $g\in L\log^{n+\varepsilon_0}L(X)$. Assume also  $I_{\gamma,\varepsilon_0}(f), I_{\gamma,\varepsilon_0}(g)\leq K$ for some constant $K>0$. Then there exists a constant $C_2=C_2(\varepsilon_0,\gamma,\Omega, K)>0$ such that 
\begin{eqnarray*}
\|\varphi-\psi \|_{_{L^{\infty}(X)}}
&\leq&
2C_2^{^{\alpha_0}}\, \left(\log 
\|\varphi-\psi \|^{-1}_{_{ L^1(X,\,\Omega)}}\right)^{^{-\alpha_0}}\,,
\\
\\
\alpha_0&:=&\frac{1}{(n+1+n^2/\varepsilon_0)}\,,
\end{eqnarray*}
provided that the inequality 
$\|\varphi-\psi \|_{_{ L^1(X,\,\Omega)}}\leq \min\{1/2,e^{-C_2}\}$ holds.
\\
{\bf (C)} Let $(\gamma_t)_{t>0}$ be a family of currents satisfying the same properties as $\gamma$, fix a finite covering $(U_{\alpha})_{\alpha}$ of coordinate starshaped open sets, and let us write $\gamma_t=i\partial\bar\partial h_{t,\alpha}$ with $h_{t,\alpha}$ plurisubharmonic over $U_{\alpha}$, normalised by $\sup_{U_{\alpha}}h_{t,\alpha}=0$ and let $C_{1,t}:=C_1(\varepsilon_0, \gamma_t, \Omega)$, $C_{2,t}=C_2(\varepsilon_0,\gamma_t,\Omega, K)$.
Assume
\\
{\bf (C1)} $\sup_{t>0} \max_{\alpha}\|h_{t,\alpha}\|_{L^{\infty}(U_{\alpha})}<+\infty$ and
\\
{\bf (C2a)} there exist a decomposition of the type $\gamma_t=\theta_t+i\partial\bar\partial u_t$, with $\theta_t$ smooth, $\min_X u_t=0$, $\sup_{t>0}\max_X u_t<+\infty$ and $\theta_t\leq (\{\gamma_t\}^n)^{1/n}\omega$ for some K\"ahler metric $\omega>0$ on $X$,
\\
or
\\
{\bf (C2b)} the distributions $\gamma^n_t/\Omega$ are represented by $L^1$-functions and
$$
\sup_{t>0}\;\;\{\gamma_t\}^{-n}\int\limits_X\log\left(e+\{\gamma_t\}^{-n}\gamma^n_t/\Omega\right)\gamma^n_t<+\infty\,.
$$
Then $\sup_{t>0}C_{j,t}<+\infty$ for $j=1,2$.
\end{theorem}
Statement (C) will follow from the arguments of the proof of Statements (A) and (B) of Theorem \ref{Kolo}.

We start by proving a few basic facts about pluripotential theory, in a way which is best adapted for the understanding of the proof of the theorem \ref{Kolo}. The reader can also consult and compare with the related results in \cite{Be-Te}, \cite{Dem1}, \cite{Dem2}, \cite{G-Z} and \cite{Sic}.
	
Let $X$ be a compact complex manifold of complex dimension $n$, let $\gamma$ be a big closed positive $(1,1)$-current with bounded local potentials. Set
$$
{\cal P}_{\gamma}[0,1]:=\{\varphi\in {\cal P}_{\gamma}\,\mid\,0\leq \varphi \leq 1 \}\,,
$$ 
$\gamma_{\varphi}:=\gamma+i\partial\bar\partial\varphi$ and 
$$
\Cap_{\gamma}(E):=\sup_{\varphi\in {\cal P}_{\gamma}[0,1]}\,\{\gamma\}^{-n}\int\limits_E \gamma^n_{\varphi}\,,
$$
for all Borel sets $E\subset X$. We remark that if $(E_j)_j$, $E_j\subset E_{j+1}\subset X$ is a family of Borel sets and $E=\bigcup_jE_j$ then clearly, we have
\begin{eqnarray}\label{CapProp}
\Cap_{\gamma}(E)=\lim_{j\rightarrow+\infty}\Cap_{\gamma}(E_j)\,.
\end{eqnarray}

\begin{lemma}\label{IntBndPSh}
Let $X$ be a compact connected complex manifold of complex dimension $n$, let $\gamma$ be a closed real $(1,1)$-current with continuous local potentials or a closed positive $(1,1)$-current with bounded local potentials and let $\Omega>0$ be a smooth volume form. Then there exist constants $\alpha=\alpha(\gamma,\Omega)>0$,  $C=C(\gamma,\Omega)>0$ such that $\int_X-\psi\,\Omega\leq C$ and $\int_X e^{-\alpha\psi}\,\Omega\leq C$ for all $\psi\in {\cal P}^0_{\gamma}$.
\end{lemma}
(We notice that the first inequality follows from the second one.) 
The first two integral estimates of Lemma \ref{IntBndPSh} are quite standard in the ele\-mentary theory of plurisubharmonic functions and the dependence of the constants $\alpha$ and $C$  on $\gamma$ is only on the $L^{\infty}$ bound of its local potentials (see e.g.\ \cite{Hor} and \cite{Skoda}). To be more precise in sight of the uniform estimate $\int_X e^{-\alpha\psi}\,\Omega\leq C$ one can make the constant $\alpha$ depending only on the cohomology class of $\gamma$ as in \cite{Ti1}, but in this case the constant $C$ will depend on the $L^{\infty}$ bound of the local potentials of $\gamma$ and on the volume form $\Omega$. One can also make $C$ depending only on the volume form $\Omega$, but in this case $\alpha$ will depend on the $L^{\infty}$ bound of the local potentials of $\gamma$ and on the volume form $\Omega$.

The following lemma is the key technical tool which allows to deduce Statement (C) of Theorem \ref{Kolo}.
%CHANGE-REF\^A{\S}\^A{\S}\^A{\S}\^A{\S}\^A{\S}\^A{\S}\^A{\S}\^A{\S}\^A{\S}\^A{\S}\^A{\S}\^A{\S}\^A{\S}\^A{\S}\^A{\S}\^A{\S}\^A{\S}\^A{\S}\^A{\S}\^A{\S}
%CHANGE-REF\^A{\S}\^A{\S}\^A{\S}\^A{\S}\^A{\S}\^A{\S}\^A{\S}\^A{\S}\^A{\S}\^A{\S}\^A{\S}\^A{\S}\^A{\S}\^A{\S}\^A{\S}\^A{\S}\^A{\S}\^A{\S}\^A{\S}\^A{\S}
\begin{lemma} \label{CapBndPSh}%CHANGE-REF\^A{\S}\^A{\S}\^A{\S}\^A{\S}\^A{\S}\^A{\S}\^A{\S}\^A{\S}\^A{\S}\^A{\S}\^A{\S}\^A{\S}\^A{\S}\^A{\S}\^A{\S}\^A{\S}\^A{\S}\^A{\S}\^A{\S}\^A{\S}
Let $X$ be a compact connected K\"{a}hler manifold of complex dimension $n$ and let $\gamma$ be a big closed positive  $(1,1)$-current with continuous local potentials.
\\
{\bf (A)} There exists a constant $C=C(\gamma)>0$ such that $\Cap_{\gamma}(\{\psi<-t\})\leq C/t$ for all $\psi\in {\cal P}^0_{\gamma}$ and $t>0$. Moreover the constant $C$ stays bounded for perturbations of $\gamma$ satisfying the hypothesis {\rm (C1)} and {\rm (C2a)} of Statement {\rm (C)} in Theorem \ref{Kolo}.
\\
{\bf (B)} If $\gamma^n/\Omega\in L\log L(X)$, for a smooth volume form $\Omega>0$ then the conclusion of Statement {\rm(A)} holds with a constant $C=C(\gamma,\Omega)>0$ which stays bounded for perturbations of $\gamma$ satisfying the hypothesis {\rm(C1)} and {\rm(C2b)} of Statement {\rm(C)} in Theorem  \ref{Kolo}.
\end{lemma}
\emph{Proof}. 
We first notice the obvious inequality
$$
\int\limits_{\psi<-t}\gamma^n_{\varphi} \;\leq \frac{1}{t}\int\limits_X -\psi\, \gamma^n_{\varphi}
$$
which implies
\begin{eqnarray}\label{ParCap0}
\Cap_\gamma(\{\psi<-t\})\leq\frac{1}{t}
\sup_{\varphi\in{\cal P}_{\gamma}[0,1]}\{\gamma\}^{-n}
\int\limits_X -\psi\, \gamma^n_{\varphi}\,,
\end{eqnarray}
and we prove the following elementary claim.
\begin{claim}\label{CapIest}
Let $\gamma$ be a closed positive  $(1,1)$-current with bounded local potentials over a compact complex manifold $X$ of complex dimension $n$ and let $\varphi,\,\psi\in {\cal P}_{\gamma}$ such that $0\leq \varphi\leq 1$ and $\psi\leq 0$. Then
\begin{eqnarray}\label{ParCap1}
\int\limits_X -\psi\,\gamma^n_{\varphi} \;\leq \int\limits_X -\psi\,\gamma^n+n\int\limits_X\,\gamma^n\,.
\end{eqnarray}
\end{claim}
\emph{Proof}. The fact that the current $\gamma$ is positive implies $\psi_c:=\max\{\psi,c\}\in {\cal P}_{\gamma}$, $c\in \R_{<0}$, so by the monotone convergence theorem it is sufficient to prove inequality \eqref{ParCap1} for  $\psi\in {\cal P}_{\gamma}\cap L^{\infty}(X)$. So assume this and let $\omega>0$ be a hermitian metric over $X$. By the regularization result of \cite{Dem3} there exists a family of functions $(\psi_{\varepsilon})_{\varepsilon>0}$, $\psi_{\varepsilon}\in {\cal P}_{\gamma+\varepsilon\omega}\cap C^{\infty}(X)$ such that $\psi_{\varepsilon}\downarrow \psi$ as $\varepsilon\rightarrow 0^+$. Consider now the integrals
$$
I_j:=\int_X -\psi\,\gamma^j\wedge \gamma_{\varphi}^{n-j}\,,
$$
for all $j=0,...,n$. Then $I_j\leq I_{j+1}+\int_X\gamma^n$. In fact by Stokes' formula
\begin{eqnarray*}
I_j&=&I_{j+1}-
\lim_{\varepsilon\rightarrow 0^+}\,\int\limits_X \psi_{\varepsilon}\,\gamma^j\wedge i\partial\bar\partial \varphi\wedge\gamma_{\varphi}^{n-j-1} 
\\
\\
&=&
I_{j+1}-
\lim_{\varepsilon\rightarrow 0^+}\,\int\limits_X i\partial\bar\partial\psi_{\varepsilon} \wedge  \varphi\,\gamma^j\wedge\gamma_{\varphi}^{n-j-1}
\\
&\leq&
I_{j+1}+\int\limits_X  \varphi\,\gamma^{j+1}\wedge\gamma_{\varphi}^{n-j-1}
\leq
I_{j+1}+\int\limits_X\gamma^n\,.
\end{eqnarray*}
In this way we deduce the required inequality $I_0\leq I_n+n\int_X\gamma^n$.\hfill $\Box$
\\
\\
The following claim will be very useful for the rest of the paper.
\begin{claim}\label{StpdCv}
Let $(X,\omega)$ be a polarized compact connected K\"ahler manifold of complex dimension $n$ and let $\gamma$, $T$ be closed positive $(1,1)$-currents with continuous $($or more generally, bounded$\,)$ local potentials. Then for all $l=0,...,n$
$$
C_l:=\sup_{\psi\in {\cal P}^0_{\gamma}}\;\int\limits_X-\psi\,T^l\wedge \omega^{n-l}<+\infty
$$ 
and $\gamma_{\psi}\wedge T^l=T^l\wedge \gamma_{\psi}$ for all $\psi\in {\cal P}_{\gamma}$. 
\end{claim}
\emph{Proof}. The proof of the convergence of the constants $C_l$ goes by induction on $l=0,...,n$. The statement is true for $l=0$ by the first integral estimate of Lemma \ref{IntBndPSh}. So we assume it is true for $l$ and we prove it for $l+1$. Let $\psi_c:=\max\{\psi,c\}\in {\cal P}_{\gamma}$, $c\in \R_{<0}$. By the regularization result of \cite{Dem3} we find $(\psi_{c,\varepsilon})_{\varepsilon>0}$, $\psi_{c,\varepsilon}\in {\cal P}_{\gamma+\varepsilon\omega}\cap C^{\infty}(X)$ such that $\psi_{c,\varepsilon}\downarrow \psi_c$ as $\varepsilon\rightarrow 0^+$. Let us write $T=\theta+i\partial\bar\partial u$, with $\theta$ smooth, $\theta\leq K\omega$ and $u$ bounded with $\inf_X u=0$.
By using the monotone convergence theorem and Stokes' formula, we expand the integral 
\begin{eqnarray*}
&&\int\limits_X-\psi\,T^{l+1}\wedge\omega^{n-l-1}
= 
\lim_{c\rightarrow -\infty}\;\lim_{\varepsilon\rightarrow 0^+}\;\int\limits_X-\psi_{c,\varepsilon}\, T^{l+1}\wedge\omega^{n-l-1}
\\
\\
&=&
\lim_{c\rightarrow -\infty}\;\lim_{\varepsilon\rightarrow 0^+}\,\left[\;\int\limits_X-\psi_{c,\varepsilon}\, \theta\wedge T^{l}\wedge\omega^{n-l-1}-\int\limits_X \psi_{c,\varepsilon}\,i\partial\bar\partial u\wedge T^{l}\wedge\omega^{n-l-1}\right]
\\
\\
&\leq&
\lim_{c\rightarrow -\infty}\;\lim_{\varepsilon\rightarrow 0^+}\,\left[\;\int\limits_X-\psi_{c,\varepsilon}\, T^{l}\wedge K\omega^{n-l}-\int\limits_Xu\,i\partial\bar\partial\psi_{c,\varepsilon}\wedge T^{l}\wedge\omega^{n-l-1}\right]
\end{eqnarray*}
\begin{eqnarray*}
&=&
\int\limits_X-\psi\,T^{l}\wedge K\omega^{n-l}
+
\lim_{c\rightarrow -\infty}\;\lim_{\varepsilon\rightarrow 0^+}\,\left[-\int\limits_Xu\,(\gamma_{\psi_{c,\varepsilon}}+\varepsilon\omega)\wedge T^l\wedge\omega^{n-l-1}\right.
\\
\\
&+&
\left. \int\limits_Xu\,(\gamma+\varepsilon\omega)\wedge T^l\wedge\omega^{n-l-1}\right]
\\
\\
&\leq&
KC_l+\sup_X u\,\int\limits_X\gamma\wedge T^{l}\wedge\omega^{n-l-1}<+\infty\,,
\end{eqnarray*}
by the inductive hypothesis.
In sight of the symmetry of the exterior product we remark that the decreasing monotone convergence theorem implies
$$
\lim_{c\rightarrow -\infty}\,
\int\limits_X(\psi_{c}-\psi)T^l\wedge \omega^{n-l}=0\,,
$$
which means the convergence of the mass $\|(\psi_{c}-\psi)T^l\|_{\omega}(X)\rightarrow 0$ as $c\rightarrow -\infty$, in particular $\psi_{c}T^l\rightarrow \psi T$ weakly as $c\rightarrow -\infty$.  So by the weak continuity of the $i\partial\bar\partial$ operator we deduce 
\begin{eqnarray}
\gamma_{\psi_{c}}\wedge T^l \longrightarrow \gamma_{\psi}\wedge T^l\,, \label{bbCnv}
\end{eqnarray}
weakly as $c\rightarrow -\infty$. Moreover the weak continuity of the $i\partial\bar\partial$ operator implies by induction on $l$
$$
T^l\wedge \gamma_{\psi_{c}}  \longrightarrow T^l\wedge \gamma_{\psi}\,,
$$ 
weakly as $c\rightarrow -\infty$. This combined with \eqref{bbCnv} implies $\gamma_{\psi}\wedge T^l=T^l\wedge \gamma_{\psi}$
.\hfill$\Box$
\\
\\
In the particular case $T=\gamma$ big, the constant
$$
0<C(\gamma):=n+\sup_{\psi\in {\cal P}^0_{\gamma}}\;\;\{\gamma\}^{-n}\int\limits_X-\psi\,\gamma^n<+\infty
$$
satisfies the capacity estimate of Statement (A) in Lemma \ref{CapBndPSh}, by inequality (\ref{ParCap0}) and Claim \ref{CapIest}. 
Thus if $(\gamma_t)_{t>0}$ is a family satisfying the hypothesis $(C1)$ and $(C2a)$ of Statement (C) in Theorem  \ref{Kolo} and $K_t=(\{\gamma_t\}^n)^{1/n}$, then the constant $C(\gamma)$ satisfies the stability properties of Statement (A) of the Lemma \ref{CapBndPSh}, and we can use
the induction in the proof of Claim \ref{StpdCv} with
$T=\gamma_t$, $\theta=\theta_t$, $u=u_t$ and $K=K_t$ to get
$$
C_1\leq K_t \int\limits_X-\psi\,\omega^n+\sup_X u_t\,\int\limits_X \gamma_t\wedge \omega^{n-1}
\leq K_t\int\limits_X-\psi\,\omega^n+RK_t\,\int\limits_X \omega^n\,,
$$
where $R\geq \sup_X u_t$ and in general
$$
C_{l+1}\leq K_t C_l+R\int\limits_X \gamma^{l+1}_t\wedge \omega^{n-l-1}\leq K_t C_l+RK_t^{l+1}\,\int\limits_X \omega^n\,.
$$
We deduce 
$$
C_n\leq K_t^n \int_X-\psi\,\omega^n+nRK_t^n\int_X\omega^n\,.
$$
We now prove Statement (B) of Lemma \ref{CapBndPSh}.
In fact let $f:=\{\gamma\}^{-n}\gamma^n/\Omega\geq 0$. Then the uniform estimate for the integral 
$$
\{\gamma\}^{-n}\int\limits_X -\psi\,\gamma^n=\frac{1}{\alpha}\,\int\limits_X-\alpha\psi f\,\Omega
$$ 
follows from the elementary inequality $-\alpha\psi f\leq e^{-\alpha\psi}-1+f\log(1+f)$ combined with the uniform estimate $\int_Xe^{-\alpha\psi}\Omega\leq C$ of Lemma \ref{IntBndPSh}. In this case the required stability properties of the constant $C(\gamma, \Omega)>0$ in the capacity estimate are obvious.
\hfill $\Box$
\begin{lemma}{\bf (Comparison Principle).}\label{comp-Prnc}
Let $X$ be a compact complex manifold of complex dimension~$n$ and let $\gamma$ be a closed real $(1,1)$-current with bounded local potentials and consider $\varphi,\,\psi\in {\cal P}_{\gamma}\cap L^{\infty}(X)$. Then
$$
\int\limits_{\varphi<\psi}\gamma^n_{\psi}\;\leq \int\limits_{\varphi<\psi}\gamma^n_{\varphi}\;.
$$
\end{lemma}
\emph{Proof.} Let $\Theta:=\left(\gamma+i\partial\bar\partial \max\{\varphi,\psi\}\right)^n$.
By the inequality of measures 
$$
\Theta\;\ge\; \I_{_{\varphi\ge \psi}}\,\gamma_{\varphi}^n\;+\;\I_{_{\varphi < \psi}}\,\gamma_{\psi}^n\,,
$$
proved in \cite{Dem1}, we infer
\begin{eqnarray*}
\int\limits_{\varphi<\psi}\Theta
\;\ge
\int\limits_{\varphi<\psi}\gamma^n_{\psi}\;,
\qquad\qquad
\int\limits_{\varphi\ge\psi}\Theta
\;\ge
\int\limits_{\varphi\ge\psi}\gamma^n_{\varphi}\,.
\end{eqnarray*}
This combined with Stokes' formula implies
\begin{eqnarray*}
\int\limits_{\varphi<\psi}\gamma_{\psi}^n
\;\le\;
\int\limits_X\Theta
\;-
\int\limits_{\varphi\ge\psi}\Theta
\;\le
\int\limits_X\gamma_{\varphi}^n 
\;-\;
\int\limits_{\varphi\ge\psi}  \gamma_{\varphi}^n\;=\int\limits_{\varphi<\psi}\gamma_{\varphi}^n\,.
\end{eqnarray*}
\hfill $\Box$
\\
\\
We recall now the following lemma due to Ko{\l}odziej \cite{Kol1}, (see also \cite{Ti-Zhu1},
\cite{Ti-Zhu2}).
\begin{lemma}\label{KolIterLm}
Let $a:(-\infty,0]\rightarrow [0,1]$, be a monotone non-decreasing function such that for some $B>0$, $\delta>0$ the inequality
$$
t\,a(s)\leq B\,a(s+t)^{1+\delta}
$$
holds for all $s\leq 0,\, t\in [0,1],\,s+t\leq 0$. Then for all $S<0$ such that $a(S)>0$ and all $D\in [0,1],\,S+D\leq 0$ we have the estimate
$$
D\leq e(3+2/\delta)B\,a(S+D)^{\delta}\,.
$$
\end{lemma}
The following lemma is a simple application of the main result in Bedford-Taylor \cite{Be-Te} and of the monotone increasing convergence theorem in pluripotential theory.
\begin{lemma}\label{comprCap}
Let $X$ be a compact connected complex manifold of complex dimension $n$, let $\gamma$ be a big closed positive $(1,1)$-current with continuous local potentials and let $\Omega>0$ be a smooth volume form. Then there exist constants $\alpha=\alpha(\gamma,\Omega)>0$,  $C=C(\gamma,\Omega)>0$ such that for all Borel sets $E\subset X$ we have
\begin{eqnarray}\label{BTVol-Cap-Est}
\int\limits_E\Omega\leq e^{\alpha}C e^{-\alpha/\Cap_{\gamma}(E)^{1/n}}\,.
\end{eqnarray}
In particular $\Cap_{\gamma}(E)=0$ implies $\int_E\Omega=0$.
\end{lemma}
\emph{Proof}. It is sufficient  to prove this estimate for an arbitrary compact set. In fact assume \eqref{BTVol-Cap-Est} for compact sets and let $(K_j)_j$, $K_j\subset K_{j+1}\subset E$ be a family of compact sets such that $\int_{K_j}\Omega\rightarrow\int_E\Omega$ as $j\rightarrow +\infty$. Set $U:=\cup_jK_j\subset E$ and take the limit in \eqref{BTVol-Cap-Est} with $E$ replaced by $K_j$. By \eqref{CapProp} we deduce 
$$
\int\limits_E\Omega\leq e^{\alpha}C e^{-\alpha/\Cap_{\gamma}(U)^{1/n}}\leq e^{\alpha}C e^{-\alpha/\Cap_{\gamma}(E)^{1/n}}\,.
$$
We prove now \eqref{BTVol-Cap-Est} for compact sets $K\subset X$.
For this purpose, consider the function introduced in \cite{Sic}, \cite{G-Z}
$$
\Psi_K(x):=\sup\{\varphi(x)\,\mid\,\varphi\in {\cal P}_{\gamma}\,,\,\varphi_{\mid_{K}}\leq 0\}\,.
$$
We remark that $\Psi_K\ge 0$ over $X$ and $(\Psi_K)_{\mid_{K}}= 0$ since $0\in {\cal P}_{\gamma}$ by the positivity assumption on $\gamma$. Assume 
$$
\int_K\Omega\not =0\,,
$$
otherwise there is nothing to prove. In this case there exists a constant $C_K>0$ such that $\sup_X\,\varphi\leq C_K$ for all 
$\varphi\in {\cal P}_{\gamma}\,,\,\varphi_{\mid_{K}}\leq 0$.
In fact
let 
$$
S_K:=\{\varphi\in {\cal P}_{\gamma}\,\mid\,\varphi_{\mid_{K}}\leq 0\}
$$ 
and set $\tilde\varphi:=\varphi-\sup_X\varphi$. By contradiction we would get a sequence $\varphi_j\in S_K$ such that $\sup_X\varphi_j\rightarrow+\infty$. This implies 
$$
\sup_{K}\tilde\varphi_j\rightarrow-\infty\,
$$ and so 
$$
\int_{K}-\tilde\varphi_j\,\Omega\geq -\left(\,\int_K\Omega\,\right)\sup_{K}\tilde\varphi_j\rightarrow+\infty\,,
$$
which contradicts the first integral estimate of Lemma \ref{IntBndPSh}.
\\
Then it follows from quite standard local arguments that the upper regularization $\Psi^*_K\in {\cal P}_{\gamma}$. (Here we use the assumption that the local potentials of $\gamma$ are continuous.)
Moreover $\Psi^*_K\in L^{\infty}(X)$, $\Psi^*_K\geq 0$ and $\Psi^*_K= 0$ over the interior $K^0$ of $K$. 
\\
\\
We recall now the following well known consequence of a result of Bedford and Taylor \cite{Be-Te}. 
\begin{theorem}\label{DirichBedfTayl}
Let $\varphi\in {\cal P}_{\gamma}\cap L^{\infty}(X)$ and let $B$ be an open coordinate ball. Then there exists $\hat\varphi\in {\cal P}_{\gamma}\cap L^{\infty}(X)$, $\hat\varphi\geq \varphi$ such that $\gamma^n_{\hat\varphi}=0$ on $B$ and $\hat\varphi=\varphi$ on $X\smallsetminus B$. Moreover if $\varphi_1\leq \varphi_2$, then  $\hat\varphi_1\leq \hat\varphi_2$.
\end{theorem}
This implies the following quite standard fact in pluripotential theory \cite{Sic}, \cite{Dem1}, \cite{G-Z}.
\begin{corollary}\label{CorBedfTay}Let $K\subset X$ be a compact set such that $\int_K\Omega\not =0$. Then the extremal function $\Psi^*_K\in {\cal P}_{\gamma}\cap L^{\infty}(X)$ satisfies $\Psi^*_K\geq 0$ over $X$, $\Psi^*_K= 0$ over the interior $K^0$ of $K$ and 
$\gamma^n_{\Psi^*_K}=0$ over $X\smallsetminus K$. 
\end{corollary}
\emph{Proof}. By the classical Choquet lemma there exists a sequence  $(\varphi_j)_j\subset S_K$, $\varphi_j\geq 0$ such that $\Psi^*_K=(\sup_j\varphi_j)^*$. We can assume that this sequence is increasing. Otherwise, set $\tilde \varphi_1:=\varphi_1$ and 
$$
\tilde \varphi_j:=\max \{\varphi_j, \tilde \varphi_{j-1}\}\in S_K\,.
$$
Let $B$ be an open coordinate ball in $X\smallsetminus K$ and let $\hat\varphi_j\in S_K$ be a solution of the Dirichlet problem $\gamma^n_{\hat\varphi_j}=0$ over $B$ as in Theorem \ref{DirichBedfTayl}. Thus the sequence $(\hat\varphi_j)_j\subset S_K$ is still increasing and $\Psi^*_K=(\sup_j\hat\varphi_j)^*$. Remember also that the plurisubharmonicity implies that $\Psi^*_K=\lim_j\hat\varphi_j$ almost everywhere. By the monotone increasing theorem from classical pluripotential theory, we infer $\gamma^n_{\Psi^*_K}=0$ on $B$, and the conclusion follows from the fact that $B$ is arbitrary.\hfill~~$\Box$
\\
\\
By using a basic fact about Lebesgue measure theory and the second integral estimate of Lemma \ref{IntBndPSh} we get
$$
\int\limits_{K}\Omega=\int\limits_{K^0}\Omega=\int\limits_{K^0}e^{-\alpha\,\Psi^*_K}\,\Omega\,\leq\,
\int\limits_Xe^{-\alpha\,\Psi^*_K}\,\Omega\,\leq \,Ce^{-\alpha\sup_X\Psi^*_K }\,.
$$
Set $A_K:=\sup_X\Psi^*_K$. If $A_K > 1$ set $\varphi:=A_K^{-1}\Psi^*_K$. Then $0\leq \gamma_{\Psi^*_K}\leq A_K\gamma_{\varphi}$ and so $\varphi\in {\cal P}_{\gamma}[0,1]$. By corollary \ref{CorBedfTay} we deduce 

\begin{eqnarray*}
\{\gamma\}^nA_K^{-n}=A_K^{-n}\int\limits_{K}\gamma^n_{\Psi^*_K}\,\leq\,\int\limits_{K}\gamma^n_{\varphi}\,\leq \,\{\gamma\}^n\Cap_{\gamma}(\,K\,)\,,
\end{eqnarray*}
thus $-\alpha A_K\leq -\alpha/\Cap_{\gamma}(\,K\,)^{1/n}$ by the bigness assumption on the current $\gamma$. If $A_K\leq 1$ then $\Psi^*_K\in {\cal P}_{\gamma}[0,1]$ and so 
\begin{eqnarray*}
1=\{\gamma\}^{-n}\int\limits_{K}\gamma^n_{\Psi^*_K}\,\leq\,\Cap_{\gamma}(\,K\,)\leq\Cap_{\gamma}(X)= 1\,.
\end{eqnarray*}
In both cases we reach the required conclusion.
\hfill $\Box$
\\
\\
{\bf Proof of Theorem \ref{Kolo}, part A}. 
\\
We can assume $\sup_X\psi=0$. 
Let $U_s:=\{\psi< s\}$, $s\leq 0$, $t\in [0,1]$, $s+t\leq 0$, 
$\varphi\in {\cal P}_{\gamma}[-1,0]$ and set 
$$
V:=\{\psi-s-t<t\varphi\}\,.
$$
Then we have inclusions 
$U_{s}\subset V\subset U_{s+t}$. By using the Comparison  Principle \eqref{comp-Prnc} we infer
$$
t^n\int\limits_{U_s}\gamma_{\varphi}^n\;\leq \int\limits_{U_s}\gamma_{t\varphi}^n\;\leq \int\limits_{V}\gamma_{t\varphi}^n\;
\leq \int\limits_{V}\gamma_{\psi}^n\;\leq \int\limits_{U_{s+t}}\gamma_{\psi}^n\,,
$$
thus combining this with H\"{o}lder inequality in Orlicz spaces \eqref{Hold-Ineq}, formula \eqref{Exp1-Vol} and Lemma \ref{comprCap} we obtain
\begin{eqnarray*}
t^n\Cap_{\gamma}(U_s)
&\leq & 
\{\gamma\}^{-n}\int\limits_{U_{s+t}}\gamma_{\psi}^n
\;=\,\{\gamma\}^{-n}\int\limits_{U_{s+t}}f\,\Omega
\\
\\
&\leq& \{\gamma\}^{-n}C_{\varepsilon_0}\|f\|_{L\log^{n+\varepsilon}L(X)}\cdot\|1\|_{\Exp^{\frac{1}{n+\varepsilon}}L(U_{s+t})}
\\
\\
&=&
\frac{\{\gamma\}^{-n}C_{\varepsilon_0}\|f\|_{L\log^{n+\varepsilon}L(X)}}{\log^{n+\varepsilon}\left(1+1/\Vol_{\Omega}(U_{s+t})\right)}
\\
\\
&\leq&
\frac{\{\gamma\}^{-n}C_{\varepsilon_0}\|f\|_{L\log^{n+\varepsilon}L(X)}}{\log^{n+\varepsilon}
\left(1+e^{-\alpha}C^{-1}e^{\alpha/\Cap_{\gamma}(U_{s+t})^{1/n}}\right)}
\\
\\
&\leq&
C_{\varepsilon_0}(k/\alpha)^{n+\varepsilon}\{\gamma\}^{-n}\|f\|_{L\log^{n+\varepsilon}L(X)}\Cap_{\gamma}(U_{s+t})^{(n+\varepsilon)/n}\,.
\end{eqnarray*}
(Here the constant $C>0$ depends on the same quantities as the constant $C_1$ in Statement~A and
$k>0$ is a constant such that 
$$
k^{-1}\alpha/x\leq \log(1+e^{-\alpha}C^{-1}e^{\alpha/x})\,,
$$ 
for all $x\in (0,1]$). So if we set $\delta:=\varepsilon/n$ and 
$$
B:=C_{\varepsilon_0}^{1/n}(k/\alpha)^{1+\varepsilon/n}I_{\gamma,\varepsilon}(f)^{1/n}\,,
$$
we deduce that the function $a(s):=\Cap_{\gamma}(U_s)^{1/n}$, $s\leq 0$, satisfies the hypothesis of Lemma \ref{KolIterLm}. (We use here the inequality \eqref{LlogL-Int}.)
Consider now the function $\kappa(t):=K_{\delta}B\,t^{\delta}$, with constant  $K_{\delta}:=e(3+2/\delta)$. Remember also the uniform capacity estimate 
$$
a(s)\leq C\,(-s)^{-1/n}\,,
$$ 
of Lemma \ref{CapBndPSh}. Let now $\eta>1$ be arbitrary. We claim that $a(S_{\eta})=0$ for 
$$
-S_{\eta}=C^n(K_{\delta}B\,\eta)^{n/\delta}+1\,.
$$
The fact that the function $a$ is left continuous (by formula \eqref{CapProp}) will imply that $a(S_1)=0$ also. Remark that $S_{\eta}$ is a solution of the equation 
$$
C(-S_{\eta}-1)^{-1/n}=\kappa^{-1}(\eta^{-1})\,,
$$
where $\kappa^{-1}$ is the inverse of the function $\kappa$. So if we assume by contradiction that $a(S_{\eta})>0$ we deduce by Lemmas \ref{KolIterLm} and \ref{CapBndPSh} 
$$
1\leq \kappa(a(S_{\eta}+1))\leq \kappa(C(-S_{\eta}-1)^{-1/n})=\eta^{-1}<1\,,
$$
which is a contradiction. Thus if we set $-I:=\max\{s\leq 0\,\mid \,a(s)=0\}$ we obtain 
$$
I\leq -S_1\leq C^n(K_{\delta}B)^{n/\delta}+1\,,
$$
which by arranging the coefficients yields the right hand side of the estimate in Statement A of Theorem \ref{Kolo}. Moreover by definition $\Cap_{\gamma}(U_{-I})=0$, thus $\Vol_{\Omega}(U_{-I})=0$ by Lemma \ref{comprCap}. The fact that the current $\gamma$ has continuous local potentials implies that the function $\psi$ is upper semicontinuous, so the set $U_{-I}$ is open, thus empty. This implies the required conclusion.\hfill $\Box$
\\
\\
{\bf Proof of part B}. 
\\
Set
$a:=\max\{\|\varphi\|_{L^{\infty}(X)},\,\|\psi\|_{L^{\infty}(X)}\}$, consider $\theta\in {\cal P}_{\gamma}[0,1]$, $s\geq 0$, $t\in [0,1]$
and set 
$$
V:=\left\{\varphi<\frac{t}{1+a}\,\theta+\left(1-\frac{t}{1+a}\right)\psi-s-t\right\}\,.
$$
Then the obvious inequality $0\leq -\frac{t}{1+a}\psi\leq \frac{at}{1+a}$ implies the inclusions 
\\
$\{\varphi-\psi< -s-t\}\subset V\subset \{\varphi-\psi< -s\}$. Thus by applying the Comparison Principle \eqref{comp-Prnc} as in \cite{Kol2} we obtain
\begin{eqnarray*}
\frac{t^n}{(1+a)^n}\int\limits_{\varphi-\psi< -s-t}\gamma_{\theta}^n
&\leq& \int\limits_V\left[\frac{t}{1+a}\,\gamma_{\theta}+\left(1-\frac{t}{1+a}\right)\gamma_{\psi}\right]^n
\\
\\
&\leq& 
\int\limits_V\gamma_{\varphi}^n\;\;\leq\;\; \int\limits_{\varphi-\psi< -s}\gamma_{\varphi}^n \,.
\end{eqnarray*}
By inverting the roles of $\varphi$ and $\psi$ in the previous inequality and by summing up we get
$$
\frac{t^n}{(1+a)^n}\,\int\limits_{ |\varphi-\psi |>s+t}\gamma_{\theta}^n
\;\;
\leq
\;\;
\int\limits_{| \varphi-\psi |>s}(f+g)\,\Omega\,.
$$
By taking the supremum over $\theta$ we obtain the capacity estimate
\begin{eqnarray}\label{CapStab}
t^n\Cap_{\gamma} (| \varphi-\psi |>s+t)\leq (1+a)^n \{\gamma\}^{-n}\int\limits_{| \varphi-\psi |>s}(f+g)\,\Omega\,,
\end{eqnarray}
for all $s\geq 0$, $t\in [0,1]$. 
Set $U_s:=\{| \varphi-\psi |>s\}\subset X$. By combining Lemma \ref{comprCap} with a computation similar to the one in the proof of part A we obtain
\begin{eqnarray*}
t^n\Cap_{\gamma} (U_{s+t})
&\leq& 
(1+a)^n \{\gamma\}^{-n}C'_{\varepsilon_0}\|f+g\|_{L\log^{n+\varepsilon_0}L(X)}\Cap_{\gamma}(U_s)^{(n+\varepsilon_0)/n}
\\
\\
&\leq& B^n\Cap_{\gamma}(U_s)^{(n+\varepsilon_0)/n}\,,
\end{eqnarray*}
where the constant $B>0$ depends on the same quantities as the constant $C_2$ in Statement (B) of Theorem \ref{Kolo}. We deduce that the function $a(s):=\Cap_{\gamma}(U_{-s})^{1/n}$, $s\leq 0$, satisfies the hypothesis of Lemma \ref{KolIterLm} with $\delta=\varepsilon_0/n$. On the other hand, the capacity estimate \eqref{CapStab} combined with H\"{o}lder's inequality in Orlicz spaces implies for all $t\in [0,1]$ the inequalities
\begin{eqnarray}
t^n\Cap_{\gamma} ( |\varphi-\psi| >2t)
&\leq&
(1+a)^n \{\gamma\}^{-n}\int\limits_{| \varphi-\psi |>t}(f+g)\,\Omega\nonumber
\\\nonumber
\\
&\leq&
\frac{(1+a)^n \{\gamma\}^{-n}}{t} \int\limits_X| \varphi-\psi| (f+g)\,\Omega\nonumber
\\\nonumber
\\
&\leq&
\frac{2(1+a)^n \{\gamma\}^{-n}}{t} \|\varphi-\psi \|_{\Exp L(X)}\|f+g\|_{L\log L(X)}\nonumber
\\\nonumber
\\
&\leq&
\frac{4K(1+a)^n}{t} \|\varphi-\psi \|_{\Exp L(X)}\,.\label{CapExpNor1}
\end{eqnarray}
\begin{claim}\label{Exp-L1}
If $\|\varphi-\psi \|_{L^1(X)}\leq 1/2$, then there exists a constant $C_a>0$ such that
$$
\|\varphi-\psi \|_{\Exp L(X)}\leq C_a/\log\|\varphi-\psi \|^{-1}_{ L^1(X)}\,.
$$
\end{claim}
\emph{Proof}. We assume $\|\varphi-\psi \|_{L^1(X)}>0$, otherwise there is nothing to prove. Set 
$$
C_{k,a}:=k(e^{2a/k}-1)/(2a)\,,
$$
$k>0$. Then for all $k>0$ and all $x\in [0,2a/k]$ the inequality $e^x-1\leq C_{k,a}\,x$ holds. Thus the inequality $| \varphi-\psi |/k\leq 2a/k$ implies 
$$
\int\limits_X \left(e^{|\varphi-\psi|/k}-1\right)\Omega\leq C_{k,a}\int\limits_X\frac{|\varphi-\psi|}{k}\,\Omega\,.
$$
We get from there the implication
\begin{eqnarray}\label{ImpExp-L1}
\|\varphi-\psi \|_{L^1(X)}=k/C_{k,a}\quad\Longrightarrow \quad \|\varphi-\psi \|_{\Exp L(X)}\leq k\,,
\end{eqnarray}
since by definition
$$
\|\varphi-\psi \|_{\Exp L(X)}:=\inf\left\{k>0\,\mid\;\int\limits_X \left(e^{|\varphi-\psi|/k}-1\right)\Omega\leq 1\,\right\}\,.
$$
So if we set $\mu(k):=k/C_{k,a}>0$ we deduce by the implication \eqref{ImpExp-L1}
\begin{eqnarray}\label{EstExp-L1}
\|\varphi-\psi \|_{\Exp L(X)}\leq \mu^{-1}\left(\|\varphi-\psi \|_{ L^1(X)}\right)\,,
\end{eqnarray}
where $\mu^{-1}:\R_{>0}\rightarrow\R_{>0}$ is the inverse function of $\mu$.
Explicitly 
$$
\mu^{-1}(y)=2a/\log(1+2a/y)\,,
$$
for all $y>0$. Now there exists a constant $C_a>0$ such that 
$$
\mu^{-1}(y)\leq C_a/\log(1/y)\,,
$$ 
for all $y\in (0,1/2]$. This combined with \eqref{EstExp-L1} implies the conclusion.\hfill$\Box$
\\
\\
Combining Claim \ref{Exp-L1} with the estimate \eqref{CapExpNor1} we infer the capacity estimate
\begin{eqnarray}\label{CapExpNor}
a(-t)\leq  \frac{C}{t^{1+1/n}}\left(\log\|\varphi-\psi \|^{-1}_{ L^1(X)}\right)^{-1/n}\,,
\end{eqnarray}
where the constant $C>0$ depends on the same quantities as the constant $C_2$ in Statement~B. Set now $C_2:=C^n(2K_{\delta}B)^{n/\delta}>0$ (with $K_{\delta}>0$ as in the proof of Statement~(A)) and define
$$
t:=C_2^{\alpha_0}\left(\log\|\varphi-\psi \|^{-1}_{ L^1(X)}\right)^{-\alpha_0}\,.
$$
The hypothesis $t\in (0,1]$ combined with the hypothesis of Claim \ref{Exp-L1} forces the condition $\|\varphi-\psi \|_{ L^1(X)}\leq \min\{1/2, e^{-C_2}\}$.
Moreover $t$ is solution of the equation
$$
\frac{C}{t^{1+1/n}}\left(\log\|\varphi-\psi \|^{-1}_{ L^1(X)}\right)^{-1/n}
=\kappa^{-1}\left(\frac{t}{2}\right)\,,
$$
where $\kappa^{-1}$ is the inverse  of the function $\kappa$ introduced in the proof of part~A. We claim that $a(-2t)=0$. Otherwise, by Lemma \ref{KolIterLm} and inequality \eqref{CapExpNor}, we infer
$$
0<t\leq \kappa(a(-t))\leq \kappa(\kappa^{-1}(t/2))=t/2\,,
$$
which is absurd. We deduce 
$$
\Vol_{\Omega}(| \varphi-\psi |>2t)=0
$$ 
by Lemma \ref{comprCap}. We prove now that the set 
$$
U_{2t}=\{| \varphi-\psi |>2t\}\subset X\,,
$$ 
is empty, which will imply the desired $L^{\infty}$-stability estimate. The fact that $|\varphi-\psi|\leq 2t$ a.e. over $X$, implies 
$$
\left|\;\dashint_{B(x,r)}(\varphi-\psi)\,d\lambda\;\right|\leq 2t\,,
$$ 
for all coordinate open balls $B(x,r)\subset X$. (The symbol $\dashint_B$ represents the mean value operator.) By elementary properties of plurisubharmonic functions follows
$$
\varphi(x)-\psi(x)=\lim_{r\rightarrow 0^+}\,\dashint_{B(x,r)}(\varphi-\psi)\,d\lambda\,,
$$
for all $x\in X$. We infer $| \varphi-\psi |\leq 2t$ over $X$.
\hfill $\Box$
%%%%%%%%%%%%%%%%%%%%%%%%%%%%%%%%%%%%%%%%%%%%%%%%%%%%%%%%%%%%%%%%%%%%%%%%%%%%%%%%%%%%%%%%%%%%%%%%%%%%%%%%%%%%%%%%%%%%%%%%%%%%%%%%%%%%%%%%%%%%%%%%%%%%%%%%%%%%%%%%%%%%%%%%%%%%%%%%%%%%%%%%%%%%%%%%%%%%%%%%%%%%%%%%%%%%%%%%%%%
%%%%%%%%%%%%%%%%%%%%%%%%%      C0 ESTIMATE %%%%%%%%%%%%%%%%%%%%%%%%%%%%%%%%%%%%%%%%%%%%%%%%%%%%%%%%%%%%%%%%%%%%%%%%%%%%%%%%%%%%%%%%%%%%%%%%%%%%%%%%%%%%%%%%%%%%%%%%%%%%%%%%%%%%%%%%%%%%%%%%%%%%%%%%%%%%%%%%%%%%%%%%%%%%%%%%%%%%%%%%%%%%%%%%%%%%%%%%%%%%%%%%%%%%%%%%%%%%%%%%%%%%%%%%%%%%%%%%%%%%%%%%%%%%%%%%%%%%%%%%%%%%%%%%%%%%%%%%%%%%%%%%%%%%%%%%%%%%%%%%%
\begin{corollary}\label{C0est}
Let $(X,\omega)$ be a polarized compact connected $n$-dimensional K\"ahler manifold, $\Omega>0$ a smooth volume form and  $\gamma\geq 0$ a big closed smooth $(1,1)$-form. Take also $f\in L\log^{n+\delta}L(X)$, $\delta>0$, such that $\int_X\gamma^n=\int_Xf\,\Omega$ and let $(f_{\varepsilon})_{\varepsilon>0}\subset C^{\infty}(X)$ be a family converging to $f$ in the $L\log^{n+\delta}L(X)$-norm as $\varepsilon\rightarrow 0^+$, satisfying the integral condition
\begin{eqnarray}\label{IntgC0est}
\int\limits_X(\gamma+\varepsilon\omega)^n=\int\limits_Xf_{\varepsilon}\,\Omega\,.
\end{eqnarray}
Then, for any real number $\lambda\geq 0$, the unique solution of the non-degenerate complex Monge-Amp\`ere equation
\begin{eqnarray}\label{MAC0est}
(\gamma+\varepsilon\omega+i\partial\bar\partial \psi_{\varepsilon})^n=f_{\varepsilon}\,e^{\lambda\,\psi_{\varepsilon}}\Omega\,,
\end{eqnarray}
given by the Aubin-Yau solution of the Calabi conjecture
$($which in the case $\lambda=0$ is normalized by $\max_X\psi_{\varepsilon}=0)$ satisfies the uniform $L^{\infty}$-estimate
$
\|\psi_{\varepsilon}\|_{L^{\infty}(X)}\leq C(\delta,\gamma,\Omega)\,I_{\gamma,\delta}(f)^{\frac{n}{\delta}}+\;1
$.
\end{corollary}
\emph{Proof}. The existence of a regularizing family $f_{\varepsilon}$ of $f$ in $L\log^{n+\delta}L(X)$ follows from \cite{Ra-Re} page 364 or \cite{Iw-Ma}, Theorem 4.12.2, page 79.
We can always assume the integral condition \eqref{IntgC0est} 
otherwise we multiply $f_{\varepsilon}$ by a constant $c_{\varepsilon}>0$ which converges to $1$ by the normalizing condition $\int_X\gamma^n=\int_Xf\,\Omega$. We distinguish two cases.
\\
\\
{\bf Case $\lambda=0$}.
The hypothesis (C1) and (C2a) of Statement (C) of Theorem \ref{Kolo} are obviously satisfied for the family $(\gamma+\varepsilon\omega)_{\varepsilon}$.
We deduce that the constant $C_1=C_1(\delta\,,\,\gamma+\varepsilon\omega\,,\,\Omega)>0$ in the Statement of Theorem \ref{Kolo}, A does not blow up as $\varepsilon\rightarrow 0^+$. Moreover the uniform estimate
\begin{eqnarray}\label{LlgLf}
\|f_{\varepsilon}\|_{L\log^{n+\delta}L(X)}\leq C'\|f\|_{L\log^{n+\delta}L(X)}=:K
\end{eqnarray}
holds for all $\varepsilon\in (0,1)$. 
Thus by Theorem \ref{Kolo},~A we obtain the required uniform estimate $\|\psi_{\varepsilon}\|_{L^{\infty}(X)}\leq C:=C(\delta,\gamma,\Omega)\,I_{\gamma,\delta}(f)^{\frac{n}{\delta}}+\;1$. 
\\
\\
%%%%%%%%%%%%%%%%%%%%%%%%%%%%%%%%%%%%%%%%%%%%%%%%%%%%%%%%%%%%%%%%%%%%%%%%%%%%%%%%%%%%%%%%%%%%%%%%%%%%%%%%%%%%%%%%%%%%%%%%%%%%%%%%%%%%%%%%%%%%%%%%%%%%%%%%%%%%%%%%%%%%%%%%%%%%%%%%%%%%%%%%%%%%%%%%%%%%%%%%%%%%%%%%%%%%%%%%%%%%%%%%%%%%%%%%%%%%%%%%%%%%%%%%%%%%%%%%%%%%%%%%%%%%%%%
{\bf Case $\lambda>0$}. We start by proving the following lemma, which is a parti\-cular case of a more general result due to Yau (see \cite{Yau}, sect.\ 6, page~376).
\begin{lemma}\label{YauIth}
Let $(X,\omega)$ be a polarized compact K\"ahler manifold of complex dimension $n$, let $h$ be a smooth function such that $\int_X\omega^n=\int_X e^h\omega^n$  and let $\varphi\in {\cal P}_{\omega}$ be the unique solution of the complex Monge-Amp\`ere equation 
\begin{eqnarray}\label{coolCMANeg}
(\omega+i\partial\bar\partial \varphi)^n=e^{h+\lambda\varphi}\omega^n\,,
\end{eqnarray}
$\lambda>0$. Consider also two solutions $\varphi'$, $\varphi''\in {\cal P}_{\omega}$ of the complex Monge-Amp\`ere equation $(\omega+i\partial\bar\partial \psi)^n=e^h\omega^n$ such that $\min_X\varphi'=0=\max_X \varphi''$. Then $\varphi''\leq \varphi\leq \varphi'$.
\end{lemma}
\emph{Proof}. The argument is a simplification, in our particular case, of Yau's original argument for the proof of Theorem\ 4, sect.\ 6 in \cite{Yau}. Set $\varphi'_0:=\varphi',\,\varphi''_0:=\varphi'' $ and consider the solutions $\varphi'_j,\,\varphi''_j$ of the complex Monge-Amp\`ere equations given by the iteration
\begin{eqnarray}
(\omega+i\partial\bar\partial\varphi'_j)^n&=&e^{h+(\lambda+1)\varphi'_j-\varphi'_{j-1}}\,\omega^n\,,\label{MAith1}
\\\nonumber
\\
(\omega+i\partial\bar\partial\varphi''_j)^n&=&e^{h+(\lambda+1)\varphi''_j-\varphi''_{j-1}}\,\omega^n\,.\label{MAith2}
\end{eqnarray}
Notice that we can solve these equations even if the terms $e^{h-\varphi'_{j-1}}$, $e^{h-\varphi''_{j-1}}$ are not normalized, see Lemma 2 page 378 in \cite{Yau}. Set $L:=\lambda+1$ and consider
$$
(\omega+i\partial\bar\partial\varphi'_1)^n=e^{h+L(\varphi'_1-\varphi'_0)+\lambda \varphi'_0}\omega^n\geq e^{L(\varphi'_1-\varphi'_0)}e^h\omega^n=e^{L(\varphi'_1-\varphi'_0)}(\omega+i\partial\bar\partial\varphi'_0)^n\,.
$$
At a maximum point of $\varphi'_1-\varphi'_0$ we have the inequality 
$$
(\omega+i\partial\bar\partial\varphi'_0)^n\geq (\omega+i\partial\bar\partial\varphi'_1)^n\,.
$$
By plugging this into the previous one, we deduce $\varphi'_1\leq \varphi'_0$. We now prove by induction the inequality $\varphi'_{j}\leq \varphi'_{j-1}$. In fact by dividing $(\ref{MAith1})_j$ with $(\ref{MAith1})_{j-1}$ we get
$$
\frac{(\omega+i\partial\bar\partial\varphi'_j)^n}{(\omega+i\partial\bar\partial\varphi'_{j-1})^n}=e^{L(\varphi'_j-\varphi'_{j-1})-(\varphi'_{j-1}-\varphi'_{j-2})}\geq e^{L(\varphi'_j-\varphi'_{j-1})}\,.
$$
At a maximum point of $\varphi'_j-\varphi'_{j-1}$ we find again the inequality
$$
(\omega+i\partial\bar\partial\varphi'_j)^n\leq (\omega+i\partial\bar\partial\varphi'_{j-1})^n\,.
$$
Combining this with the previous one we deduce $\varphi'_{j}\leq \varphi'_{j-1}$. By applying a quite similar argument to \eqref{MAith2} we obtain also $\varphi''_{j-1}\leq \varphi''_j$.
We also prove by induction the inequality $\varphi''_{j}\leq \varphi'_{j}$, which is true by definition in the case $j=0$. By dividing $(\ref{MAith1})_j$ with $(\ref{MAith2})_j$ we get
$$
\frac{(\omega+i\partial\bar\partial\varphi'_j)^n}{(\omega+i\partial\bar\partial\varphi''_j)^n}=e^{L(\varphi'_j-\varphi''_j)-(\varphi'_{j-1}-\varphi''_{j-1})}\leq e^{L(\varphi'_j-\varphi''_j)}\,,
$$
by the induction hypothesis $\varphi''_{j-1}\leq \varphi'_{j-1}$. At a minimum point of $\varphi'_j-\varphi''_j$ we get
$$
(\omega+i\partial\bar\partial\varphi'_j)^n\geq (\omega+i\partial\bar\partial\varphi''_j)^n\,,
$$
hence $\varphi''_{j}\leq \varphi'_{j}$. As a conclusion, we have proved the sequence of inequalities
\begin{eqnarray}\label{IneqYauIther}
\varphi''_{0}\leq \varphi''_{j-1}\leq \varphi''_{j}\leq \varphi'_{j}\leq \varphi'_{j-1}\leq \varphi'_{0}\,.
\end{eqnarray}
We now prove a uniform estimate for the Laplacian of the potentials $\varphi'_{j}$.
The inequalities \ref{IneqYauIther} imply $0<2n+\Delta_{\omega}\varphi'_j\leq C\,B_j$, where $B_j>0$ satisfies the uniform estimate
\begin{eqnarray}\label{SecndOrdB}
0\geq C_1\,B_j^{\frac{1}{n-1}}-\left(2n+\max_X \Delta_{\omega}\varphi'_{j-1}\right)B_j^{-1}-C_0\,,
\end{eqnarray}
$C_0, \,C_1>0$, which is obtained by applying the maximum principle in a similar way as in Yau's proof of the second order estimate for the solution of the Calabi conjecture \cite{Yau}. (It can also be obtained by setting $\delta=l=h=0$ and $\tilde\omega_{\varepsilon}=\omega$ in step (B) in the proof of Theorem \ref{HgRegMA}, (see Appendix B). In the case $n=1$ the uniform estimate $0<2n+\Delta_{\omega}\varphi'_j\leq C'$ follows immediately from the inequalities \eqref{IneqYauIther}.) Fix now a constant $C_2>0$ such that the inequality
$$
C_1\,x^{1+\frac{1}{n-1}}\geq (C_0+2C)x-C_2\,,
$$
holds for all $x\geq 0$.
This implies by \eqref{SecndOrdB} the estimate
$$
2(2n+\Delta_{\omega}\varphi'_j)\leq 2C\,B_j\leq \left(2n+\max_X \Delta_{\omega}\varphi'_{j-1}\right)+C_2\,,
$$
thus
$$
2n+\max_X \Delta_{\omega}\varphi'_j\leq 2^{-j} \left(2n+\max_X \Delta_{\omega}\varphi'_0\right)+C_2\,,
$$
by iteration. By taking the derivative in the Green Formula (see \cite{Aub}, Th. 4.13 page 108) we get the identity
$$
d_x\varphi'_j=-\int\limits_Xd_xG_{\omega}(x,\cdot)\,\Delta_{\omega}\varphi'_j\,\omega^n\,,
$$
which implies the estimate 
$$
|\nabla_{\omega}\varphi'_j|_{\omega}\leq C_{\omega}\max_X \Delta_{\omega}\varphi'_j\leq K\,.
$$
By applying the complex version of the Evans-Krylov theory \cite{Ti2} we deduce the uniform estimate $\|\varphi'_j\|_{C^{2,\alpha}(X)}\leq K'$. This combined with \eqref{MAith1} implies that the monotone sequence $(\varphi'_j)_j$ converges in the $C^{2,\alpha}$-topology to the unique solution $\varphi$ of the complex Monge-Amp\`ere equation \eqref{coolCMANeg}. Then the conclusion follows from the inequalities \eqref{IneqYauIther}.
\hfill$\Box$
\\
\\
Consider now the solutions $\psi'_{\varepsilon}$, $\psi''_{\varepsilon}$, $\min_X\psi'_{\varepsilon}=0=\max_X\psi''_{\varepsilon}$ of the complex Monge-Amp\`ere equation \eqref{MAC0est} for $\lambda=0$.
By applying Lemma \ref{YauIth} we deduce $\psi''_{\varepsilon}\leq \psi_{\varepsilon}\leq \psi'_{\varepsilon}$ for all $\varepsilon>0$. By the argument in the case $\lambda=0$, we infer $\|\psi'_{\varepsilon}\|_{L^{\infty}(X)},\,\|\psi''_{\varepsilon}\|_{L^{\infty}(X)}\leq C$, thus  $\|\psi_{\varepsilon}\|_{L^{\infty}(X)}\leq C$. 
\hfill $\Box$
%%%%%%%%%%%%%%%%%%%%%%%%%%%%%%%%%%%%%%%%%%%%%%%%%%%%%%%%%%%%%%%%%%%%%%%%%%%%%%%%%%%%%%%%%%%%%%%%%%%%%%%%%%%%%%%%%%%%%%%%%%%%%%%%%%%%%%%%%%%%%%%%%%%%%%%%%%%%%%%%%%%%%%%%%%%%%%%%%%%%%%%%%%%%%%%%%%%%%%%%%%
%%%%%%%%%%%%%%%%%%%%%%%%%%%%%%%%%%%%%%%%%%%%%%%%%%%%%%%%%%%%%%%%%%%%%%%%%%%%%%%%%%%%%%%%%%%%%%%%%%%%%%%%%%%%%%%%%%%%%%%%%%%%%%%%%%%%%%%%%%%%%%%%%%%%%%%%%%%%%%%%%%%%%%%%%%%%%%%%%%%%%%%%%%%%%%%%%%%%%%%%%%
\section{Currents with Bedford-Taylor type singularities}
%%%%%%%%%%%%%%%%%%%%%%%%%%%%%%%%%%%%%%%%%%%%%%%%%%%%%%%%%%%%%%%%%%%%%%%%%%%%%%%%%%%%%%%%%%%%%%%%%%%%%%%%%%%%%%%%%%%%%%%%%%%%%%%%%%%%%%%%%%%%%%%%%%%%%%%%%%%%%%%%%%%%%%%%%%%%%%%%%%%%%%%%%%%%%%%%%%%%%%%%%%
%%%%%%%%%%%%%%%%%%%%%%%%%%%%%%%%%%%%%%%%%%%%%%%%%%%%%%%%%%%%%%%%%%%%%%%%%%%%%%%%%%%%%%%%%%%%%%%%%%%%%%%%%%%%%%%%%%%%%%%%%%%%%%%%%%%%%%%%%%%%%%%%%%%%%%%%%%%%%%%%%%%%%%%%%%%%%%%%%%%%%%%%%%%%%%%%%%%%%%%%%%
In the situation we have to consider, the relevant class of currents which can be used as the input of Monge-Amp\`ere operators is defined as follows.

\begin{definition} On a complex manifold, we consider the class
$\BT$ of closed positive $(1,1)$-currents $\Theta$ whose exterior products $\Theta^k$, $0\leq k\leq n$, can be defined 
inductively in the sense
of Bedford-Taylor, namely, if $\Theta=i\partial\bar\partial\psi$ on 
any open set, then $\psi\Theta^k$ is locally of finite mass and
$\Theta^{k+1}=i\partial\bar\partial(\psi\Theta^k)$ for~$k<n$.
\end{definition}

Notice that the local finiteness of the mass of $\psi\Theta^k$ is independent
of the choice of the psh potential~$\psi$, and that this assumption allows indeed to compute inductively $i\partial\bar\partial(\psi\Theta^k)$ in the sense of currents. Now, if $\chi$ is a $(1,1)$-cohomology class, we set 
\begin{eqnarray}\label{BTchi}
\BT_\chi=\BT\cap\;\chi.
\end{eqnarray}
Let $\gamma\geq 0$ be a closed positive $(1,1)$-current with continuous local potentials. We define corresponding classes of potentials
\begin{eqnarray*}
{\cal P}\BT_{\gamma}&:=&\left\{\varphi\in {\cal P}_{\gamma}\,|\,\gamma+i\partial\bar\partial \varphi\in \BT_{\{\gamma\}}\right\}\,,\\
\noalign{\vskip5pt}
{\cal P}\BT^0_{\gamma}&:=&\{\varphi\in {\cal P}\BT_{\gamma}\,\mid\,\sup_X\varphi=0 \}\,.
\end{eqnarray*}
Let $\varphi\in {\cal P}\BT_{\gamma}$ with zero Lelong numbers. It is well known from the work of the first author \cite{Dem4} (which becomes drastically simple in this particular case), that there exists a family $(\varphi_{\varepsilon})_{\varepsilon>0}$, $\varphi_{\varepsilon}\in {\cal P}_{\gamma+\varepsilon\omega}\cap C^{\infty}(X)$, such that $\varphi_{\varepsilon}\downarrow \varphi$ as $\varepsilon\downarrow 0^+$. 
%In the case $\varphi\in {\cal P}_{\gamma}\cap L^{\infty}(X)$ the convergence of $\varphi_{\varepsilon}$ is also uniform. 
In the case the Lelong numbers of $\varphi$ are not zero we can chose $R>0$ sufficiently big such that $0\leq \gamma+R\omega+i\partial\bar\partial \varphi_{\varepsilon}$ for all $\varepsilon\in (0,1)$ and $\varphi_{\varepsilon}\downarrow \varphi$ as $\varepsilon\rightarrow 0^+$.
We have the following crucial result.
\begin{theorem} {\bf(Degenerate monotone convergence result).}\label{MontnConvPP}
\\
Let $(X,\omega)$ be a polarized compact K\"ahler manifold of complex dimension $n$ and let $\gamma$, $T$ be closed positive $(1,1)$-currents with bounded local potentials. Then the following statements hold true.
\\
{\bf (A)} For all $\varphi\in {\cal P}\BT_{\gamma}$, $\varphi\leq 0$ and $k,l\geq 0$, $k+l\leq n$, $k\leq n-1$
$$
\int\limits_X-\varphi\,\gamma_{\varphi}^k\wedge T^l\wedge \omega^{n-k-l}<+\infty\,,\quad\mbox{and}
\quad 
\gamma_{\varphi}^{k+1}\wedge T^l=T^l\wedge \gamma_{\varphi}^{k+1}\,.
$$
{\bf (B)} Let $\varphi\in {\cal P}\BT_{\gamma}$, $\varphi\leq 0$ with zero Lelong numbers and $\varphi_{\varepsilon}\in {\cal P}_{\gamma+\varepsilon\omega}\cap C^{\infty}(X)$, such that $\varphi_{\varepsilon}\downarrow \varphi$ as
$\varepsilon\rightarrow 0^+$. Then for all $k,l\geq 0$, $k+l\leq n$, $k\leq n-1$
\begin{eqnarray}
\varphi_{\varepsilon}\,(\gamma_{\varphi_{\varepsilon}}+\varepsilon\omega)^k\wedge T^l \longrightarrow \varphi\,\gamma_{\varphi}^k\wedge T^l\,,\label{Mcv1}
\\\nonumber
\\
(\gamma_{\varphi_{\varepsilon}}+\varepsilon\omega)^{k+1}\wedge T^l
\longrightarrow \gamma_{\varphi}^{k+1}\wedge T^l\,,\label{Mcv2}
%\\\nonumber
%\\
%\varphi\,(\gamma_{\varphi_{\varepsilon}}+\varepsilon\omega)^k\wedge T^l \longrightarrow \varphi\,\gamma_{\varphi}^k\wedge T^l\,,
%\\\nonumber
%\\
%\gamma_{\varphi}\wedge (\gamma_{\varphi_{\varepsilon}}+\varepsilon\omega)^k\wedge T^l \longrightarrow \gamma_{\varphi}^{k+1}\wedge T^l\,,
\end{eqnarray}
weakly as $\varepsilon\rightarrow 0^+$. %Moreover $\gamma_{\varphi}^k\wedge T^l=T^l\wedge \gamma_{\varphi}^k$ for all $\varphi\in {\cal P}\BT_{\gamma}$ and $k,l\geq 0$, $k+l\leq n$.
\\
{\bf (C)} Let $\varphi\in {\cal P}\BT_{\gamma}$, $\varphi\leq 0$ and $\varphi_{\varepsilon}\in {\cal P}_{\gamma+R\omega}\cap C^{\infty}(X)$ such that $\varphi_{\varepsilon}\downarrow \varphi$ as
$\varepsilon\rightarrow 0^+$. Then for all $k,l\geq 0$, $k+l\leq n$, $k\leq n-1$
\begin{eqnarray}
\varphi_{\varepsilon}\,(\gamma_{\varphi_{\varepsilon}}+R\omega)^k\wedge T^l \longrightarrow \varphi\,(\gamma_{\varphi}+R\omega)^k\wedge T^l\,,\label{Mcv1B}
\\\nonumber
\\
(\gamma_{\varphi_{\varepsilon}}+R\omega)^{k+1}\wedge T^l
\longrightarrow (\gamma_{\varphi}+R\omega)^{k+1}\wedge T^l\,,\label{Mcv2B}
\end{eqnarray}
weakly as $\varepsilon\rightarrow 0^+$.
\end{theorem}
As follows immediately from the proof, the statement of this theorem still holds if we replace $T^l$ with a product $T_1\wedge ....\wedge T_l$, where the currents $T_j$ have the same properties as~$T$.  As a matter of fact, we wrote the statement in the previous special case only for the sake of notation simplicity. However, in the course of the proof, it is useful to notice that statements concerning terms involving $T^l$ are still valid if we replace  $T^l$ with $\gamma^r\wedge T^{l-r}$.
\\
\\
\emph{Proof}. Statement \eqref{Mcv2} follows from \eqref{Mcv1} by using the weak continuity of the $i\partial\bar\partial$ operator. The argument for Statement (B) is the same as for~(C).
\\
\\
{\bf Proof of (A).} 
We denote by $A_{k,l}$ the special case of Statement (A) in the theorem for the relative indices $(k,l)$.
We prove Statements $A_{k,l}$, $l=0,...,n-k$ by using an induction on $k=0,...,n-1$. We remark that Claim \ref{StpdCv} asserts Statement (A) in full generality for $k=0$. So we assume Statement $A_{k-1,\bullet}$ and we prove $A_{k,l}$, $l=0,...,n-k$ by using an induction on $l$. We remark that $A_{k,0}$ holds by the hypothesis $\varphi\in {\cal P}\BT_{\gamma}$. So we assume $A_{k,l}$ and we prove $A_{k,l+1}$. In fact let $\varphi_c:=\max\{\varphi,c\}\in {\cal P}_{\gamma}$, $c\in \R_{<0}$. 
By the regularization result in \cite{Dem4} let $(\varphi_{c,\varepsilon})_{\varepsilon>0}$, $\varphi_{c,\varepsilon}\in {\cal P}_{\gamma+\varepsilon\omega}\cap C^{\infty}(X)$ such that $\varphi_{c,\varepsilon}\downarrow \varphi_c$ as $\varepsilon\rightarrow 0^+$ and write $T=\theta+i\partial\bar\partial u$, with $\theta$ smooth, $\theta\leq K\omega$ and $u$ bounded with $\inf_X u=0$.
By using the monotone convergence theorem, the symmetry of the wedge product provided by the inductive hypothesis in $k$ and Stokes' formula, we expand the integral 
\begin{eqnarray*}
&&\int\limits_X-\varphi\,\gamma_{\varphi}^k\wedge T^{l+1}\wedge\omega^{n-k-l-1}
\\
\\
&=& 
\lim_{c\rightarrow -\infty}\;\lim_{\varepsilon\rightarrow 0^+}\;\int\limits_X-\varphi_{c,\varepsilon}\,\gamma_{\varphi}^k\wedge T^{l+1}\wedge\omega^{n-k-l-1}
\\
\\
&=&
\lim_{c\rightarrow -\infty}\;\lim_{\varepsilon\rightarrow 0^+}\;\int\limits_X-\varphi_{c,\varepsilon}\,T^{l+1}\wedge\gamma_{\varphi}^k\wedge \omega^{n-k-l-1}
\\
\\
&=&
\lim_{c\rightarrow -\infty}\;\lim_{\varepsilon\rightarrow 0^+}\,
\int\limits_X-\varphi_{c,\varepsilon}\, \theta\wedge T^{l}\wedge \gamma_{\varphi}^k\wedge \omega^{n-k-l-1}
\\
\\
&-&
\lim_{c\rightarrow -\infty}\;\lim_{\varepsilon\rightarrow 0^+}\,
\int\limits_X \varphi_{c,\varepsilon}\,i\partial\bar\partial u\wedge T^{l}\wedge \gamma_{\varphi}^k\wedge\omega^{n-k-l-1}
\\
\\
&\leq&
\lim_{c\rightarrow -\infty}\;\lim_{\varepsilon\rightarrow 0^+}\,\int\limits_X-\varphi_{c,\varepsilon}\,\gamma_{\varphi}^k\wedge T^{l}\wedge K\omega^{n-k-l}
\\
\\
&-&
\lim_{c\rightarrow -\infty}\;\lim_{\varepsilon\rightarrow 0^+}\,
\int\limits_Xu\,i\partial\bar\partial\varphi_{c,\varepsilon}\wedge \gamma_{\varphi}^k\wedge T^{l}\wedge\omega^{n-k-l-1}
\\
\\
&=&
K\int\limits_X-\varphi\,\gamma_{\varphi}^k\wedge T^{l}\wedge \omega^{n-k-l}
\\
\\
&-&
\lim_{c\rightarrow -\infty}\;\lim_{\varepsilon\rightarrow 0^+}\,
\int\limits_Xu\,(\gamma_{\varphi_{c,\varepsilon}}+\varepsilon\omega)\wedge \gamma_{\varphi}^k\wedge T^{l}\wedge\omega^{n-k-l-1}
\\
\\
&+&
\lim_{c\rightarrow -\infty}\;\lim_{\varepsilon\rightarrow 0^+}\,\int\limits_Xu\,(\gamma+\varepsilon\omega) \wedge \gamma_{\varphi}^k\wedge T^{l}\wedge\omega^{n-l-1}
\\
\\
&\leq&
K\int\limits_X-\varphi\,\gamma_{\varphi}^k\wedge T^{l}\wedge \omega^{n-k-l}
+
\sup_X u\,\int\limits_X\gamma\wedge \gamma_{\varphi}^k\wedge T^{l}\wedge\omega^{n-l-1}<+\infty\,,
\end{eqnarray*}
by the inductive hypothesis in $l$. We now prove the symmetry relation 
\begin{eqnarray}\label{SymExt}
\gamma_{\varphi}^{k+1}\wedge T^l=T^l\wedge \gamma_{\varphi}^{k+1}\,.
\end{eqnarray}
The decreasing monotone convergence theorem implies
$$
\lim_{c\rightarrow -\infty}\,
\int\limits_X(\varphi_{c}-\varphi)\,\gamma_{\varphi}^k\wedge T^l\wedge \omega^{n-k-l}=0\,,
$$
which means the convergence of the mass $\|(\varphi_{c}-\varphi)\,\gamma_{\varphi}^k\wedge T^l\|_{\omega}(X)\rightarrow 0$ as $c\rightarrow -\infty$. In particular 
$$
\varphi_{c}\,\gamma_{\varphi}^k\wedge T^l\longrightarrow \varphi \,\gamma_{\varphi}^k\wedge T\,,
$$ 
weakly as $c\rightarrow -\infty$.  So by the weak continuity of the $i\partial\bar\partial$ operator we deduce 
\begin{eqnarray}
\gamma_{\varphi_{c}}\wedge \gamma_{\varphi}^k\wedge T^l \longrightarrow \gamma_{\varphi}^{k+1}\wedge T^l\,, \label{CnvMS1}
\end{eqnarray}
weakly as $c\rightarrow -\infty$. 
The symmetry of the wedge product provided by the inductive hypothesis in $k$ implies 
$$
\gamma_{\varphi_{c}}\wedge \gamma_{\varphi}^k\wedge T^l
\;=\;
\gamma_{\varphi_{c}}\wedge T^l \wedge \gamma_{\varphi}^k
\;=\;
T^l \wedge \gamma_{\varphi_{c}}\wedge \gamma_{\varphi}^k 
$$
By the other hand $\eqref{CnvMS1}_{k,0}$ combined with the weak continuity of the $i\partial\bar\partial$ operator implies, by an induction on $l$
$$
T^l\wedge \gamma_{\varphi_{c}} \wedge \gamma_{\varphi}^k \longrightarrow T^l\wedge \gamma_{\varphi}^{k+1}\,,
$$ 
weakly as $c\rightarrow -\infty$. This combined with $\eqref{CnvMS1}$ implies the required symmetry \eqref{SymExt}.
\\
\\
{\bf Proof of (B).}
For all $k=0,...,n-1$ and $l=0,...,n-k$ we consider the following statement $B_{k,l}$: for all $p=0,...,k$
\begin{eqnarray}
&&\varphi_{\varepsilon}\,\gamma_{\varphi}^p\wedge(\gamma_{\varphi_{\varepsilon}}+\varepsilon\omega)^{k-p} \wedge T^l \longrightarrow \varphi\,\gamma_{\varphi}^k\wedge T^l\,,\label{Mcv3}
\\\nonumber
\\
&&i\partial\bar\partial\varphi_{\varepsilon}\wedge \gamma_{\varphi}^p\wedge(\gamma_{\varphi_{\varepsilon}}+\varepsilon\omega)^{k-p} \wedge T^l \longrightarrow i\partial\bar\partial\varphi\wedge\gamma_{\varphi}^k\wedge T^l\,,\label{Mcv31}
\\\nonumber
\\
&&\gamma_{\varphi}^p\wedge(\gamma_{\varphi_{\varepsilon}}+\varepsilon\omega)^{k-p+1}\wedge T^l
\longrightarrow \gamma_{\varphi}^{k+1}\wedge T^l\,,\label{Mcv32}
\\\nonumber
\\
&&\varphi\,\gamma_{\varphi}^p\wedge (\gamma_{\varphi_{\varepsilon}}+\varepsilon\omega)^{k-p}\wedge T^l
\longrightarrow \varphi\,\gamma_{\varphi}^k\wedge T^l\,,\label{Mcv4}
%\\\nonumber
%\\
%&&i\partial\bar\partial\varphi\wedge\gamma_{\varphi}^p\wedge(\gamma_{\varphi_{\varepsilon}}+\varepsilon\omega)^{k-p}\wedge T^l
%\longrightarrow i\partial\bar\partial\varphi\wedge\gamma_{\varphi}^k\wedge T^l\,,\label{Mcv41}
%\\\nonumber
%\\
%&&\gamma_{\varphi}^{p+1}\wedge (\gamma_{\varphi_{\varepsilon}}+\varepsilon\omega)^{k-p}\wedge T^l
%\longrightarrow \gamma_{\varphi}^{k+1}\wedge T^l \,,\label{Mcv42}
\end{eqnarray}
weakly as $\varepsilon\rightarrow 0^+$. 
We remark that \eqref{Mcv31} follows from \eqref{Mcv3} %and \eqref{Mcv41},  \eqref{Mcv42} follows from \eqref{Mcv4} 
by the weak continuity of the $i\partial\bar\partial$ operator. By combining \eqref{Mcv31} with the weak continuity of the $i\partial\bar\partial$ operator we obtain
$$
(\gamma_{\varphi_{\varepsilon}}+\varepsilon\omega)\wedge
\gamma_{\varphi}^{p}\wedge (\gamma_{\varphi_{\varepsilon}}+\varepsilon\omega)^{k-p}\wedge T^l
\longrightarrow \gamma_{\varphi}^{k+1}\wedge T^l \,,
$$
weakly as $\varepsilon\rightarrow 0^+$.
On the other hand the symmetry of the wedge product proved in part (A) of the theorem implies
\begin{eqnarray*}
(\gamma_{\varphi_{\varepsilon}}+\varepsilon\omega)\wedge
\gamma_{\varphi}^{p}\wedge (\gamma_{\varphi_{\varepsilon}}+\varepsilon\omega)^{k-p}\wedge T^l
&=&
(\gamma_{\varphi_{\varepsilon}}+\varepsilon\omega)^{k-p+1}\wedge T^l\wedge
\gamma_{\varphi}^{p}
\\
\\
&=&
\gamma_{\varphi}^{p}\wedge (\gamma_{\varphi_{\varepsilon}}+\varepsilon\omega)^{k-p+1}\wedge T^l\,.
\end{eqnarray*}
In this way we deduce \eqref{Mcv32}.
The statements $B_{0,\bullet}$ are true by the proof of Claim \ref{StpdCv}. 
We now prove by induction on $k=0,...,n-1$ that Statements $B_{k,l}$, $l=0,...,n-k$ hold true. In fact we prove the following claim.
\begin{claim}
If $B_{j,\bullet}$ holds true for all $j=0,...,k-1$, then $B_{k,l}$ holds also true for all $l=0,...,n-k$.
\end{claim}
As pointed out before in order to prove $B_{k,l}$ it is sufficient to show \eqref{Mcv3} and \eqref{Mcv4}. The proof of \eqref{Mcv4} is quite similar to the proof of \eqref{Mcv3} that we now explain. We first prove by induction on $s=0,...,k-p$ the inequality
\begin{eqnarray}\label{CvIndI}
&&\int\limits_X-\varphi_{\varepsilon} \,\gamma_{\varphi}^p\wedge (\gamma_{\varphi_{\varepsilon}}+\varepsilon\omega)^{k-p}\wedge T^l\wedge \omega^{n-k-l}\nonumber
\\\nonumber
\\
&\leq&
\int\limits_X-\varphi \,\gamma_{\varphi}^{p+s}\wedge (\gamma_{\varphi_{\varepsilon}}+\varepsilon\omega)^{k-p-s}\wedge T^l\wedge \omega^{n-k-l}\nonumber
\\\nonumber
\\
&+&
\sum_{r=0}^{s-1}\;\int\limits_X(\varphi_{\varepsilon}-\varphi) \,\gamma_{\varphi}^{p+r}\wedge (\gamma_{\varphi_{\varepsilon}}+\varepsilon\omega)^{k-p-r-1}\wedge \gamma\wedge T^l\wedge \omega^{n-k-l}\nonumber
\\\nonumber
\\
&-&
\sum_{r=0}^{s-1}\;\int\limits_X\varepsilon\varphi \,\gamma_{\varphi}^{p+r}\wedge (\gamma_{\varphi_{\varepsilon}}+\varepsilon\omega)^{k-p-r-1}\wedge T^l\wedge \omega^{n-k-l+1}.
\end{eqnarray}
Inequality \eqref{CvIndI} is obviously true for $s=0$. (Here we adopt the usual convention of neglecting a sum when it runs over an empty set of indices.) Before proceding to the proof of the inequality \eqref{CvIndI}, we need to point out two useful facts.
\smallskip

\noindent
{\bf 1)} Let $\alpha$ be a smooth closed real $(q,q)$-form, $R$ be a closed positive $(r,r)$-current, $v\geq 0$ be a measurable function such that $\int_XvR\wedge \omega^{n-r}<+\infty$. This implies that the currents $i\partial\bar\partial v\wedge R:=i\partial\bar\partial (v\, R)$ and 
$i\partial\bar\partial v\wedge \alpha\wedge R:=i\partial\bar\partial (v \alpha\wedge R)$
are well defined. Then the Leibniz formula implies
\begin{eqnarray}\label{MltSymWeg}
\alpha\wedge i\partial\bar\partial v\wedge R=i\partial\bar\partial v\wedge \alpha\wedge R\,.
\end{eqnarray}
{\bf 2)} Thanks to part (A) of the theorem we have
$$
\int\limits_X-\varphi \,\gamma_{\varphi}^{p+r}\wedge \gamma^h\wedge T^l\wedge \omega^{n-p-r-h-l}<+\infty
$$
for all $h=0,...,k-p-r-1$. By \eqref{MltSymWeg} this implies
$$
\int\limits_X-\varphi \,\gamma_{\varphi}^{p+r}\wedge (\gamma_{\varphi_{\varepsilon}}+\varepsilon\omega)^{k-p-r-1}\wedge T^l\wedge \omega^{n-k-l+1}<+\infty\,,
$$
so the current 
$$
S:=\varphi\,\gamma_{\varphi}^{p+r}\wedge (\gamma_{\varphi_{\varepsilon}}+\varepsilon\omega)^{k-p-r-1}\wedge T^l
$$ 
is well defined and we can define the current 
$$
i\partial\bar\partial \varphi\wedge 
\gamma_{\varphi}^{p+r}\wedge (\gamma_{\varphi_{\varepsilon}}+\varepsilon\omega)^{k-p-r-1}\wedge T^l
:=i\partial\bar\partial S\,.
$$ 
Then the integration by parts formula
$$
\int\limits_X i\partial\bar\partial\varphi_{\varepsilon} \wedge S\wedge \omega^{n-k-l}\;=\;\int\limits_X \varphi_{\varepsilon}\, i\partial\bar\partial S \wedge \omega^{n-k-l}
$$
can be written explicitly as
\begin{eqnarray}\label{IntPart}
&&\int\limits_X i\partial\bar\partial\varphi_{\varepsilon} \wedge 
\varphi\,\gamma_{\varphi}^{p+r}\wedge (\gamma_{\varphi_{\varepsilon}}+\varepsilon\omega)^{k-p-r-1}\wedge T^l
\wedge\omega^{n-k-l}\nonumber
\\\nonumber
\\
&=&
\;\int\limits_X \varphi_{\varepsilon}\, i\partial\bar\partial \varphi\wedge 
\gamma_{\varphi}^{p+r}\wedge (\gamma_{\varphi_{\varepsilon}}+\varepsilon\omega)^{k-p-r-1}\wedge T^l
\wedge \omega^{n-k-l}\,.
\end{eqnarray}
We suppose now the inequality \eqref{CvIndI} true for $s$ and we prove it for $s+1$. We start by expanding, thanks to formula \eqref{MltSymWeg}, the integral 
%of the main $(k-1)$-inductive hypothesis $\eqref{SymWeg}_{j,\,\bullet}$ , $j\leq k-1$ and 
\begin{eqnarray*}
I
&:=&\int\limits_X-\varphi \,\gamma_{\varphi}^{p+s}\wedge (\gamma_{\varphi_{\varepsilon}}+\varepsilon\omega)^{k-p-s}\wedge T^l\wedge \omega^{n-k-l}
\\
\\
&=&
\int\limits_X-\varphi \,\gamma_{\varphi}^{p+s}
\wedge
(\gamma+\varepsilon\omega)
\wedge (\gamma_{\varphi_{\varepsilon}}+\varepsilon\omega)^{k-p-s-1}\wedge T^l\wedge \omega^{n-k-l}
\\
\\
&+&
\int\limits_X-\varphi\,\gamma_{\varphi}^{p+s}\wedge 
i\partial\bar\partial \varphi_{\varepsilon}\wedge
(\gamma_{\varphi_{\varepsilon}}+\varepsilon\omega)^{k-p-s-1}\wedge T^l\wedge \omega^{n-k-l}
\\
\\
&=&
\int\limits_X-\varepsilon\varphi \,\gamma_{\varphi}^{p+s}\wedge (\gamma_{\varphi_{\varepsilon}}+\varepsilon\omega)^{k-p-s-1}\wedge T^l\wedge \omega^{n-k-l+1}
\\
\\
&-&
\int\limits_X\varphi \,\gamma_{\varphi}^{p+s}\wedge (\gamma_{\varphi_{\varepsilon}}+\varepsilon\omega)^{k-p-s-1}\wedge \gamma\wedge T^l\wedge \omega^{n-k-l}
\\
\\
&-&
\int\limits_X i\partial\bar\partial \varphi_{\varepsilon}\wedge\varphi\,\gamma_{\varphi}^{p+s}\wedge 
(\gamma_{\varphi_{\varepsilon}}+\varepsilon\omega)^{k-p-s-1}\wedge T^l\wedge \omega^{n-k-l}\,.
\end{eqnarray*}
By applying the integration by parts formula \eqref{IntPart} to the last integral we deduce
\begin{eqnarray*}
I
&=&
\int\limits_X-\varphi_{\varepsilon} \,\gamma_{\varphi}^{p+s+1}\wedge (\gamma_{\varphi_{\varepsilon}}+\varepsilon\omega)^{k-p-s-1}\wedge T^l\wedge \omega^{n-k-l}
\\
\\
&+&
\int\limits_X\varphi_{\varepsilon} \,\gamma\wedge\gamma_{\varphi}^{p+s}\wedge (\gamma_{\varphi_{\varepsilon}}+\varepsilon\omega)^{k-p-s-1}\wedge T^l\wedge \omega^{n-k-l}
\\
\\
&-&
\int\limits_X\varphi \,\gamma_{\varphi}^{p+s}\wedge (\gamma_{\varphi_{\varepsilon}}+\varepsilon\omega)^{k-p-s-1}\wedge\gamma\wedge T^l\wedge \omega^{n-k-l}
\\
\\
&-&
\int\limits_X\varepsilon\varphi \,\gamma_{\varphi}^{p+s}\wedge (\gamma_{\varphi_{\varepsilon}}+\varepsilon\omega)^{k-p-s-1}\wedge T^l\wedge \omega^{n-k-l+1}\,.
\end{eqnarray*}
By combining the symmetry of the wedge product proved in part (A)
with formula \eqref{MltSymWeg} we get
\begin{eqnarray*}
\gamma\wedge\gamma_{\varphi}^{p+s}\wedge (\gamma_{\varphi_{\varepsilon}}+\varepsilon\omega)^{k-p-s-1}\wedge T^l
&=&
\gamma\wedge(\gamma_{\varphi_{\varepsilon}}+\varepsilon\omega)^{k-p-s-1}\wedge T^l\wedge\gamma_{\varphi}^{p+s}
\\
\\
&=&
(\gamma_{\varphi_{\varepsilon}}+\varepsilon\omega)^{k-p-s-1}\wedge 
\gamma\wedge T^l\wedge\gamma_{\varphi}^{p+s}
\\
\\
&=&
\gamma_{\varphi}^{p+s}\wedge 
(\gamma_{\varphi_{\varepsilon}}+\varepsilon\omega)^{k-p-s-1}\wedge 
\gamma\wedge T^l\,.
\end{eqnarray*}
By plugging this into the previous expression of $I$ we obtain
\begin{eqnarray*}
I
&=&
\int\limits_X-\varphi_{\varepsilon} \,\gamma_{\varphi}^{p+s+1}\wedge (\gamma_{\varphi_{\varepsilon}}+\varepsilon\omega)^{k-p-s-1}\wedge T^l\wedge \omega^{n-k-l}
\\
\\
&+&
\int\limits_X(\varphi_{\varepsilon}-\varphi) \,\gamma_{\varphi}^{p+s}\wedge (\gamma_{\varphi_{\varepsilon}}+\varepsilon\omega)^{k-p-s-1}\wedge\gamma\wedge T^l\wedge \omega^{n-k-l}
\\
\\
&-&
\int\limits_X\varepsilon\varphi \,\gamma_{\varphi}^{p+s}\wedge (\gamma_{\varphi_{\varepsilon}}+\varepsilon\omega)^{k-p-s-1}\wedge T^l\wedge \omega^{n-k-l+1}\,,
\end{eqnarray*}
which implies inequality \eqref{CvIndI} for $s+1$. For $s=k-p$ the inequality \eqref{CvIndI} rewrites as
\begin{eqnarray*}
&&\int\limits_X-\varphi_{\varepsilon} \,\gamma_{\varphi}^p\wedge (\gamma_{\varphi_{\varepsilon}}+\varepsilon\omega)^{k-p}\wedge T^l\wedge \omega^{n-k-l}
\;
\leq 
\;
\int\limits_X-\varphi \,\gamma_{\varphi}^k\wedge T^l\wedge \omega^{n-k-l}
\\
\\
&+&
\sum_{r=p}^{k-1}\;\int\limits_X(\varphi_{\varepsilon}-\varphi) \, \gamma_{\varphi}^r\wedge (\gamma_{\varphi_{\varepsilon}}+\varepsilon\omega)^{k-r-1}\wedge \gamma\wedge T^l\wedge \omega^{n-k-l}
\\
\\
&-&
\sum_{r=p}^{k-1}\;\int\limits_X\varepsilon\varphi \, \gamma_{\varphi}^r\wedge (\gamma_{\varphi_{\varepsilon}}+\varepsilon\omega)^{k-r-1}\wedge T^l\wedge \omega^{n-k-l+1}\,.
\end{eqnarray*}
By using the inductive convergence hypothesis  $\eqref{Mcv3}_{j,\bullet}$, $\eqref{Mcv4}_{j,\bullet}$ in $B_{j,\bullet}$ for $j\leq k-1$ we deduce
\begin{eqnarray}\label{CmpctMass}
&&\limsup_{\varepsilon\rightarrow 0^+}\int\limits_X-\varphi_{\varepsilon} \,\gamma_{\varphi}^p\wedge (\gamma_{\varphi_{\varepsilon}}+\varepsilon\omega)^{k-p}\wedge T^l\wedge \omega^{n-k-l}\nonumber
\\\nonumber
\\
&&\leq 
\int\limits_X-\varphi \,\gamma_{\varphi}^k\wedge T^l\wedge \omega^{n-k-l}<+\infty\,,
\end{eqnarray}
by Statement A. (We can always arrange $\varphi_{\varepsilon}\leq 0$ for all $\varepsilon\in (0,1)$ by changing $\varphi$ into $\varphi-C$.) Thus by weak compactness of the mass there exists a sequence $(\varepsilon_j)_j$, $\varepsilon_j\downarrow 0^+$ and a current of order zero $\Theta\in {\cal D}'_{n-k-l,n-k-l}(X)$ such that
$$
\varphi_{\varepsilon_j}\, \gamma_{\varphi}^p\wedge (\gamma_{\varphi_{\varepsilon_j}}+\varepsilon_j\,\omega)^{k-p}\wedge T^l\longrightarrow \Theta\,,
$$
weakly as $j\rightarrow +\infty$. So for any smooth and strongly positive form $\alpha$ of bidegree $(n-k-l,n-k-l)$, we have
$$
\varphi_{\varepsilon_j} \,\gamma_{\varphi}^p\wedge (\gamma_{\varphi_{\varepsilon_j}}+\varepsilon_j\,\omega)^{k-p}\wedge T^l\wedge\alpha
\longrightarrow \Theta\wedge\alpha\,,
$$
weakly as $j\rightarrow +\infty$. The fact that $\varphi_{\varepsilon_j}\downarrow\varphi$ and 
$$
\gamma_{\varphi}^p\wedge (\gamma_{\varphi_{\varepsilon_j}}+\varepsilon_j\,\omega)^{k-p}\wedge T^l\wedge\alpha
\longrightarrow 
\gamma_{\varphi}^k\wedge T^l\wedge\alpha\,,
$$
weakly as $j\rightarrow +\infty$, by the convergence inductive hypothesis $\eqref{Mcv32}_{k-1,l}$, implies 
$$
\Theta\wedge\alpha\leq \varphi\,\gamma_{\varphi}^k\wedge T^l\wedge\alpha\,,
$$ 
thanks to Lemma (3.9), page 189 in \cite{Dem2}. Thus $\Theta\leq \varphi\,\gamma_{\varphi}^k\wedge T^l$ . Combining this with the inequality  \eqref{CmpctMass} we obtain
\begin{eqnarray*}
\int\limits_X\Theta\wedge \omega^{n-k-l}
&\leq &
\int\limits_X\varphi\,\gamma_{\varphi}^k\wedge T^l\wedge \omega^{n-k-l}
\\
\\
&\leq&
\liminf_{\varepsilon\rightarrow 0^+}\int\limits_X\varphi_{\varepsilon} \,\gamma_{\varphi}^p\wedge (\gamma_{\varphi_{\varepsilon}}+\varepsilon\omega)^{k-p}\wedge T^l\wedge \omega^{n-k-l}
\\
\\
&\leq&
\lim_{j\rightarrow +\infty}\int\limits_X\varphi_{\varepsilon_j} \,\gamma_{\varphi}^p\wedge (\gamma_{\varphi_{\varepsilon_j}}+\varepsilon_j\,\omega)^{k-p}\wedge T^l\wedge \omega^{n-k-l}
\\
\\
&=&\int\limits_X\Theta\wedge \omega^{n-k-l}\,.
\end{eqnarray*}
We deduce $\Tr_{\omega}(\varphi\,\gamma_{\varphi}^k\wedge T^l-\Theta)=0$, which implies $\varphi\,\gamma_{\varphi}^k\wedge T^l=\Theta$ since $0\leq \varphi\,\gamma_{\varphi}^k\wedge T^l-\Theta$. This proves Statement $B_{k,l}$.\hfill$\Box$
\\
\\
We introduce also the subsets 
\begin{eqnarray*}
\hat{\cal P}\BT_{\gamma}&:=&\left\{\varphi\in {\cal P}\BT^0_{\gamma}\,\mid\,\int_X-\varphi\,\gamma_{\varphi}^n<+\infty \right\}+\R\subset {\cal P}\BT_{\gamma}\,,
\\
\\
\hat{\cal P}\BT^0_{\gamma}&:=&\{\varphi\in \hat{\cal P}\BT_{\gamma}\,\mid\,\sup_X\varphi=0 \}\,.
\end{eqnarray*}
Without changes in the proof of Theorem \ref{MontnConvPP} we get the following corollary.
\begin{corollary}\label{Cor-MontnConvPP}
For all  $\varphi\in \hat{\cal P}\BT_{\gamma}$, $\varphi\leq 0$, the assertions A$)$, B$)$ and C$)$ of Theorem \ref{MontnConvPP} hold for all $k=0,...,n$.
\end{corollary}
Let now $\Theta$ be a closed positive $(n-1,n-1)$-current and consider the $L^2$-space 
\begin{eqnarray*}
L^2(X,\Theta):=\left\{\alpha\in \Gamma(X,\Lambda^{1,0}T_X^*)\;\mid\; \int\limits_X i\alpha\wedge \bar\alpha\wedge \Theta<+\infty \right\}_{\Big/\Theta-a.e} \,,
\end{eqnarray*}
equipped with the hermitian product $\left<\alpha,\beta\right>_{\Theta}:=\int_X i\alpha\wedge \bar\beta\wedge \Theta$, which is well defined by the polarization identity. The $\Theta$-almost everywhere equivalence relation is defined by$\,$: $\alpha\sim\beta$ iff 
$$
\int\limits_X i(\alpha-\beta)\wedge \overline{(\alpha-\beta)}\wedge \Theta=0\,.
$$
The subscript ''$\Theta$-a.e.'' in the definition of $L^2(X,\Theta)$ above is ''$\Theta$-almost everywhere.
Let $\alpha_k,\,\alpha\in L^2(X,\Theta)$. We say that the sequence $\alpha_k$ converges $L^2(X,\Theta)$-weakly to $\alpha$ if 
$$
\int_X i\alpha\wedge \bar\beta\wedge \Theta=\lim_{k\rightarrow +\infty}\int_X i\alpha_k\wedge \bar\beta\wedge \Theta\,,
$$
for all $\beta\in L^2(X,\Theta)$. Let $\varphi\in {\cal P}^0_{\gamma}$ such that $\int_X-\varphi\,\Theta\wedge\, \omega<+\infty$. Then one can define $\partial \varphi\wedge \Theta:=\partial (\varphi\Theta)$. We write $\partial \varphi\in L^2(X,\Theta)$ if there exists $\alpha\in L^2(X,\Theta)$ such that $\partial (\varphi\Theta)=\alpha\wedge \Theta$ in the sense of currents. In this case we write
$$
\int\limits_X i\partial \varphi\wedge \bar\partial \varphi\wedge \Theta:=\int\limits_X i\alpha\wedge \bar\alpha\wedge \Theta\,.
$$
With these notations we have the following corollary.
%% of Theorem \ref{MontnConvPP}.

\begin{corollary}\label{ConvGradMA}
Let $(X,\omega)$ be a polarized compact K\"ahler manifold of complex dimension $n$ and let $\gamma$, $T$ be closed positive $(1,1)$-currents with bounded local potentials, let $\Theta$ be a closed positive $(n-1,n-1)$-current and consider $\varphi\in \hat{\cal P}\BT_{\gamma}$, $\varphi\leq 0$, $\psi\in {\cal P}_{\gamma}\cap L^{\infty}(X)$, $\psi\leq 0$. Then for all $k,l\geq 0$, $k+l\leq n-1$, 
\begin{eqnarray}
&&\int\limits_X i\partial\varphi\wedge \bar\partial\varphi\wedge \gamma_{\varphi}^k\wedge T^l\wedge \omega^{n-k-l-1}<+\infty\,,\label{CvGrad}
\\\nonumber
\\
&&
\int\limits_X i\partial\psi\wedge \bar\partial\psi\wedge\Theta<+\infty\,.\label{C0-CvGrad}
\end{eqnarray}
Moreover let $(\varphi_{\varepsilon})_{\varepsilon>0}$, $(\psi_{\varepsilon})_{\varepsilon>0}\subset C^{\infty}(X)$, 
$\varphi_{\varepsilon}\in {\cal P}_{\gamma+R\omega}$, $\psi_{\varepsilon}\in {\cal P}_{\gamma+\varepsilon\omega}$ such that $\varphi_{\varepsilon}\downarrow \varphi$, $\psi_{\varepsilon}\downarrow \psi$ as $\varepsilon\rightarrow 0^+$. Then
\begin{eqnarray}
&&\lim_{\varepsilon\rightarrow 0^+}\,\int\limits_X i\partial(\varphi_{\varepsilon}-\varphi)\wedge \bar\partial(\varphi_{\varepsilon}-\varphi)\wedge \gamma_{\varphi}^k\wedge T^l\wedge \omega^{n-k-l-1}=0\,,\label{L2-CVGrad}
\\\nonumber
\\
&&
\lim_{\varepsilon\rightarrow 0^+}\,\int\limits_X i\partial(\psi_{\varepsilon}-\psi)\wedge \bar\partial(\psi_{\varepsilon}-\psi)\wedge\Theta=0\,.\label{C0-L2-CVGrad}
\end{eqnarray}
\end{corollary}
%By puting $k=l=0$ in the \eqref{CvGrad} we deduce ${\cal P}\BT_{\gamma}\subset H^1(X)$.
\emph{Proof}. By integrating by parts we obtain
\begin{eqnarray*}
&&\int\limits_X i\partial\varphi_{\varepsilon}\wedge \bar\partial\varphi_{\varepsilon}\wedge \gamma_{\varphi}^k\wedge T^l\wedge \omega^{n-k-l-1} 
\\
\\
&=&
-\int\limits_X \varphi_{\varepsilon}\, i\partial\bar\partial\varphi_{\varepsilon}\wedge \gamma_{\varphi}^k\wedge T^l\wedge \omega^{n-k-l-1} 
\end{eqnarray*}
\begin{eqnarray*}
&=&
\int\limits_X \varphi_{\varepsilon}\, (\gamma+R\omega)\wedge \gamma_{\varphi}^k\wedge T^l\wedge \omega^{n-k-l-1} 
\\
\\
&-&
\int\limits_X \varphi_{\varepsilon}\,(\gamma_{\varphi_{\varepsilon}}+R\omega)\wedge \gamma_{\varphi}^k\wedge T^l\wedge \omega^{n-k-l-1} \,.
\end{eqnarray*}
By the proof of Theorem \ref{MontnConvPP}, B we can take the limit, so
\begin{eqnarray}
0&\leq&
\lim_{\varepsilon\rightarrow 0^+}\,\int\limits_X i\partial\varphi_{\varepsilon}\wedge \bar\partial\varphi_{\varepsilon}\wedge \gamma_{\varphi}^k\wedge T^l\wedge \omega^{n-k-l-1} \nonumber
\\\nonumber
\\
&=&
\int\limits_X \varphi\,(\gamma-\gamma_{\varphi})\wedge \gamma_{\varphi}^k\wedge T^l\wedge \omega^{n-k-l-1} <+\infty\,.\label{MasCvGrad}
\end{eqnarray}
On the other hand the weak convergence of the sequence 
$$
\varphi_{\varepsilon}\,\gamma_{\varphi}^k\wedge T^l\wedge \omega^{n-k-l-1}\longrightarrow \varphi\,\gamma_{\varphi}^k\wedge T^l\wedge \omega^{n-k-l-1}\,,
$$ 
combined with the weak continuity of the $\partial$ operator implies
$$
\partial\varphi_{\varepsilon}\wedge \gamma_{\varphi}^k\wedge T^l\wedge \omega^{n-k-l-1}\longrightarrow \partial\varphi\wedge\gamma_{\varphi}^k\wedge T^l\wedge \omega^{n-k-l-1}\,,
$$ 
weakly as $\varepsilon\rightarrow 0^+$. Then the $L^2(X,\gamma_{\varphi}^k\wedge T^l\wedge \omega^{n-k-l-1})$-weak compactness 
provided by \eqref{MasCvGrad}
implies \eqref{CvGrad} and 
the $L^2(X,\gamma_{\varphi}^k\wedge T^l\wedge \omega^{n-k-l-1})$-weak convergence $\partial\varphi_{\varepsilon}\rightarrow \partial\varphi$ as $\varepsilon\rightarrow 0^+$. This implies
\begin{eqnarray*}
&&\int\limits_X i\partial\varphi\wedge \bar\partial\varphi\wedge \gamma_{\varphi}^k\wedge T^l\wedge \omega^{n-k-l-1} 
\\
\\
&=&
\lim_{\varepsilon\rightarrow 0^+}\,\int\limits_X i\partial\varphi_{\varepsilon}\wedge \bar\partial\varphi\wedge \gamma_{\varphi}^k\wedge T^l\wedge \omega^{n-k-l-1} 
\end{eqnarray*}
\begin{eqnarray*}
&=&
\lim_{\varepsilon\rightarrow 0^+}\,\int\limits_X -\varphi_{\varepsilon}\,i\partial\bar\partial\varphi\wedge \gamma_{\varphi}^k\wedge T^l\wedge \omega^{n-k-l-1}
\\
\\
&=&
\lim_{\varepsilon\rightarrow 0^+}\,\int\limits_X -\varphi_{\varepsilon}\,(\gamma-\gamma_{\varphi})\wedge \gamma_{\varphi}^k\wedge T^l\wedge \omega^{n-k-l-1}
\\
\\
&=&
\int\limits_X -\varphi\,(\gamma-\gamma_{\varphi})\wedge \gamma_{\varphi}^k\wedge T^l\wedge \omega^{n-k-l-1}
\\
\\
&=&
\lim_{\varepsilon\rightarrow 0^+}\,\int\limits_X i\partial\varphi_{\varepsilon}\wedge \bar\partial\varphi_{\varepsilon}\wedge \gamma_{\varphi}^k\wedge T^l\wedge \omega^{n-k-l-1} \,,
\end{eqnarray*}
by identity \eqref{MasCvGrad}. This implies \eqref{L2-CVGrad} by elementary facts about Hilbert spaces. The proof of \eqref{C0-CvGrad} and \eqref{C0-L2-CVGrad} is quite similar.\hfill$\Box$
\\
\\
The conclusion of the corollary \ref{ConvGradMA} still holds true if we replace the current $\gamma_{\varphi}^k\wedge T^l\wedge \omega^{n-k-l-1}$ with a sum of currents
$$
\Xi:=\sum_{k+l\leq n-1}C_{k,l}\,\gamma_{\varphi}^k\wedge T^l\wedge \omega^{n-k-l-1}\,,
$$
where $C_{k,l}\in \R$ such that $\Xi\geq 0$. We infer the linearity formula
\begin{eqnarray*}
\int\limits_X i\partial\varphi\wedge \bar\partial\varphi\wedge\Xi\;=\;\sum_{k+l\leq n-1}C_{k,l}\,\int\limits_X i\partial\varphi\wedge \bar\partial\varphi\wedge\gamma_{\varphi}^k\wedge T^l\wedge \omega^{n-k-l-1}\,.
\end{eqnarray*}
%%%%%%%%%%%%%%%%%%%%%%%%%%%%%%%%%%%%%%%%%%%%%%%%%%%%%%%%%%%%%%%%%%%%%%%%%%%%%%%%%%%%%%%%%%%%%%%%%%%%%%%%%%%%%%%%%%%%%%%%%%%%%%%%%%%%%%%%%%%%%%%%%%%%%%%%%%%%%%%%%%%%%%%%%%%%%%%%%%%%%%%%%%%%%%%%%%%%%%%%%%%%%%%%%%%%%%%%%%%%%%%%%%%%%%%%%%%%%%%%%%%%%%%%%%%%%%%%%%%%%%%%%%%%%%%%%%%%%%%%%%%%%%%%%%%%%%%%%%%%%%%%
%%%%%%%%%%%%%%%%%%%%%%%%%%%%%%%%%%%%%%%%%%%%%%%%%%%%%%%%%%%%%%%%%%%%%%%%%%%%%%%%%%%%%%%%%%%%%%%%%%%%%%%%%%%%%%%%%%%%%%%%%%%%%%%%%%%%%%%%%%%%%%%%%%%%%%%%%%%%%%%%%%%%%%%%%%%%%%%%%%%%%%%%%%%%%%%%%%%%%%%%%%%%%%%%%%%%%%%%%%%%%%%%%%%%%%%%%%%%%%%%%%%%%%%%%%%%%%%%%%%%%%%%%%%%%%%%%%%%%%%%%%%%%%%%%%%%%%%%%%%%%%%%
\section{Uniqueness of the solutions}
%%%%%%%%%%%%%%%%%%%%%%%%%%%%%%%%%%%%%%%%%%%%%%%%%%%%%%%%%%%%%%%%%%%%%%%%%%%%%%%%%%%%%%%%%%%%%%%%%%%%%%%%%%%%%%%%%%%%%%%%%%%%%%%%%%%%%%%%%%%%%%%%%%%%%%%%%%%%%%%%%%%%%%%%%%%%%%%%%%%%%%%%%%%%%%%%%%%%%%%%%%%%%%%%%%%%%%%%%%%%%%%%%%%%%%%%%%%%%%%%%%%%%%%%%%%%%%%%%%%%%%%%%%%%%%%%%%%%%%%%%%%%%%%%%%%%%%%%%%%%%%%%
%%%%%%%%%%%%%%%%%%%%%%%%%%%%%%%%%%%%%%%%%%%%%%%%%%%%%%%%%%%%%%%%%%%%%%%%%%%%%%%%%%%%%%%%%%%%%%%%%%%%%%%%%%%%%%%%%%%%%%%%%%%%%%%%%%%%%%%%%%%%%%%%%%%%%%%%%%%%%%%%%%%%%%%%%%%%%%%%%%%%%%%%%%%%%%%%%%%%%%%%%%%%%%%%%%%%%%%%%%%%%%%%%%%%%%%%%%%%%%%%%%%%%%%%%%%%%%%%%%%%%%%%%%%%%%%%%%%%%%%%%%%%%%%%%%%%%%%%%%%%%%%%
We start with a renormalization result for the density volume form of a big and nef $(1,1)$-cohomology class. This uses \cite{De-Pa} in a crucial way.

\begin{lemma}\label{C0-Kahl}
Let $X$ be a compact K\"ahler manifold of complex dimension~$n$, let $T$ be a big closed positive $(1,1)$-current with continuous local potentials. Then there exist a big closed positive $(1,1)$-current $\gamma$ with continuous local potentials, cohomologous to $T$ and a complex analytic subset $Z\subset X$ such that $\gamma$ is a smooth K\"ahler metric over $X\smallsetminus Z$.
\end{lemma}

\noindent
\emph{Proof}. Let $\alpha$ be a smooth closed $(1,1)$-form representing the cohomology class of~$T$. The assumption on $T$ means that we can write $T=\alpha+i\partial\bar\partial\psi\geq 0$ where $\psi$ is a continuous quasi-plurisubharmonic function. By the approximation theorem of \cite{Dem4}, there exists a decreasing sequence $\psi_j$ of smooth quasi-plurisubharmonic functions converging to $\psi$ such that 
$$\alpha+i\partial\bar\partial\psi_j\geq -\frac{1}{j}\,\omega,$$
in particular the class $\{T\}=\{\alpha\}$ is nef (i.e.\ numerically effective in the sense of \cite{Dem4}) and big. By Theorem 0.5 of \cite{De-Pa}, there also exists a K\"ahler current $\Theta=\alpha+i\partial\bar\partial\varphi\in\{\alpha\}$, with $\Theta\ge\varepsilon\omega$ (in the sense of currents) and $\varepsilon>0$, such that $\varphi$ has logarithmic poles on some analytic subset $Z\subset X$ and $\varphi$ is smooth on~$X\smallsetminus Z$. If we consider $\varphi_j:=\max(\varphi,\psi-j)$, then $\varphi_j$ is a continuous quasi-subharmonic function which coincides with $\varphi$ on a neighborhood of the compact set 
$$
K_j=\{z\in X\smallsetminus Z\,;\; \psi(z)-\varphi(z)\le j-1\}\,.
$$ 
Clearly we have $\alpha+i\partial\bar\partial\varphi_j\ge 0$ on $X$ and
$$
\alpha+i\partial\bar\partial\varphi_j=\alpha+i\partial\bar\partial\varphi=\Theta\ge \varepsilon\omega\quad\hbox{on a neighborhood of~$K_j$.}
$$
Therefore if we put $\Phi=\sum_{j\ge 1}2^{-j}\varphi_j$, then $\Phi$ is a continuous quasi-plurisub\-harmonic potential on $X$ (notice that there is uniform convergence since\break $-j-C_1\le\varphi_j\le C_2$ on $X$), and there exists a continuous function $\lambda\ge 0$ on $X$ such that
$\alpha+i\partial\bar\partial\Phi\ge\lambda$ and $\lambda(x)>0$ on $\bigcup K_j=X\smallsetminus Z$. By Richberg's approximation theorem applied on $X\smallsetminus Z$ (invoking e.g.\ \cite{Dem2}, Theorem (I.5.21) with the error function $\frac{1}{2}\lambda(x)>0$ on~$X\smallsetminus Z$), we can find a function 
$$
\Psi\in {\cal P}_{\alpha}\cap C^{\infty}(X\smallsetminus Z)\,,
$$
such that $\Phi\le\Psi\le\Phi+\frac{1}{2}\lambda$ and $\alpha+i\partial\bar\partial\Psi\ge\frac{1}{2}\lambda>0$ on $X\smallsetminus Z\,$; this implies that $\Psi$ has a continuous extension to $X$ such that $\Psi=\Phi$ on $Z$, and also that the extension satisfies $\gamma:=\alpha+i\partial\bar\partial\Psi\ge\frac{1}{2}\lambda\ge 0$ everywhere on $X$ by standard arguments of potential theory.\hfill$\Box$

\begin{theorem}\label{UniqSol}
Let $X$ be a compact connected K\"ahler manifold of complex dimension~$n$ and let $\gamma$ be a big closed positive $(1,1)$-current with continuous local potentials.
\\
{\bf (A)}  Let $\psi\in{\cal P}^0_{\gamma}\cap L^{\infty}(X)$ and $\varphi\in{\cal P}\BT^0_{\gamma}$ such that $$(\gamma+i\partial\bar\partial \psi)^n=(\gamma+i\partial\bar\partial \varphi)^n\,.$$ 
Then $\psi=\varphi$.
\\
{\bf (B)} Let $\psi, \varphi\in{\cal P}_{\gamma}\cap L^{\infty}(X)$ such that
$$
e^{-\lambda \psi}(\gamma+i\partial\bar\partial \psi)^n=e^{-\lambda \varphi}(\gamma+i\partial\bar\partial \varphi)^n\,.
$$ 
with  $\lambda>0$. Then $\psi=\varphi$.
\end{theorem}
\emph{Proof of A}. By the $\partial\bar\partial$-lemma and by the previous statement \ref{C0-Kahl} we can assume that the current $\gamma$ is a smooth K\"ahler metric in the complement $X\smallsetminus Z$ of an analytic set. The identity $\gamma_{\varphi}^n=\gamma_{\psi}^n$ implies $\varphi\in \hat{\cal P}\BT^0_{\gamma}$ by Claim~\ref{StpdCv}. 
Let $\varphi_{\varepsilon}$, $\psi_{\varepsilon}$ be as in the statement of corollary \ref{ConvGradMA} and set $u:=\psi-\varphi$, $u_{\varepsilon}:=\psi_{\varepsilon}-\varphi_{\varepsilon}$. Let us also recall the formula
$$
\alpha^k-\beta^k=(\alpha-\beta)\wedge \sum_{l=0}^{k-1}\alpha^l\wedge \beta^{k-l-1}\,.
$$
From this we deduce
\begin{eqnarray}
0&=&
\int\limits_X -u(\gamma^n_{\psi}-\gamma^n_{\varphi})
=
\lim_{\varepsilon\rightarrow 0^+}\int\limits_X -u_{\varepsilon}(\gamma^n_{\psi}-\gamma^n_{\varphi})
\nonumber
\\\nonumber
\\
&=&
\lim_{\varepsilon\rightarrow 0^+}\,\sum_{l=0}^{n-1}\,\int\limits_X -u_{\varepsilon}\,i\partial\bar\partial u\wedge \gamma^l_{\psi}\wedge \gamma^{n-l-1}_{\varphi}\nonumber
\\\nonumber
\\
&=&
\lim_{\varepsilon\rightarrow 0^+}\,\sum_{l=0}^{n-1}\,\int\limits_X
i\partial u_{\varepsilon}\wedge\bar\partial u\wedge \gamma^l_{\psi}\wedge \gamma^{n-l-1}_{\varphi}\nonumber
\\\nonumber
\\
&=&\sum_{l=0}^{n-1}\,\int\limits_X i\partial u\wedge\bar\partial u\wedge \gamma^l_{\psi}\wedge \gamma^{n-l-1}_{\varphi}=:I\label{Uniq}\,,
\end{eqnarray}
since  $\partial u_{\varepsilon}\rightarrow \partial u$ in $L^2(X,\gamma^l_{\psi}\wedge \gamma^{n-l-1}_{\varphi})$  by corollary \ref{ConvGradMA}.
Inspired by an idea of S.~B{\l}ocki \cite{Blo1}, we will prove by induction on $k=0,...,n-1$ that
\begin{eqnarray}\label{UniqIndc}
\int\limits_X i\partial u\wedge\bar\partial u\wedge \gamma^r_{\psi}\wedge \gamma^s_{\varphi}\wedge \gamma^k=0
\end{eqnarray}
for all $r,s\geq 0$, $r+s=n-k-1$. For $k=0$ this follows from \eqref{Uniq}. So we assume \eqref{UniqIndc} for $k-1$ and we prove it for $k$. In fact consider the identity
$$
\gamma^k=\gamma^k_{\psi}-i\partial\bar\partial \psi \wedge \sum_{l=0}^{k-1}\gamma^l_{\psi}\wedge \gamma^{k-l-1}
\quad\mbox{and set }\quad
\Xi:=\gamma^r_{\psi}\wedge \gamma^s_{\varphi}\wedge \sum_{l=0}^{k-1}\gamma^l_{\psi}\wedge \gamma^{k-l-1}\,.
$$
By applying several times corollary \ref{ConvGradMA} and by integrating by parts we derive
\begin{eqnarray*}
&&\int\limits_X i\partial u\wedge\bar\partial u\wedge \gamma^r_{\psi}\wedge \gamma^s_{\varphi}\wedge \gamma^k
\;
=
\;
\lim_{\varepsilon\rightarrow 0^+}
\int\limits_X i\partial u_{\varepsilon}\wedge\bar\partial u\wedge \gamma^r_{\psi}\wedge \gamma^s_{\varphi}\wedge \gamma^k\nonumber
%\\\nonumber
%\\
%&=&
%\lim_{\varepsilon\rightarrow 0^+}\,
%\int\limits_X i\partial u_{\varepsilon}\wedge\bar\partial (u \gamma^r_{\psi}\wedge \gamma^s_{\varphi}\wedge \gamma^k)
%\nonumber
\\\nonumber
\\
&=&
\lim_{\varepsilon\rightarrow 0^+}\,\left[\;
\int\limits_X i\partial u_{\varepsilon}\wedge\bar\partial (u \gamma^{r+k}_{\psi}\wedge \gamma^s_{\varphi})
\;
-
\;
\int\limits_X i\partial u_{\varepsilon}\wedge\bar\partial (u\, i\partial\bar\partial \psi\wedge \Xi)\;\right]\nonumber
%\\\nonumber
\end{eqnarray*}
\begin{eqnarray}
&=&
\lim_{\varepsilon\rightarrow 0^+}\,\left[\;
\int\limits_X i\partial u_{\varepsilon}\wedge\bar\partial u\wedge \gamma^{r+k}_{\psi}\wedge \gamma^s_{\varphi}
\;
+
\;
\int\limits_X  u_{\varepsilon}\, i\partial\bar\partial u\wedge i\partial\bar\partial \psi\wedge \Xi\;\right]
\nonumber
\\\nonumber
\\
&=&
\int\limits_X i\partial u\wedge\bar\partial u\wedge \gamma^{r+k}_{\psi}\wedge \gamma^s_{\varphi}
\;
-
\;
\lim_{\varepsilon\rightarrow 0^+}\,
\int\limits_X   u_{\varepsilon}\, i\partial\bar\partial \psi\wedge (\gamma_{\varphi}-\gamma_{\psi}) \wedge\Xi\nonumber
\\\nonumber
\\
&\leq&
I+\lim_{\varepsilon\rightarrow 0^+}\,\int\limits_X i\partial u_{\varepsilon}\wedge\bar\partial \left[\,\psi\,(\gamma_{\varphi}-\gamma_{\psi}) \wedge \Xi\,\right]\nonumber
\\\nonumber
\\
%&=&
%\lim_{\varepsilon\rightarrow 0^+}\,\left[\;
%%\int\limits_X i\partial u_{\varepsilon}\wedge\bar\partial (\psi\, \gamma_{\varphi}\wedge \Xi)
%\;
%-
%\;
%\int\limits_X i\partial u_{\varepsilon}\wedge \bar\partial ( \psi\,\gamma_{\psi}\wedge %\Xi)\;\right]
%\nonumber
%\\\nonumber
%\\
&=&
\lim_{\varepsilon\rightarrow 0^+}\,\left[\;
\int\limits_X i\partial u_{\varepsilon}\wedge\bar\partial \psi\wedge \gamma_{\varphi}\wedge \Xi
\;
-
\;
\int\limits_X i\partial u_{\varepsilon}\wedge \bar\partial  \psi\wedge\gamma_{\psi}\wedge \Xi\;\right]
\nonumber
\\\nonumber
\\
&=&
\int\limits_X i\partial u\wedge\bar\partial \psi\wedge \gamma_{\varphi} \wedge \Xi
\;-\;
\int\limits_X i\partial u\wedge\bar\partial \psi\wedge \gamma_{\psi} \wedge \Xi\,.\label{UniII}
\end{eqnarray}
Set $\chi=\varphi$ or $\chi=\psi$. Then the Cauchy-Schwarz inequality implies
\begin{eqnarray*}
&&\left|\;\int\limits_X i\partial u\wedge\bar\partial \psi\wedge \gamma_{\chi} \wedge \Xi\;\right|
\\
\\
&\leq &
\left(\;\int\limits_X i\partial u\wedge\bar\partial u\wedge \gamma_{\chi} \wedge \Xi\;\right)^{1/2}
\left(\;\int\limits_X i\partial \psi\wedge\bar\partial \psi\wedge \gamma_{\chi} \wedge \Xi\;\right)^{1/2}=0\,,
\end{eqnarray*}
by the inductive hypothesis. This combined with \eqref{UniII} implies \eqref{UniqIndc} for $k$. So at the end of the induction we get $0=\int_Xi\partial u\wedge\bar\partial u\wedge \gamma^{n-1}$. This implies $\varphi=\psi$ by elementary properties of  plurisubharmonic functions.
\\
\\
%%%%%%%%%%%%%%%%%%%%%%%%%%%%%%%%%%%%%%%%%%%%%%%%%%%%%%%%%%%%%%%%%%%%%%%%%%%%%%%%%%%%%%%%%%%%%%%%%%%%%%%%%%%%%%%%%%%%%%%%%%%%%%%%%%%%%%%%%%%%%%%%%%%%%%%%%%%%%%%%%%%%%%%%%%%%%%%%%%%%%%%%%%%%%%%%%%%%%%%%%%%%%%%%%%%%%%%%%%%%%%%%%%%%%%%%%%%%%%%%%%%%%%%%%%%%%%%%%%%%%%%%%%%%%%%
\emph{Proof of B}.
By applying the comparison principle \eqref{comp-Prnc} as in \cite{E-G-Z1} we get
$$
\int\limits_{\varphi<\psi}\gamma^n_{\psi}\,\leq\, \int\limits_{\varphi<\psi}\gamma^n_{\varphi}
\,=\,
\int\limits_{\varphi<\psi}e^{\lambda(\varphi-\psi)}\gamma^n_{\psi}\,,
$$
which implies $\int_{\varphi<\psi}\gamma^n_{\psi}=0$ since $e^{\lambda(\varphi-\psi)}<1$. This implies that the inequality $\varphi\geq \psi$ holds $\gamma^n_{\psi}$-almost everywhere, thus the inequality 
$$
\gamma^n_{\varphi}=e^{\lambda(\varphi-\psi)}\gamma^n_{\psi}\geq \gamma^n_{\psi}\,,
$$ 
holds
$\gamma^n_{\psi}$-almost everywhere. By symmetry we also deduce that $\gamma^n_{\psi}\geq \gamma^n_{\varphi}$ holds $\gamma^n_{\varphi}$-almost everywhere. The fact that the potentials $\varphi$ and $\psi$ satisfies
\begin{eqnarray}\label{UnqNeg}
\gamma^n_{\varphi}=e^{\lambda(\varphi-\psi)}\gamma^n_{\psi}
\end{eqnarray}
implies that a property holds $\gamma^n_{\psi}$-almost everywhere if and only if it holds $\gamma^n_{\varphi}$-almost everywhere. We infer $\gamma^n_{\psi}=\gamma^n_{\varphi}$, 
hence $\psi-\varphi=\hbox{Const}$ by part (A), and equality \eqref{UnqNeg} now implies $\psi=\varphi$.
\hfill $\Box$
\\
\\
We now show a uniqueness result in the non-nef case. 
We denote by $\UB_{\varphi}\subset X$ the unbounded locus of a quasi-plurisubharmonic function $\varphi$. Let us first recall the following well known lemma \cite{Dem1}, \cite{Be-Bo}.
\begin{lemma}\label{fin-Mass}
Let $(X,\omega)$ be a compact K\"ahler manifold of complex dimension~$n$, $T$ a closed positive $(q,q)$-current on~$X$, $\theta$ 
a smooth closed real $(1,1)$-form
and $\varphi$ a quasi-plurisubharmonic function such that $\theta+i\partial\bar\partial \varphi\ge 0$ over~$X$. Then the following holds
\\
{\bf (A)} For all $k=1,...,n-q$ 
$$
\int\limits_{X\smallsetminus \UB_{\varphi}}(\theta+i\partial\bar\partial \varphi)^k\wedge T\wedge \omega^{n-k-q}\;<\;+\infty\,.
$$
\\
{\bf (B)} If in addition $\varphi$ has zero Lelong numbers, then  
$$
\int\limits_{X\smallsetminus \UB_{\varphi}}(\theta+i\partial\bar\partial \varphi)^k\wedge T\wedge \omega^{n-k-q}\;\le\; \int\limits_X\theta^k\wedge T\wedge \omega^{n-k-q}\,,
$$ 
for all $k=1,...,n-q$.
\end{lemma}
$Proof$. Set $\Theta:=\theta+i\partial\bar\partial \varphi\ge 0$. Let $(\varphi_{\varepsilon})_{\varepsilon>0}\subset C^{\infty}(X,\R)$  such that $\varphi_{\varepsilon}\downarrow \varphi$ as $\varepsilon\rightarrow 0$ and let $C>0$ be a sufficiently big constant such that $\Theta_{\varepsilon}:=\theta+i\partial\bar\partial \varphi_{\varepsilon}\ge -C\omega$ for all $\varepsilon\in (0,1)$. By the monotone decreasing theorem in pluripotential theory we infer that 
$$
(\Theta_{\varepsilon}+C\omega)^k\wedge T\longrightarrow (\Theta+C\omega)^k\wedge T\,,
$$
weakly over the open set $U:=X\smallsetminus \UB_{\varphi}$ as $\varepsilon\rightarrow 0$. We infer
\begin{eqnarray*}
\int\limits_U\Theta^k\wedge T\wedge \omega^{n-k-q}
&\le&
\int\limits_U(\Theta+C\omega)^k\wedge T\wedge \omega^{n-k-q}
\\
\\
&\le&
\liminf_{\varepsilon\rightarrow 0}\;\int\limits_U(\Theta_{\varepsilon}+C\omega)^k\wedge T\wedge \omega^{n-k-q}
\\
\\
&\le&
\int\limits_X(\theta+C\omega)^k\wedge T\wedge\omega^{n-k-q}\;<\;+\infty\,,
\end{eqnarray*}
which concludes the proof of statement (A). Statement (B) follows from the fact that, thanks to the work in \cite{Dem4}, we can replace the loss of positivity constant $C$ with constants $C_{\varepsilon}>0$ such that $C_{\varepsilon}\downarrow 0$  as $\varepsilon\rightarrow 0$.\hfill$\Box$
\\
\\
The following lemma can be found in \cite{Be-Bo} and is based on a simple but efficient increasing singularity approach  introduced by the first named author.
\begin{lemma}\label{ztmt-0}
Let $X$ be a compact K\"ahler manifold of complex dimension $n$, let $T$ be a closed positive $(n-q,n-q)$-current, $q\ge 1$, let
  $\theta$ be a smooth closed real $(1,1)$-form and consider
  $\varphi,\,\psi\in {\cal P}_{\theta}$ such that  $\varphi\ge \psi$ over $X$. Then
$$
\int\limits_{X\smallsetminus \UB_{\psi}}
(\theta+i\partial\bar\partial\varphi)^q\wedge T
\;\ge 
\int\limits_{X\smallsetminus \UB_{\psi}}
(\theta+i\partial\bar\partial\psi)^q\wedge T\,.
$$ 
\end{lemma}
$Proof$. Consider the closed positive current $\Theta:=(\theta+i\partial\bar\partial\psi)^{q-1}\wedge T$ over $X\smallsetminus \UB_{\psi}$. In order to conclude, it is sufficient to prove the inequality
\begin{eqnarray}\label{ztmt-1}
\int\limits_{X\smallsetminus \UB_{\psi}}
(\theta+i\partial\bar\partial\varphi)\wedge \Theta
\;\ge 
\int\limits_{X\smallsetminus \UB_{\psi}}
(\theta+i\partial\bar\partial\psi)\wedge\Theta\,,
\end{eqnarray}
thanks to the symmetry of the wedge product and to an obvious induction. 
(Notice that the integral on the left hand side of the inequality \eqref{ztmt-1} is also finite by the same type of argument as in the proof of Lemma \ref{fin-Mass}.)
Let $C>0$ be a sufficiently big constant such that $i\partial\bar\partial\psi\ge -C\omega$ and set $\psi_{\varepsilon}:=(1+\varepsilon)\psi\in {\cal P}_{\theta_{\varepsilon}}$, with $\theta_{\varepsilon}:=\theta+\varepsilon C\omega$. Then the inequality \eqref{ztmt-1} will follow by letting $\varepsilon\rightarrow 0$ in the inequality
\begin{eqnarray}\label{ztmt-2}
\int\limits_{X\smallsetminus \UB_{\psi}}
(\theta_{\varepsilon}+i\partial\bar\partial\varphi)\wedge \Theta
\;\ge 
\int\limits_{X\smallsetminus \UB_{\psi}}
(\theta_{\varepsilon}+i\partial\bar\partial\psi_{\varepsilon})\wedge\Theta\,,
\end{eqnarray}
that we prove now. Let  $k>0$ be an arbitrary constant. The fact that $\varphi-k>\psi_{\varepsilon}$ over the open set $\{\psi<-k/\varepsilon\}$ implies
$$
\int\limits_{X\smallsetminus \UB_{\psi}}
(\theta_{\varepsilon}+i\partial\bar\partial\varphi)\wedge \Theta
\;=
\int\limits_{X\smallsetminus \UB_{\psi}}
\left(\theta_{\varepsilon}+i\partial\bar\partial\max\{\varphi-k \,,\psi_{\varepsilon}\}\right)\wedge \Theta\;=:\;I\,,
$$
by Stokes' formula. Let $L\subset X\smallsetminus \UB_{\psi}$ be an arbitrary compact set and let $U\subset \overline{U}\subset X\smallsetminus \UB_{\psi}$  be an open set such that $L\subset U\subset \{\psi_{\varepsilon}\ge -R\}$ for a sufficiently big constant $R>0$. We infer
$$
I\;\ge\; \int\limits_U
\left(\theta_{\varepsilon}+i\partial\bar\partial\max\{\varphi-k \,, \psi_{\varepsilon}\}\right)\wedge \Theta
\;\ge\; 
\int\limits_L(\theta_{\varepsilon}+i\partial\bar\partial\psi_{\varepsilon})\wedge\Theta\,,
$$
for $k$ such that $\varphi-k<-R$ over $X$. Then the inequality \eqref{ztmt-2} follows by taking the supremum over $L$.\hfill$\Box$
\\
\\
We now define the potential with minimal singularities
$$
\varphi_{\theta}(x):=\sup\left\{\psi(x)\,\mid\,\psi\in {\cal P}_{\theta}^0\right\}\,,
$$ 
and we observe that $\varphi_{\theta}\in {\cal P}_{\theta}^0$. Let $\theta'\in \{\theta\}$ be another element in the cohomology class of $\theta$. We write $\theta=\theta'+i\partial\bar\partial u$. By definition, we infer $\varphi_{\theta}+u-C\le \varphi_{\theta'}$ and $\varphi_{\theta'}-u-C'\le \varphi_{\theta}$ for some constants $C,C'>0$. By Lemmas \ref{fin-Mass} and \ref{ztmt-0} we infer the equality
$$
\int\limits_{X\smallsetminus \UB_{\theta}}(\theta+i\partial\bar\partial
\varphi_{\theta})^q\wedge T
\;=\;
\max_{\psi\in {\cal P}_{\theta}}\int\limits_{X\smallsetminus \UB_{\psi}}(\theta+i\partial\bar\partial
\psi)^q\wedge T<+\infty\,.
$$
Notice that the closed set $\UB_{\theta}$ depends only on the cohomology class $\{\theta\}$. In the case $\UB_{\theta}$ is contained in a complete pluripolar set $E\subset X$, the trivial extension  of the current 
$$
\I_{X\smallsetminus E}\,(\theta+i\partial\bar\partial\varphi_{\theta})^q\,,
$$ 
over $X$ is closed and positive by the Skoda-El Mir extension theorem, which applies thanks to Lemma \ref{fin-Mass}. Moreover this extension is independent of $E$ by the definition of $\UB_{\theta}$. In fact the current 
$
(\theta+i\partial\bar\partial\varphi_{\theta})^q\,,
$
does not carry any mass on pluripolar sets contained in the open set $X\smallsetminus \UB_{\theta}$, since $\varphi_{\theta}$ is locally bounded over this set.
In this case we will still denote by $(\theta+i\partial\bar\partial\varphi_{\theta})^q$ the extension over $X$. In this setting we can define the cohomology invariant 
$$
\{\theta\}^{[q]}\cdot \{T\}:=\int\limits_{X}(\theta+i\partial\bar\partial
\varphi_{\theta})^q\wedge T
\;=\;
\max_{\psi\in {\cal P}_{\theta}}\int\limits_{X\smallsetminus \UB_{\psi}}(\theta+i\partial\bar\partial
\psi)^q\wedge T<+\infty\,.
$$
In general the number $\alpha^{[q]}\cdot \{\omega\}^{n-q}$ associated to a pseudoeffective class $\alpha\in H^{1,1}(X,\R)$ over a compact K\"{a}hler manifold $(X,\omega)$ of complex dimension $n$ is not a cohomology invariant, so we will denote it by $\alpha^{[q]}\cdot\omega^{n-q}$. However the numerical dimension of $\alpha$, namely the number
$$
\nu(\alpha):=\max\{q\in \{0,...,n\}\,\mid\,\alpha^{[q]}\cdot \omega^{n-q}>0\,\}\,.
$$
is well defined. In fact it is independent of the choice of the K\"{a}hler metric $\omega$ since the trace operator controls the mass of a positive  $(q,q)$-current.
We prove now the following degenerate version of the Comparison Principle which is also based on the increasing singularity approach previously used. 
(Compare with the statement and the proof of corollary 1.4 in \cite{Be-Bo}).
\begin{lemma}{\bf (Degenerate Comparison Principle).}\label{Cmp-Prnc}
Let $X$ be a compact K\"ahler manifold of complex dimension $n$, let
$\theta$ be a smooth closed real $(1,1)$-form and consider
$\varphi,\,\psi\in {\cal P}_{\theta}$ such that 
$\varphi\ge \psi-K$ for some constant $K>0$. Then 
$$
\int\limits_{\{\varphi<\psi\}\smallsetminus \UB_{\psi}} (\theta+i\partial\bar\partial \psi)^n
\;\leq \int\limits_{\{\varphi<\psi\}\smallsetminus \UB_{\psi}}(\theta+i\partial\bar\partial
\varphi)^n\,.
$$
\end{lemma}
\emph{Proof.} For any set $E\subset X$ we put $E_{\psi}:=E\smallsetminus \UB_{\psi}$ and define the closed positive current $\Theta:=\left(\theta+i\partial\bar\partial \max\{\varphi,\psi\}\right)^n$ over $X_{\psi}$.
We start by proving the inequality
\begin{eqnarray}\label{Cmp-P1}
\int\limits_{X_{\psi}}(\theta+i\partial\bar\partial
\varphi)^n 
\;\ge\; 
\int\limits_{X_{\psi}}\Theta\,.
\end{eqnarray}
Let $R>0$ be a sufficiently big constant such that $i\partial\bar\partial\psi\ge -R\omega$ and set $\psi_{\varepsilon}:=(1+\varepsilon)\psi\in {\cal P}_{\theta_{\varepsilon}}$, with $\theta_{\varepsilon}:=\theta+\varepsilon R\omega$. The fact that $\varphi>\psi_{\varepsilon}$ on the open set $\{\psi<-K/\varepsilon\}$ implies
\begin{eqnarray*}
\int\limits_{X_{\psi}}(\theta_{\varepsilon}+i\partial\bar\partial
\varphi)^n 
\;=\; 
\int\limits_{X_{\psi}}\left(\theta_{\varepsilon}+i\partial\bar\partial \max\{\varphi,\psi_{\varepsilon}\}\right)^n\,,
\end{eqnarray*}
by Stokes' formula. We infer
\begin{eqnarray*}%\label{Cmp-P1}
\int\limits_{X_{\psi}}(\theta+i\partial\bar\partial
\varphi)^n 
\;=\;
\liminf_{\varepsilon\rightarrow 0}\int\limits_{X_{\psi}}\left(\theta_{\varepsilon}+i\partial\bar\partial \max\{\varphi,\psi_{\varepsilon}\}\right)^n
\;\ge\;
\int\limits_{X_{\psi}}\Theta\,.
\end{eqnarray*}
by the weak convergence 
\begin{eqnarray}\label{Cv-Cmp}
\Theta_{\varepsilon}:=\left(\theta_{\varepsilon}+i\partial\bar\partial \max\{\varphi,\psi_{\varepsilon}\}\right)^n\longrightarrow \Theta\,,
\end{eqnarray}
as $\varepsilon\rightarrow 0$ over the open set $X_{\psi}$. In order to prove the convergence \eqref{Cv-Cmp} we restrict our considerations to an arbitrary open set $U\subset \overline{U}\subset X_{\psi}$. Let $C>0$ be a constant such that $\psi\ge -C$ over $U$. Then the function
$$
\Phi_{\varepsilon}:=\max\{\varphi+\varepsilon C,\psi_{\varepsilon}+\varepsilon C\}\in {\cal P}_{\theta_{\varepsilon}}\,,
$$
decreases to $\max\{\varphi,\psi\}$ over $U$ as $\varepsilon\rightarrow 0$ and satisfies $\Theta_{\varepsilon}=\left(\theta_{\varepsilon}+i\partial\bar\partial\Phi_{\varepsilon} \right)^n$. Then the convergence \eqref{Cv-Cmp} over $U$ follows from the monotone decreasing theorem in pluripotential theory. 
On the other hand the inequality of measures 
$$
\Theta\;\ge\; \I_{_{\varphi\ge \psi}}\,\theta_{\varphi}^n\;+\;\I_{_{\varphi < \psi}}\,\theta_{\psi}^n\,,
$$
over the open set $X_{\psi}$ (see \cite{Dem1}), implies
\begin{eqnarray*}
\int\limits_{\{\varphi<\psi\}_{\psi}}\Theta
\;\ge\;
\int\limits_{\{\varphi<\psi\}_{\psi}}\theta^n_{\psi}\;,
\qquad\qquad
\int\limits_{\{\varphi\ge\psi\}_{\psi}}\Theta
\;\ge\;
\int\limits_{\{\varphi\ge\psi\}_{\psi}}\theta^n_{\varphi}\,.
\end{eqnarray*}
This combined with the inequality \eqref{Cmp-P1} implies
\begin{eqnarray*}
\int\limits_{\{\varphi<\psi\}_{\psi}}\theta_{\psi}^n
\;\le\;
\int\limits_{X_{\psi}}\Theta
\;-
\int\limits_{\{\varphi\ge\psi\}_{\psi}}\Theta
\;\le\;
\int\limits_{X_{\psi}}\theta_{\varphi}^n 
\;-
\int\limits_{\{\varphi\ge\psi\}_{\psi}}  \theta_{\varphi}^n
\;=
\int\limits_{\{\varphi<\psi\}_{\psi}}\theta_{\varphi}^n\,.
\end{eqnarray*}
\hfill $\Box$
\begin{corollary}\label{cmp-KE}
Let $X$ be a compact K\"ahler manifold of complex dimension $n$, $\Omega>0$ a smooth volume form and $\theta$ a smooth closed real $(1,1)$-form. Assume that $\varphi_j\in {\cal P}_{\theta}$, $j=1,2$ is such that $\UB_{\varphi_j}$ is a zero measure set and
$$
(\theta+i\partial\bar\partial
\varphi_j)^n=e^{\varphi_j}\Omega
$$
over $X\smallsetminus \UB_{\varphi_j}$.
If $\varphi_1\ge \varphi_2-K$ over $X$, for some constant $K>0$ then $\varphi_1\ge \varphi_2$ over $X$. In particular if there exists a K\"ahler-Einstein current  
$\omega_{_{E}}\in 2\pi c_1(K_X)$, then this current is unique in the class of currents with minimal singularities in $2\pi c_1(K_X)$.
\end{corollary}

{\it Proof}. We set $E:=\{\varphi_1<\varphi_2\}\smallsetminus \UB_{\varphi_2}$ and 
we apply the degenerate comparison principle \ref{Cmp-Prnc} as before. We obtain
\begin{eqnarray*}
\int\limits_{E}(\theta+i\partial\bar\partial\varphi_2)^n
\;\le\;
\int\limits_{E}(\theta+i\partial\bar\partial\varphi_1)^n
\;=\;
\int\limits_{E}e^{\varphi_1-\varphi_2}(\theta+i\partial\bar\partial\varphi_2)^n\,.
\end{eqnarray*}
We infer $0=\int_Ee^{\varphi_2}\Omega$, and so $\varphi_1\ge \varphi_2$ almost everywhere over $X$, thus everywhere by elementary properties of quasi-plurisubharmonic functions.\hfill$\Box$
%%%%%%%%%%%%%%%%%%%%%%%%%%%%%%%%%%%%%%%%%%%%%%%%%%%%%%%%%%%%%%%%%%%%%%%%%%%%%%%%%%%%%%%%%%%%%%%%%%%%%%%%%%%%%%%%%%%%%%%%%%%%%%%%%%%%%%%%%%%%%%%%%%%%%%%%%%%%%%%%%%%%%%%%%%%%%%%%%%%%%%%%%%%%%%%%%%%%%%%%%%%%%%%%%%%%%%%%%%%%
%%%%%%%%%%%%%%%%%%%%%%%%%%%%%%
%%%%%%%%%%%%%%%%%%%%%%%%%%%%%%%%%%%%%%%%%%%%%%%%%%%%%%%%%%%%%%%%%%%%%%%%%%%%%%%%%%%%%%%%%%%%%%%%%%%%%%%%%%%%%%%%%%%%%%%%%%%%%%%%%%%%%%%%%%%%%%%%%%%%%%%%%%%%%%%%%%%%%%%%%%%%%%%%%%%%%%%%%%%%%%%%%%%%%%%%%%%%%%%%%%%%%%%%%%%%%%%%%%%%%%%%%%%%%%%%%%%%%%%%%%%%%%%%%%%%%%%%%%%%%%%%%%%%%%%%%%%%%%%%%%%%%%%%%%%%%%%%
\section{Generalized Kodaira lemma}\label{KOAIRLEMM}
%%%%%%%%%%%%%%%%%%%%%%%%%%%%%%%%%%%%%%%%%%%%%%%%%%%%%%%%%%%%%%%%%%%%%%%%%%%%%%%%%%%%%%%%%%%%%%%%%%%%%%%%%%%%%%%%%%%%%%%%%%%%%%%%%%%%%%%%%%%%%%%%%%%%%%%%%%%%%%%%%%%%%%%%%%%%%%%%%%%%%%%%%%%%%%%%%%%%%%%%%%%%%%%%%%%%%%%%%%%%%%%%%%%%%%%%%%%%%%%%%%%%%%%%%%%%%%%%%%%%%%%%%%%%%%%%%%%%%%%%%%%%%%%%%%%%%%%%%%%%%%%%
We first recall a few standard definitions of algebraic and analytic geometry which will be useful in our situation. 
\begin{definition}
Let $(X,\omega)$ be a compact K\"{a}hler manifold. 
\\
{\bf (A)} A modification of $X$ is  a bimeromorphic morphism $\mu:\tilde X\rightarrow X$ of compact complex manifolds with connected fibers. Then there is a
smallest analytic set $Z\subset X$ such that the restriction $\mu:\tilde X\smallsetminus\mu^{-1}(Z)\to X\smallsetminus Z$ is a biholomorphism; we say that
$\Exc(\mu)=\mu^{-1}(Z)$ is the exceptional locus of $\mu$.
\\
{\bf (B)}
A class $\chi\in H^{1,1}(X,\R)$ is called big if there exist a current $T\in \chi$ such that $T\geq \varepsilon\omega$, for some $\varepsilon>0$.
\end{definition}
By a result of \cite{De-Pa}, a nef class $\chi$ on a compact K\"ahler manifold is big if and only if $\int_X\chi^n>0$. By the proof of Theorem 3.4 in \cite{De-Pa} we obtain the following generalization of Kodaira's lemma.
\begin{lemma}\label{Kod-Lm}
Let $X$ be a compact K\"{a}hler manifold and  $\chi\in H^{1,1}(X,\R)$ be a big class. Then there exist a modification $\mu:\tilde X\rightarrow X$ of $X$, an effective integral divisor $D$ on $\tilde X$ with support $|D|\supset \Exc(\mu)$ and a number $\delta\in \Q_{>0}$, such that the class $\mu^*\chi-\delta\{D\}$ is K\"{a}hler.

\end{lemma}
We associate to $\chi$ the set $I_{\chi}$ of triples $(\mu,D,\delta)$ satisfying the generalized Kodaira lemma \ref{Kod-Lm}, and a complex analytic set $\Sigma_\chi$ which we call the \emph{augmented singular locus} of $\chi$, defined as
\begin{eqnarray}\label{Sigma-chi}
\Sigma_{\chi}:=\bigcap_{(\mu,D,\delta)\in I_{\chi}}\mu(|D|)\,.
\end{eqnarray}
A trivial approximation argument shows that the set $\Sigma_{\chi}$ depends only on the half line $\R_{>0}\chi$. In the case the class $\chi$ is K\"{a}hler, $(\id_X, 0, 1)\in I_{\chi}$, thus~$\Sigma_{\chi}=\emptyset$. Conversely, if $\Sigma_\chi=\emptyset$, it is clear that the class $\chi$ must be K\"ahler$\,$: in fact, if $\tilde\omega_{\mu,D,\delta}$ is a K\"ahler metric in $\mu^*\chi-\delta\{D\}$, then \hbox{$\Theta=\mu_*(\tilde\omega_{\mu,D,\delta}+\delta[D])$} is a K\"ahler current contained in the class $\chi$, which is smooth on $X\smallsetminus\mu(|D|)$ and possesses logarithmic poles on $\mu(|D|)\,$; by taking the regularized upper envelope of a finite family of potentials of such currents $\Theta_j$ with \hbox{$\bigcap\mu(|D_j|)=\emptyset$}, we obtain a (smooth) K\"ahler metric on $X$. In the case the class $\chi$ is integral or rational, the set $\Sigma_{\chi}$ can be characterized as follows. 

\begin{lemma}\label{SBlocsLm}
Let $L$ be a big line bundle over a compact K\"ahler manifold. Then the class $\chi:=c_1(L)$ satisfies
\begin{eqnarray}\label{SBlocs}
SB(L)\subset \Sigma_{\chi}
=\bigcap_{E\in\Div^+(X),\;\delta\in\Q_{>0},
\atop
\chi-\delta\{E\}\,{\rm ample}}|E|\,,
\end{eqnarray}
where $SB(L)$ is the stable base locus of $L$, i.e.\ the intersection of the base loci of all line bundles $kL$, and $E$ runs over all effective integral divisors of~$X$.
\end{lemma}
\emph{Proof}. First notice that the existence of a big line bundle implies that $X$ is Moishezon. This combined with the assumption that $X$ is K\"ahler shows that $X$ must in fact be projective (see \cite{Moi}, and also \cite{Pet1}, \cite{Pet2} for a simple proof). The inclusion $SB(L)\subset \Sigma_{\chi}$ in \eqref{SBlocs} is quite easy: Let $(\mu, D, \rho)\in I_{\chi}$. Then Kodaira's theorem implies that $\{\alpha\}:=\mu^* \chi-\rho\{D\}$ is a $\Q$-ample class on $\tilde X$ and so the integer multiples $k\alpha$ are base point free for $k$ large enough. Therefore the base locus of $k\mu^*L$ is contained in $|D|$. This shows that $SB(L)$ is contained in the intersection of the sets $\mu(|D|)$, which is precisely equal to $\Sigma_\chi$ by definition. Now, if $H$ is an ample divisor on $X$, we have
$$
\mu^*(\chi-\varepsilon\{H\})=\rho\{D\}+\{\alpha\}-\varepsilon\{\mu^*H\}
$$
and, again, $\alpha-\varepsilon\mu^*H$ is ample for $\varepsilon\in\Q_{>0}$ small. We infer that the base locus of $k(L-\varepsilon H)$ is contained in~$\Sigma_\chi$ for $k$ large and sufficiently divisible. If we pick any divisor $E$ in the linear system of $k(L-\varepsilon H)$, then $L-\frac{1}{k}E\equiv \varepsilon H$
is an ample class, and the intersection of all these divisors $E$ is contained in~$\Sigma_\chi$. Therefore
$$
\bigcap_{
E\in\Div^+(X),\;\delta\in\Q_{>0},
\atop
\chi-\delta\{E\}\,{\rm ample}}|E|
\subset \Sigma_{\chi}.
$$
The opposite inclusion is obvious.\hfill$\Box$
\\
\\
The following lemma gives us an important class of densities which will be allowable as the right hand side of degenerate complex Monge-Amp\`ere equations.

\begin{lemma}\label{IntegLp}
Let $X$ be a compact complex manifold, let $\Omega>0$ be a smooth volume form and let $\sigma_j\in H^0(X,E_j)$,  $\tau_r\in H^0(X,F_r)$, $j=1,...,N$, $r=1,...,M$ be, non identically zero, holomorphic sections of some holomorphic vector bundles over $X$ such that the integral condition
\begin{eqnarray*}
\int\limits_X\,\prod\limits_{j=1}^N|\sigma_j|^{2l_j}\cdot \prod\limits_{r=1}^M|\tau_r|^{-2h_r}\,\Omega<+\infty
\end{eqnarray*}
holds for some real numbers $l_j\geq 0, \,h_r\geq 0$. Then the integrand function belongs to some $L^p$ space, $p>1$, and for $A\geq A_0\geq 0$ large enough, the family of functions 
\begin{eqnarray*}
G_{\varepsilon}:=\prod\limits_{j=1}^N(|\sigma_j|^2+\varepsilon^A)^{l_j}\cdot \prod\limits_{r=1}^M(|\tau_r|^2+\varepsilon)^{-h_r}\,,\qquad \varepsilon\in [0,1)
\end{eqnarray*}
converges in $L^p$-norm to the function $G_0$ when 
$\varepsilon\rightarrow 0$. In fact, for $N\neq 0$ and $l_j>0$, one can take
$A_0:=(\sum_r h_r)/(\min_j l_j)$.
\end{lemma}

\noindent
\emph{Proof}. By blowing-up the coherent ideals generated by the components
of any of the sections $\sigma_j$, $\tau_r$, we obtain a modification
$\mu: \tilde X\to X$ such that the pull-back of these ideals by $\mu$
is a divisorial ideal. Using Hironaka's desingularization theorem, we can
even assume that all divisors obtained in this way form a family
of normal crossing divisors in $\tilde X$. Then each square 
$|\sigma_j\circ\mu|^2$ (resp.\ $|\tau_r\circ\mu|^2$)  can be written as the
square $|z^\alpha|^2$ (resp.\ $|z^\beta|^2$) of a monomial in suitable 
local coordinates $U$ on a neighborhood of any point of~$\tilde X$,
up to invertible factors. The Jacobian of $\mu$ can also be assumed to 
be equal to
a monomial $z^\gamma$, up to an invertible factor. In restriction
to such a neighborhood~$U$, the convergence of the integral is equivalent 
to that of
$$
\int\limits_U|z^\gamma|^2
\prod_{j=1}^N|z^{\alpha_j}|^{2l_j}\prod_{r=1}^M|z^{\beta_r}|^{-2h_r}\,dz.
$$
Notice also that $\tilde X$ can be covered by finitely many such 
neighborhoods, by compactness. Now it is clear that if the integral 
is convergent, then the integrand must be in some $L^p$, $p>1$, because 
the integrability
condition is precisely that each coordinate $z_j$ appears with an 
exponent${}>-1$ in the $n$-tuple 
$$
\gamma+\sum l_j\alpha_j-\sum h_r\beta_r\,,
$$
(so that we can still replace $l_j$, $h_r$ with $pl_j$, $ph_r$ with
$p$ close to $1$). In order to prove the convergence of the functions $G_{\varepsilon}$ in the $L^p$ norm  we distinguish two cases. In the case where $l_j=0$ for all $j$, the claim follows immediately from the  monotone convergence theorem. The other possible case is $l_j>0$ for all $j$. In this case the convergence statement 
will follow if we can prove that for $A$ large enough the functions
$$
|z^\gamma|^2\prod_{j=1}^N(|z^{\alpha_j}|^2+\varepsilon^A)^{l_j}
\prod_{r=1}^M(|z^{\beta_r}|^2+\varepsilon)^{-h_r}$$
converge in $L^p$-norm as $\varepsilon\to 0$. This is trivial my
monotonicity when $N=0$. When $N>0$ and $l_j>0$, we have
$$
\prod_{j=1}^N(|z^{\alpha_j}|^2+\varepsilon^A)^{l_j}\le 
C\Big(\prod_{j=1}^N(|z^{\alpha_j}|^{2l_j}+\varepsilon^{A\min l_j}\Big),
\quad
\prod_{r=1}^M(|z^{\beta_r}|^2+\varepsilon)^{-h_r}\le \varepsilon^{-\sum h_r},
$$
so it is sufficient to take $A\ge(\sum h_r)/(\min l_j)$ to obtain the
desired uniform $L^p$-integrability in~$\varepsilon$.\hfill$\Box$

%%%%%%%%%%%%%%%%%%%%%%%%%%%%%%%%%%%%%%%%%%%%%%%%%%%%%%%%%%%%%%%%%%%%%%%%%%%%%%%%%%%%%%%%%%%%%%%%%%%%%%%%%%%%%%%%%%%%%%%%%%%%%%%%%%%%%%%%%%%%%%%%%%%%%%%%%%%%%%%%%%%%%%%%%%%%%%%%%%%%%%%%%%%%%%%%%%%%%%%%%%
%%%%%%%%%%%%%%%%%%%%%%%%%%%%%%%%%%%%%%%%%%%%%%%%%%%%%%%%%%%%%%%%%%%%%%%%%%%%%%%%%%%%%%%%%%%%%%%%%%%%%%%%%%%%%%%%%%%%%%%%%%%%%%%%%%%%%%%%%%%%%%%%%%%%%%%%%%%%%%%%%%%%%%%%%%%%%%%%%%%%%%%%%%%%%%%%%%%%%%%%%%
\section{Existence and higher order regularity of solutions}
%%%%%%%%%%%%%%%%%%%%%%%%%%%%%%%%%%%%%%%%%%%%%%%%%%%%%%%%%%%%%%%%%%%%%%%%%%%%%%%%%%%%%%%%%%%%%%%%%%%%%%%%%%%%%%%%%%%%%%%%%%%%%%%%%%%%%%%%%%%%%%%%%%%%%%%%%%%%%%%%%%%%%%%%%%%%%%%%%%%%%%%%%%%%%%%%%%%%%%%%%%
We are ready to prove the following fundamental existence theorem.
\begin{theorem}\label{HgRegMA}
Let $X$ be a compact connected K\"ahler manifold of complex dimension \hbox{$n\geq 2$}, let $\omega\geq 0$ be a big closed smooth $(1,1)$-form and let $\Omega>0$ be a smooth volume form. Consider also $\sigma_j\in H^0(X,E_j)$,  $\tau_r\in H^0(X,F_r)$, $j=1,...,N$, $r=1,...,M$ be non identically zero holomorphic sections of some holomorphic vector bundles over $X$, such that the integral condition
\begin{eqnarray}\label{IntgCond}
\int\limits_X\,\prod\limits_{j=1}^N|\sigma_j|^{2l_j}\cdot \prod\limits_{r=1}^M|\tau_r|^{-2h_r}\,\Omega= \int\limits_X\omega^n
\end{eqnarray}
holds for certain real numbers $l_j\geq 0, \,h_r\geq 0$. 
Then there exists a unique solution $\varphi\in {\cal P}\BT_{\omega}$ of the degenerate complex Monge-Amp\`ere equation
\begin{eqnarray}\label{MainMA}
(\omega+i\partial\bar\partial \varphi)^n= \prod\limits_{j=1}^N|\sigma_j|^{2l_j}\cdot \prod\limits_{r=1}^M|\tau_r|^{-2h_r}\,e^{\lambda\varphi}\,\Omega\,,\quad \lambda\geq 0\,,
\end{eqnarray}
which in the case $\lambda=0$ is normalized by $\sup_X\varphi=0$.
Moreover let $\Sigma_{\{\omega\}}$ be the augmented singular locus of the $(1,1)$-cohomology class $\{\omega\}$ as defined in~\eqref{Sigma-chi}, which is empty if the class $\{\omega\}$ is K\"ahler, and consider the complex analytic sets
$$
S':=\Sigma_{\{\omega\}}\cup \left(\bigcup_r\{\tau_r=0\}\right) \,,\quad S:=S'\cup\left(\bigcup_j\{\sigma_j=0\}\right).
$$
Then 
$\varphi\in {\cal P}_{\omega}\cap L^{\infty}(X)\cap C^0(X\smallsetminus \Sigma_{\{\omega\}})\cap C^{1,1}(X\smallsetminus S')\cap C^{\infty}(X\smallsetminus S)\,$.
\end{theorem}
%A very interesting particular case of the previous result is given by the equation $(\omega+i\partial\bar\partial \varphi)^n=|\kappa|^2\Omega$, with $\kappa$ a meromorphic section of the canonical bundle.
\emph{Proof}. We first assume the existence of an effective divisor $D$ in $X$ and of a small number \hbox{$\delta>0$} such that $\{\omega\}-\delta\{D\}$ is a K\"{a}hler class on~$X$. (We will later be able to remove this assumption thanks to Lemma \ref{Kod-Lm}). By using the Lelong-Poincar\'{e} formula we infer the existence of a smooth hermitian metric on ${\cal O}(D)$ such that
$$
0<\omega_{\delta}:=\omega- 2\pi\delta[D]+\delta\, i\partial\bar\partial\log|s|^2\,,
$$
with $\div(s)=D$. By  convention we will put $\delta=0$ if $\omega$ is a K\"ahler metric (so that $\omega_\delta=\omega$ in that case), and in general we will denote by $|D|$ the support of the divisor $D$.
\\
\\
{\bf (A) Setup.}
\\
For the sake of simplicity of notation we assume $N=M=1$. The general case would be entirely similar and we leave it to the reader. Let $\alpha>0$ be a K\"ahler metric, let $\varepsilon\in (0,1)$ and let $c_{\varepsilon}$ be a normalizing constant for the integral condition
\begin{eqnarray}\label{IntCondEps}
e^{c_{\varepsilon}}\,\int\limits_X\frac{(|\sigma|^2+\varepsilon^A)^l}{(|\tau|^2+\varepsilon)^h}\;\Omega=\int\limits_X(\omega+\varepsilon\alpha)^n\,,
\end{eqnarray}
with $A:=h/l$. %(We apply the same convention as in Lemma \ref{IntegLp}.)
Condition \eqref{IntgCond} combined with Lemma  \ref{IntegLp} implies  $c_{\varepsilon}\rightarrow 0$, when $\varepsilon\rightarrow 0^+$. Observe that here $\omega+\varepsilon\alpha$ is a K\"ahler metric for every $\varepsilon>0$. Consider the standard solutions $\varphi_{\varepsilon}\in C^{\infty}(X)$ of the complex Monge-Amp\`ere equations
\begin{eqnarray}\label{MAsingEps}
(\omega+\varepsilon\alpha+i\partial\bar\partial \varphi_{\varepsilon})^n=e^{c_{\varepsilon}}\,\frac{(|\sigma|^2+\varepsilon^A)^l}{(|\tau|^2+\varepsilon)^h}\;e^{\lambda\varphi_{\varepsilon}}\,\Omega\,,
\end{eqnarray}
given by the Aubin-Yau solution of the Calabi conjecture. As usual, in the case $\lambda=0$, we normalize the solution  $\varphi_{\varepsilon}$ with the condition $\max_X\varphi_{\varepsilon}=0$.
Notice that the integral condition \eqref{IntCondEps} implies that a non identically zero solution $\varphi_{\varepsilon}$ changes signs in the case $\lambda>0$. By combining Lemma \ref{IntegLp} with the estimate of corollary  \ref{C0est} we obtain a uniform bound for the oscillations,
$\Osc(\varphi_{\varepsilon})\leq C$. Set now 
$$
\tilde\omega_{\varepsilon}:=\omega_{\delta}+\varepsilon\alpha\,,
$$ 
and 
$$
\psi_{\varepsilon}:=\varphi_{\varepsilon}-\delta\log|s|^2\,.
$$ 
As the notation indicates, we will keep $\delta$ fixed (until step (E)).
Then we get a K\"ahler metric defined over $X$
\begin{eqnarray}\label{SmoothOVERdivisor}
\hat\omega_\varepsilon:=\omega+\varepsilon\alpha+i\partial\bar\partial \varphi_{\varepsilon}= \tilde\omega_{\varepsilon}+i\partial\bar\partial\psi_{\varepsilon}+2\pi\delta[D],
\end{eqnarray}
thus $\hat\omega_\varepsilon=\tilde\omega_{\varepsilon}+i\partial\bar\partial\psi_{\varepsilon}$ over $X\smallsetminus |D|$. In this setting, equation \eqref{MAsingEps} can be rewritten as
\begin{eqnarray}\label{MAsingEps2}
(\tilde\omega_{\varepsilon}+i\partial\bar\partial\psi_{\varepsilon})^n
=e^{F^{\varepsilon}+\lambda\delta\log|s|^2+\lambda\psi_{\varepsilon}}\,\tilde\omega_{\varepsilon}^n
\end{eqnarray}
on $X\smallsetminus |D|$, with
$
F^{\varepsilon}:=f^{\varepsilon}+l\cdot a^{\varepsilon}-h\cdot b^{\varepsilon}
$, 
and with 
$$
f^{\varepsilon}:=c_{\varepsilon}+\log(\Omega/\tilde\omega_{\varepsilon}^n)\,,\quad a^{\varepsilon}:=\log(|\sigma|^2+\varepsilon^A)\,, \quad
b^{\varepsilon}:=\log(|\tau|^2+\varepsilon)\,.
$$
(Here the superscripts in $\varepsilon$ are indices and not powers.) 
Let ${\cal C}_{\tilde\omega_{\varepsilon}}$ be the Chern curvature form  of the K\"{a}hler metric $\tilde\omega_{\varepsilon}>0$ and let
%Consider now the function $\gamma_1^{\delta, \varepsilon}:X\rightarrow \R$ defined by the formula
$$
\gamma_{\varepsilon}:=\min_{x\in X}\,\min_{\xi,\eta \in T_{X,x}\smallsetminus 0_x} 
{\cal C}_{\tilde\omega_{\varepsilon}}(\xi\otimes \eta ,\xi\otimes \eta)|\xi|_{\tilde\omega_{\varepsilon}} ^{-2}|\eta|_{\tilde\omega_{\varepsilon}} ^{-2}\,. 
$$
%So $\gamma_1^{\delta, \varepsilon}(x)$ is the smallest eigenvalue of the Chern curvature form  of the metric $\tilde\omega_{\varepsilon}>0$. 
(We remark that the minimum is always achieved by an easy compactness argument, see e.g.\ \cite{Kat}, Chap II, Sect.\ 5.1, Theorem 5.1, page 107.) 
%that the function  $\gamma_1^{\delta, \varepsilon}$ is continuous. 
We observe that the family of metrics $(\tilde\omega_{\varepsilon})_{\varepsilon}$ has bounded geometry for $\delta$ fixed and $\varepsilon\in [0,1]$ arbitrary. In particular, for all $\varepsilon\in [0,1]$
$$
\gamma_{\varepsilon}\geq \Gamma\,,\quad |f^{\varepsilon}|\leq K_{0}\,,
\quad 
\lambda(\omega_{\delta}-\omega)+
i\partial\bar\partial f^{\varepsilon}\geq -K_{0}\,\tilde\omega_{\varepsilon}\,.
$$
Moreover we can assume $i\partial\bar\partial a^{\varepsilon}\,,\,i\partial\bar\partial b^{\varepsilon} \geq -K_{0}\,\tilde\omega_{\varepsilon}$ , (see Appendix A.)
\\
\\
{\bf (B) The Laplacian estimate.}
\\
This estimate is obtained as a combination of ideas of Yau, B{\l}ocki and Tsuji, \cite{Yau}, \cite{Blo2}, \cite{Ts}.
Consider the continuous function  $\Lambda_{\varepsilon}:X\rightarrow (0,+\infty)$ given by the maximal eigenvalue of $\hat\omega_\varepsilon=\tilde\omega_{\varepsilon}+i\partial\bar\partial\psi_{\varepsilon}$ with respect to the K\"{a}hler metric~$\tilde\omega_{\varepsilon}$,
$$
\Lambda_{\varepsilon} (x):=\max_{\xi \in T_{X,x}\smallsetminus 0_x } (\tilde\omega_{\varepsilon}+i\partial\bar\partial\psi_{\varepsilon})(\xi ,J\xi)|\xi|_{\tilde\omega_{\varepsilon}} ^{-2}\,, 
$$
i.e.\ we extend $\Lambda_{\varepsilon}$ over $|D|$ by continuity, as is permitted by \eqref{SmoothOVERdivisor}. 
Consider also the continuous function over $X\smallsetminus |D|$,
$$
A_{\varepsilon}:=\log \Lambda_{\varepsilon}-k\cdot\psi_{\varepsilon}+h\cdot b^{\varepsilon}\,,
$$
with $0<k:=2(1+h\,K_{0}/2-K_1)$  and
$$
K_1:=\min\{-[\lambda+(1+ l)K_{0}/(2n)]\,,\,\Gamma\}<-\lambda\,.
$$
The reason for this crucial choice will be clear in a moment.
The singularity of the function $\psi_{\varepsilon}$ implies the existence of a maximum of the function $A_{\varepsilon}$ at a certain point $x_{\varepsilon}\in X\smallsetminus |D|$. Let $g_{\varepsilon}$ be a smooth real valued function in a neighborhood of $x_{\varepsilon}$ in $X\smallsetminus |D|$ such that $\tilde\omega_{\varepsilon}=\frac{i}{2}\partial\bar\partial g_{\varepsilon}$, and let $u_{\varepsilon}:=g_{\varepsilon}+2\psi_{\varepsilon}$. Then
$$
\tilde\omega_{\varepsilon}+i\partial\bar\partial\psi_{\varepsilon}=\frac{i}{2}\partial\bar\partial u_{\varepsilon}\,.
$$
For the simplicity of notation, we just put $g=g_{\varepsilon}$ and $u=u_{\varepsilon}$ from now on, and we also set
$u_{l,\bar{m}}:=\smash{\frac{\partial^2 u}{\partial z_l\partial\bar{z}_m}}$.
Let $(z_1,\ldots,z_n)$ be $\tilde\omega_{\varepsilon}$-geodesic holomorphic coordinates centered at the point $x_{\varepsilon}$, such that the metric \hbox{$\hat\omega_\varepsilon=\tilde\omega_{\varepsilon}+i\partial\bar\partial\psi_{\varepsilon}$} can be written in diagonal form in $x_{\varepsilon}$.  Explicitly, we have the local expression
$\tilde\omega_{\varepsilon} =\frac{i}{2}\sum_{l,m}g_{l,\bar{m}}\,dz_l\wedge d\bar{z}_m$, with
\begin{eqnarray*}
&\displaystyle{
g_{l,\bar{m}}=\delta _{l,m}-\sum_{j,k}C^{j,\bar{k}}_{l, \bar m}z_j\bar{z}_k+O(|z|^3)\,,\quad g_{j,\bar{k},l,\bar{m}}(x_{\varepsilon})=-C^{j,\bar{k}}_{l, \bar m}\,,}&
\\
%&\displaystyle{
%{\cal C}_{\omega }(T_{_{X,J} } )(x_0) =
%\sum_{j,k,l,m} C^{j,\bar{k}}_{l,m}\,(dz_j\wedge d\bar{z}_k)\otimes dz_m\otimes_{_{J}}\frac{\partial }{\partial z_l}\,,}&
%\\
&\displaystyle{
{\cal C}_{\tilde\omega_{\varepsilon}}(x_{\varepsilon})=\sum_{j,k,l,m}C^{j,\bar{k}}_{l, \bar m}\,dz_j\otimes dz_l\otimes d\bar{z}_k\otimes d\bar{z}_m\,.}&
\end{eqnarray*}
and $\frac{i}{2}\partial\bar\partial u=\frac{i}{2}\sum_{l}u_{l,\bar{l}}\,dz_l\wedge d\bar{z}_l$, with $0< u_{1,\bar 1}\leq ...\leq u_{n,\bar n}$ at the point $x_{\varepsilon}$. For every $\zeta\in \C^n$ we set $g_{\zeta,\bar\zeta}:=\sum_{l,m}g_{l,\bar{m}}\,\zeta_l\,\bar\zeta_m$. Then
$$
\Lambda_{\varepsilon}(x)=\max_{\xi \in T_{X,x}\smallsetminus 0_x }\, \frac{\partial\bar\partial u(\xi^{1,0},\xi^{0,1})}{\partial\bar\partial g(\xi^{1,0},\xi^{0,1})}=\max_{|\zeta|=1 }\,\frac{u_{\zeta,\bar\zeta}}{g_{\zeta,\bar\zeta}}\,,
$$
and so $\Lambda_{\varepsilon}(x_{\varepsilon})=u_{n,\bar n}(x_{\varepsilon})$, with $\frac{u_{n,\bar n}}{g_{n,\bar n}}\leq \Lambda_{\varepsilon}$. We also set
$$
\tilde A_{\varepsilon}:=\log \frac{u_{n,\bar n}}{g_{n,\bar n}}-k\cdot \psi_{\varepsilon}+h\cdot b^{\varepsilon}\,.
$$
Then $\tilde A_{\varepsilon}\leq A_{\varepsilon}$, with $\tilde A_{\varepsilon}(x_{\varepsilon})=A_{\varepsilon}(x_{\varepsilon})$. This implies that the function $\tilde A_{\varepsilon}$ also reaches a maximum at the point $x_{\varepsilon}$, thus $\Delta_{\hat\omega_{\varepsilon}}\tilde A_{\varepsilon}(x_{\varepsilon})\leq 0$. All the subsequent computations in this part of the proof will be made at the point $x_{\varepsilon}$. 
By the local expressions for the Ricci tensor we obtain
\begin{eqnarray*}
\partial^2_{n,\bar n}\log\det(u_{j,\bar{k}})&=&\sum_{l,p}
\Big(u_{n,\bar n,l,\bar{p}}-\sum_{s,t}u_{n,l,\bar{s} }\,u^{s,\bar{t}}\,u_{\bar n,t,\bar{p}}\Big)u^{p,\bar{l}}
\\
&=&\sum_p\frac{u_{n,\bar n,p,\bar{p}}}{u_{p,\bar{p}}}  -\sum_{p,q}
\frac{\,|u_{n,p,\bar{q} }|^2}{u_{p,\bar{p}}\,u_{q,\bar{q}}}\,,
\end{eqnarray*}
and in a similar way $\partial^2_{n,\bar n}\log\det(g_{j,\bar{k}})=\sum_p\,g_{n,\bar n,p,\bar{p}}$. 
Then by differentiating with respect to the operator $\partial^2_{n,\bar n}$ the identity \eqref{MAsingEps2}, which can be rewritten as
$$
\log\det(u_{j,\bar{k}})=F^{\varepsilon}+\lambda\delta\log|s|^2 +\lambda(u-g)/2+\log\det(g_{j,\bar{k}})\,,
$$
we obtain
\begin{eqnarray*}
\sum_p\frac{u_{n,\bar n,p,\bar{p}}}{u_{p,\bar{p}}}
-
\sum_{p,q}
\frac{\,|u_{n,p,\bar{q} }|^2}{u_{p,\bar{p}}\,u_{q,\bar{q}}}
&=&
f^{\varepsilon}_{n,\bar n}+\lambda[(\omega_{\delta})_{n,n}-\omega_{n,n}]/2
\\
\\
&+&
l\cdot a^{\varepsilon}_{n,\bar n}-h\cdot b^{\varepsilon}_{n,\bar n}
\\
\\
&+&
\lambda(u_{n,\bar n}-1)/2+\sum_p\,g_{n,\bar n,p,\bar{p}}\,.
\end{eqnarray*}
Combining this with the inequality $\Delta_{\hat\omega_{\varepsilon}}\tilde A_{\varepsilon}(x_{\varepsilon})\leq 0$, we get
\begin{eqnarray*}
0&\geq& \sum_p\frac{\tilde A_{p,\bar{p}}}{u_{p,\bar{p}}}
\\
\\
&=&\sum_p\left(\frac{u_{n,\bar n,p,\bar{p}}}{u_{p,\bar{p}}\,u_{n,\bar{n}}}-\frac{\,|u_{n,\bar{n},p }|^2}{u_{p,\bar{p}}\,u^2_{n,\bar{n}}}+\frac{k/2+h\cdot b^{\varepsilon}_{p,\bar p}-g_{n,\bar n,p,\bar{p}}}{u_{p,\bar{p}}}\right)-nk/2
\\
\\
&=&\sum_{p,q}
\frac{\,|u_{n,p,\bar{q} }|^2}{u_{p,\bar{p}}\,u_{q,\bar{q}}\,u_{n,\bar{n}}}-\sum_p\frac{\,|u_{n,\bar{n},p }|^2}{u_{p,\bar{p}}\,u^2_{n,\bar{n}}}
\\
\\
&+&
\frac{f_{n,\bar{n}}+\lambda[(\omega_{\delta})_{n,n}-\omega_{n,n}-1]/2+l\cdot a^{\varepsilon}_{n,\bar n}-h\cdot b^{\varepsilon}_{n,\bar n}}{u_{n,\bar{n}}}
\\
\\
&+&\sum_p\left(\frac{g_{n,\bar n,p,\bar{p}}}{u_{n,\bar{n}}}+\frac{k/2+h\cdot b^{\varepsilon}_{p,\bar p}-g_{n,\bar n,p,\bar{p}}}{u_{p,\bar{p}}}\right)-(nk-\lambda)/2\,.
\end{eqnarray*}
We observe at this point that the sum of the two first terms following the second equality is nonnegative and the trivial inequality
$$
-\frac{h\cdot b^{\varepsilon}_{n,\bar n}}{u_{n,\bar{n}}}+\sum_p\frac{h\cdot b^{\varepsilon}_{p,\bar p}}{u_{p,\bar{p}}}\;\geq\; \sum_p\frac{-h\cdot K_{0}/2}{u_{p,\bar{p}}}\,.
$$
By plugging these inequalities in the previous computations and by using the definition of the constants $k$ and $K_1$, we get
\begin{eqnarray*}
0&\geq&
\sum_p\left(\frac{K_1-C_{p,\bar p}^{n,\bar n} }{u_{n,\bar{n}}}+
\frac{-K_1+ C_{p,\bar p}^{n,\bar n}}{u_{p,\bar{p}}}+\frac{1}{u_{p,\bar{p}}}\right)-(nk-\lambda)/2
\\
\\
&\geq&
\sum_p\frac{ (C_{p,\bar p}^{n,\bar n}-K_1)(u_{n,\bar{n}}-u_{p,\bar{p}})}{u_{p,\bar{p}}\,u_{n,\bar{n}}}+\sum_p\frac{1}{u_{p,\bar{p}}}-C_{0}\,,\qquad\qquad\qquad\qquad\;\,
\end{eqnarray*}
where $C_{0}>0$ and all the following constants are independents of $\varepsilon$. Let denote by $(x_1,...,x_n)$ the real part of the complex coordinates $(z_1,...,z_n)$. Then the inequality
$
C^{n,\bar{n}}_{p,\bar p}={\cal C}_{\tilde\omega_{\varepsilon}}
(\frac{\partial}{\partial x_n}\otimes \frac{\partial }{\partial x_p},\frac{\partial}{\partial x_n}\otimes \frac{\partial }{\partial x_p})(x_{\varepsilon})\geq \gamma_{\varepsilon}\geq \Gamma
$
implies
\begin{eqnarray*}
0\geq  \sum_p\frac{1}{u_{p,\bar{p}}}-C_{0}
&\geq&
\left(\frac{u_{n,\bar{n}}}
{\prod_p u_{p,\bar{p}}}  \right)^{\frac{1}{n-1}}-C_{0}
\\
\\
&=&
e^{\frac{-\lambda \psi_{\varepsilon} -\lambda\delta\log|s|^2 -F^{\varepsilon}}{n-1}}\;  u_{n,\bar{n}}^{\frac{1}{n-1}}-C_{0}\,.
\end{eqnarray*}
Consider now the function $B_{\varepsilon}:=e^{A_{\varepsilon}}=\Lambda_{\varepsilon} \,e^{-k\cdot\psi_{\varepsilon}+h\cdot b^{\varepsilon}}$. Then $x_{\varepsilon}$ is also a maximum point for $B_{\varepsilon}$ over $X\smallsetminus |D|$ and the previous inequality can be written as
\begin{eqnarray*}
0&\geq &
e^{\frac{(k-\lambda)\psi_{\varepsilon}-\lambda\delta\log|s|^2  -h\cdot b^{\varepsilon} -F^{\varepsilon}}{n-1}(x_{\varepsilon})}\,\,B_{\varepsilon}(x_{\varepsilon})^{\frac{1}{n-1}}-C_{0}
\\
\\
&=& 
e^{\frac{(k-\lambda)\varphi_{\varepsilon}-\delta k\log|s|^2-l\cdot a^{\varepsilon} -f^{\varepsilon}}{n-1}(x_{\varepsilon})}\,\,B_{\varepsilon}(x_{\varepsilon})^{\frac{1}{n-1}}-C_{0}\,.
\end{eqnarray*}
Then by the inequalities $k-\lambda>0,\;|s|^2\leq C$, $a^{\varepsilon}\leq C$ and $|f^{\varepsilon}|\leq K_{0}$, we get the estimate
$$
0\geq C_1\,e^{\frac{(k-\lambda)}{n-1}\min_X\varphi_{\varepsilon}}\,\,B_{\varepsilon}(x_{\varepsilon})^{\frac{1}{n-1}}-C_{0}\,.
$$ 
In conclusion we have found over $X\smallsetminus |D|$ the estimates
\begin{eqnarray*}
0&<&2n+\Delta _{\tilde\omega_{\varepsilon} }\varphi_{\varepsilon}-
\delta\Delta _{\tilde\omega_{\varepsilon} }\log|s|^2
=\Tr_{\tilde\omega_{\varepsilon}}(\tilde\omega_{\varepsilon}+i\partial\bar\partial\psi_{\varepsilon})
\\
\\
&\leq&
2n\Lambda_{\varepsilon}\leq 2n\,e^{k\cdot\psi_{\varepsilon}-h\cdot b^{\varepsilon}}B_{\varepsilon}(x_{\varepsilon})
\\
\\
&\leq&
\frac{C_2\,e^{k\cdot \varphi_{\varepsilon}-(k-\lambda)\min_X\varphi_{\varepsilon}}}{|s|^{2\delta k}(|\tau|^2+\varepsilon)^h} 
\leq
\frac{C_2\,e^{k\Osc(\varphi_{\varepsilon})}}{|s|^{2\delta k}\,|\tau|^{2h}}
\,. 
\end{eqnarray*}
(Here  $\Tr_{\tilde\omega_{\varepsilon}}$ is the trace operator with respect to the K\"{a}hler  metric $\tilde\omega_\varepsilon$.) The last inequality follows from the fact that $\lambda\min_X\varphi_{\varepsilon}\leq 0$, since a non identically zero solution $\varphi_{\varepsilon}$ changes signs in the case $\lambda>0$. Then using the inequality 
$$
\left|\delta\Delta _{\tilde\omega_{\varepsilon} }\log|s|^2\right|
=
|\Tr_{\tilde\omega_{\varepsilon} }(\omega-\omega_{\delta})|\leq C_3
$$ 
over $X\smallsetminus |D|$, we deduce the singular Laplacian estimate
$$
-C_3<2n+\Delta _{\tilde\omega_{\varepsilon} }\varphi_{\varepsilon}
\leq
\frac{C_2\,e^{k\Osc(\varphi_{\varepsilon})}}{|s|^{2\delta k}\,|\tau|^{2h}}+C_3\,.
$$
{\bf (C) Higher order estimates.}
\\
By the previous estimates we infer 
$0<2u_{l,\bar l}<\Tr_{\tilde\omega_{\varepsilon}}(\omega+\varepsilon\alpha+i\partial\bar\partial \varphi_{\varepsilon})\leq 2\Upsilon_{\varepsilon}$, for all $l=1,...,n$, with 
$$
\Upsilon_{\varepsilon}:=\frac{C_2\,e^{k\Osc(\varphi_{\varepsilon})}}{|s|^{2\delta k}(|\tau|^2+\varepsilon)^h}\,.
$$
The equation \eqref{MAsingEps} can be rewritten as 
$$
(\omega+\varepsilon\alpha+i\partial\bar\partial \varphi_{\varepsilon})^n=e^{F^{\varepsilon}+\lambda\varphi_{\varepsilon}}\,\tilde\omega_{\varepsilon}^n\,.
$$
We infer
$$
e^{F^{\varepsilon}+\lambda\varphi_{\varepsilon}}=\prod_l u_{l,\bar l}\leq \Upsilon_{\varepsilon}^{n-1}u_{m,\bar m}\,,
$$
for all $m=1,...,n$. The fact that a non identically zero solution $\varphi_{\varepsilon}$ changes signs in the case $\lambda>0$ implies $\lambda \min_X \varphi_{\varepsilon}\geq -\lambda \Osc(\varphi_{\varepsilon})$. Thus
$$
e^{F^{\varepsilon}-\lambda\Osc(\varphi_{\varepsilon})}\,\Upsilon_{\varepsilon}^{1-n}\,\tilde\omega_{\varepsilon} \leq \omega+\varepsilon\alpha+i\partial\bar\partial \varphi_{\varepsilon}\,.
$$
Then an elementary computation yields the singular  estimate
\begin{eqnarray}\label{unifMetrc}
&&\kern-30pt C_4^{-1}\,|s|^{2\delta k(n-1)}\,|\sigma|^{2l}\,|\tau|^{2h(n-2)}\,e^{-kn\Osc(\varphi_{\varepsilon})}\,\tilde\omega_{\varepsilon} \nonumber
\\\nonumber
\\
&\leq&
\omega+\varepsilon\alpha+i\partial\bar\partial \varphi_{\varepsilon}
\;\leq\;
\frac{C_4\,e^{k\Osc(\varphi_{\varepsilon})}}{|s|^{2\delta k}\,|\tau|^{2h}}
\,\tilde\omega_{\varepsilon}\,.
\end{eqnarray}
Moreover the fact that $\varphi_{\varepsilon}\in {\cal P}_{\omega+\varepsilon\alpha}$ implies
$$
2|\partial\bar\partial \varphi_{\varepsilon}|_{\tilde\omega_{\varepsilon}}\leq \Delta_{\tilde\omega_{\varepsilon}}\varphi_{\varepsilon}+2\Tr_{\tilde\omega_{\varepsilon} }(\omega+\varepsilon\alpha)\,.
$$
At this step of the proof we consider
\begin{eqnarray*}%\label{SBlocus}
S':=|D|\cup \Big(\bigcup_r\{\tau_r=0\}\Big)\,,\qquad
S:=S'\cup\Big(\bigcup_j\{\sigma_j=0\}\Big).
\end{eqnarray*}
By the Interpolation Inequalities \cite{Gi-Tru} we find that for any coordinate compact set $K\subset X\smallsetminus S'$ there are uniform constants $C_K>0$ such that
$$
\max_K|\nabla_{_{\C^n}} \varphi_{\varepsilon}|\leq C_K
\left(\max_K \Delta_{_{\C^n}} \varphi_{\varepsilon}+\max_K|\varphi_{\varepsilon}|\right)\,.
$$
Therefore, we can apply the complex version of Evans-Krylov theory \cite{Ti2} on every compact set $K\subset X\smallsetminus S$ to get uniform constants $C_{2,K}>0$ such that $\|\varphi_{\varepsilon}\|_{C^{2,\eta}(K)}\leq C_{2,K}$ for some $\eta\in (0,1)$. Now, let $U\subset X\smallsetminus S$ be an open set and let $\xi\in {\cal O}(T_{X})(U)$. We rewrite the complex Monge-Amp\`ere equation \eqref{MAsingEps} under the form
\begin{eqnarray*}
(\omega+\varepsilon\alpha+i\partial\bar\partial \varphi_{\varepsilon})^n=e^{H_{\varepsilon}+\lambda\varphi_{\varepsilon}}\alpha^n\,,
\end{eqnarray*}
with 
$$
H_{\varepsilon}:=c_{\varepsilon}+\log(\Omega/\alpha^n)
+l\cdot a^{\varepsilon}-h\cdot b^{\varepsilon}\,.
$$
By taking the derivative  with respect to the complex vector field $\xi$, we obtain (see the proof of formula 11 in \cite{Pal})
\begin{eqnarray}\label{linEllipticReg}
\Delta_{\hat\omega_{\varepsilon}} (\xi\,.\,\varphi_{\varepsilon})-2\lambda\,\xi\,.\,\varphi_{\varepsilon}
=
-\Tr_{\hat\omega_{\varepsilon}}
L_{\xi}(\omega+\varepsilon\alpha)+\Tr_{\alpha}L_{\xi}\alpha
+
2\xi\,.\,H_{\varepsilon}\,,
\end{eqnarray}
By the uniform estimates \eqref{unifMetrc} and $\|\varphi_{\varepsilon}\|_{C^{2,\eta}(K)}\leq C_{2,K}$ it follows that the operator $\Delta_{\hat\omega_{\varepsilon}}$ is uniformly elliptic with coefficients uniformly bounded in $C^{\eta}$-norm at least, over any compact set $K\subset U$.
The right hand side of equation \eqref{linEllipticReg} is also uniformly bounded in $C^{\eta}$-norm at least, over $K$. By the standard regularity theory for linear elliptic equations \cite{Gi-Tru} we deduce $\|\xi\,.\,\varphi_{\varepsilon}\|_{{C^{2,\eta}(K)}}\leq C'_K$ for all $\varepsilon>0$. %By conjugation the same holds for $\bar\xi\,.\,\varphi_{\varepsilon}$. 
We infer the uniform estimate $\|\varphi_{\varepsilon}\|_{C^{3,\eta}(K)}\leq C_{3,K}$.
\\
In its turn, this estimate implies that the coefficients of the Laplacian $\Delta_{\hat\omega_{\varepsilon}}$ and the right hand side of equation \eqref{linEllipticReg} are uniformly bounded in $C^{1,\eta}$-norm at least. By iteration we get the uniform estimates $\|\varphi_{\varepsilon}\|_{C^{r,\eta}(K)}\leq C_{r,K}$ for all $\varepsilon>0$ and $r\in \N$. 
We infer that the family $(\varphi_{\varepsilon})_{\varepsilon>0}\subset C^{\infty}(X\smallsetminus S)$ is precompact in the smooth topology. 
\\
\\
{\bf (CI) The smooth regularity.} 
\\
By elementary properties of plurisubharmonic functions (see \cite{Dem2}, chapter 1), the uniform estimate $\|\varphi_{\varepsilon}\|_{L^{\infty}(X)}\leq C$ implies the existence of a $L^1$-convergent sequence $(\varphi_j)_j$, $\varphi_j:=\varphi_{\varepsilon_j}$, $\varepsilon_j\downarrow 0$ with limit $\varphi\in {\cal P}_{\omega}\cap L^{\infty}(X)$. We can assume that a.e.-convergence holds also.
The precompactness of the family $(\varphi_{\varepsilon})_{\varepsilon>0}\subset C^{\infty}(X\smallsetminus S)$ in the smooth topology implies the convergence of the limits
\begin{eqnarray}\label{conv-MA-limit}
(\omega+i\partial\bar\partial \varphi)^n=\lim_{l\rightarrow +\infty}(\omega+i\partial\bar\partial \varphi_j)^n=\lim_{j\rightarrow +\infty}G_j\,e^{\lambda\,\varphi_j+c_j}\Omega=G_0\,e^{\lambda\,\varphi}\Omega
\end{eqnarray}
over $X\smallsetminus S$, with $c_j:=c_{\varepsilon_j}$, $G_j:=G_{\varepsilon_j}$ and
$$
G_{\varepsilon}:=(|\sigma|^2+\varepsilon^A)^l\,(|\tau|^2+\varepsilon)^{-h}\,.
$$
The fact that $\varphi$ is a bounded potential implies that the global complex Monge-Amp\`ere measure $(\omega+i\partial\bar\partial \varphi)^n$ does not carry any mass on complex analytic sets. We infer that $\varphi$ is a global bounded solution of the complex Monge-Amp\`ere equation \eqref{MainMA} which belongs to the class ${\cal P}_{\omega}\cap L^{\infty}(X)\cap C^{\infty}(X\smallsetminus S)$.
\\
\\
{\bf (CII) The $C^{1,1}$-regularity.} 
\\
Let $U\subset \subset X\smallsetminus S'$ be a coordinate open set. By a classical result in \cite{Gi-Tru} (see Theorem 8.32, page 210) for all open sets $U'\subset\subset U$ there exists a constant $C=C(U',U)>0$ such that for all $\eta\in (0,1)$ the uniform estimate
\begin{eqnarray}\label{C1hold-est-Gi-Tr}
\|\varphi_{\varepsilon}\|_{C^{1,\eta}(U')}\leq C\left(\|\varphi_{\varepsilon}\|_{L^{\infty}(U)}+\|\Delta_{_{\C^n}}\varphi_{\varepsilon}\|_{L^{\infty}(U)}\right)
\end{eqnarray}
holds. By applying the Ascoli-Arzela theorem to the sequence $(\varphi_j)_j$, we infer
the uniform estimate
$\|\varphi\|_{C^{1,\eta}(U')}\leq C'$ for all $\eta\in (0,1)$, thus $\varphi\in C^{1,1}(X\smallsetminus S')$. 
%The existence and regularity statement in the case $n=1$ follows immediately from the properties of the Green function.
\\
\\
{\bf (D) Uniqueness of the solution.}
\\
We now prove the uniqueness of the solution $\varphi$ in the class ${\cal P}\BT_{\omega}$. In the case $\lambda=0$ this follows immediately from Theorem \ref{UniqSol}. In the case $\lambda>0$ let $\psi\in {\cal P}\BT_{\omega}$ be an other solution. The fact that $\psi\in {\cal P}_{\omega}$ implies that we can solve the degenerate complex Monge-Amp\`ere equation 
\begin{eqnarray}\label{UnqMAEx}
(\omega+i\partial\bar\partial  u)^n=G_0\,e^{\lambda \psi}\,\Omega\,,
\end{eqnarray}
with the methods so far explained, so as to obtain a solution $u\in {\cal P}^0_{\omega}\cap L^{\infty}(X)$. 
In fact we consider the solutions $u_{\varepsilon}$ of the non-degenerate complex Monge-Amp\`ere equations
$$
(\omega+\varepsilon\alpha+i\partial\bar\partial  u_{\varepsilon})^n=G_{\varepsilon}\,e^{\lambda \psi_{\varepsilon}+c'_{\varepsilon}}\,\Omega\,,
$$
with $\psi_{\varepsilon}\downarrow \psi$, $\psi_{\varepsilon}\in C^{\infty}(X)$, $\psi_{\varepsilon}\leq C$, $i\partial\bar\partial  \psi_{\varepsilon}\geq -K_0\,\tilde\omega_{\varepsilon}$ and $c'_{\varepsilon}$ being a normalizing constant converging to $0$ as $\varepsilon\rightarrow 0$. By combining Lemma  \ref{IntegLp} with the dominated convergence theorem we infer that the family $G_{\varepsilon}\,e^{\lambda \psi_{\varepsilon}+c'_{\varepsilon}}$ converges in $L^p$-norm  to $G_0\,e^{\lambda \psi}$. These conditions are sufficient to provide the singular Laplacian estimate of step (B). Thus by the $C^{1,\eta}$-compactness argument of step (CII) we infer the existence of the solution $u$ of the degenerate complex Monge-Amp\`ere equation \eqref{UnqMAEx}.
\\
By the uniqueness result in the case $\lambda=0$ we infer $u=\psi-\sup_X\psi$, thus $\psi\in L^{\infty}(X)$. Then the required uniqueness follows immediately from Theorem~\ref{UniqSol} (B).
\\
\\
{\bf (E) Eliminating the assumption on the existence of divisors $D$~in~$X$.}
\\
By section 5, the divisors $D$ which we have assumed to exist in $X$ up to now, can only be constructed (at least, in the non-projective case) by applying a blow-up process to $X$, i.e.\ we can find a modification $\mu:\tilde X\to X$ of~$X$, a divisor $D$ in $\tilde X$ with $|D|\supset \Exc(\mu)$ and a number $\delta>0$ such that the class $\{\mu^*\omega\}-\delta\{D\}$ is K\"ahler on~$\tilde X$. For this reason, we use pull-back the Monge-Amp\`ere equation by $\mu$ so as to transform equation \eqref{MainMA} into
$$
(\mu^*\omega+i\partial\bar\partial \Phi)^n= \prod\limits_{j=1}^N|\sigma_j\circ\mu|^{2l_j}\cdot \prod\limits_{r=1}^M|\tau_r\circ\mu|^{-2h_r}\,e^{\lambda\Phi}\,\mu^*\Omega\,,\quad \lambda\geq 0\,.
$$
Here $\mu^*\Omega$ is no longer a positive volume form on $\tilde X$ but we have \hbox{$\mu^*\Omega=|J\mu|^2\tilde\Omega$} where $\tilde\Omega$ is such a volume form, and $|J\mu|^2$ is the square of the Jacobian of~$\mu$ expressed with respect to the pair $(\Omega,\tilde\Omega)$. Observe that $J\mu$ is just a section of the relative canonical divisor $\smash{K_{\tilde X/X}}$ and that $|J\mu|^2$ is its norm with respect to the metric induced by $(\Omega,\smash{\tilde\Omega)}$. Thus our equation again takes the form
$$
(\mu^*\omega+i\partial\bar\partial \Phi)^n= |J\mu|^2\prod\limits_{j=1}^N|\sigma_j\circ\mu|^{2l_j}\cdot \prod\limits_{r=1}^M|\tau_r\circ\mu|^{-2h_r}\,e^{\lambda\Phi}\,\tilde\Omega\,,
$$
and it is clear that the analogue of condition \eqref{IntgCond} 
$$
0<\int\limits_{\tilde X}\mu^*\omega^n
=
\int\limits_{\tilde X}|J\mu|^2\prod\limits_{j=1}^N|\sigma_j\circ\mu|^{2l_j}\cdot \prod\limits_{r=1}^M|\tau_r\circ\mu|^{-2h_r}\,\tilde\Omega\,,
$$
holds. By steps (A)--(D), we obtain a unique solution
$$
\Phi\in {\cal P}_{\mu^*\omega}\cap L^{\infty}(\tilde X)\cap C^{1,1}(\tilde X\smallsetminus \tilde S'_{\mu,D,\delta})\cap C^{\infty}(\tilde X\smallsetminus \tilde S_{\mu,D,\delta})\,,
$$
with
\begin{eqnarray*}
\tilde S'_{\mu,D,\delta}&=&|D|\cup \Big(\bigcup_r\{\tau_r\circ\mu=0\}\Big)\,,
\\
\\
\tilde S_{\mu,D,\delta}&=&\tilde S'_{\mu,D,\delta}\cup\Big(\bigcup_j\{\sigma_j\circ\mu=0\}\Big)\cup\Exc(\mu).
\end{eqnarray*}
Actually, taking the union with $\Exc(\mu)$ will not be needed since $|D|\supset \Exc(\mu)$. 
Moreover $j_q^*\mu^*\omega=0$, where $j_q:\mu^{-1}(q)\hookrightarrow \tilde X$, $q\in \mu(\Exc(\mu))$ is the inclusion map. Thus 
$$
\Phi\circ j_q \in \Psh(\mu^{-1}(q))
$$ 
since $\Phi\in {\cal P}_{\mu^*\omega}\cap L^{\infty}(\tilde X)$.
By hypothesis $\mu^{-1}(q)$ is compact and connected, which implies that $\Phi$ is constant along the fibers $\mu^{-1}(q)$. Therefore we can define $\varphi:=\pi_*\Phi\in {\cal P}_{\omega}\cap L^{\infty}(X)$. The fact that $\varphi$ is bounded implies that the current 
$(\omega+i\partial\bar\partial  \varphi)^n$ does not carry any mass on complex analytic sets. This, combined with the fact that
$$
\mu:\tilde X\smallsetminus\Exc(\mu)\to X\smallsetminus\mu(\Exc(\mu))
$$ 
is a biholomorphism, implies (see Theorem \ref{UniqSol}) that $\varphi$ is the unique solution in ${\cal P}_{\omega}\cap L^{\infty}(X)$ of the complex Monge-Amp\`ere equation \eqref{MainMA} with the required $C^{1,1}$, $C^{\infty}$-regularity over the adequate subsets of $X\smallsetminus \mu ( |D| )$. We set finally
\begin{eqnarray*}
\Sigma_{\{\omega\}}&=&\bigcap_{(\mu,D,\delta)\in I_{\{\omega\}}}\mu(|D|)\,,
\\
\\
S'&=&\Sigma_{\{\omega\}}\cup \Big(\bigcup_r\{\tau_r=0\}\Big)\,,\qquad
S=S'\cup\Big(\bigcup_j\{\sigma_j=0\}\Big)\,.
\end{eqnarray*}
Then the conclusion about the ${\cal P}_{\omega}\cap L^{\infty}(X)\cap C^{1,1}(X\smallsetminus S')\cap C^{\infty}(X\smallsetminus S)$ regularity of the solution $\varphi$ 
follows by letting $(\mu, D, \delta)\in I_{\{\omega\}}$ vary. 
The proof of the uniqueness of the solution $\varphi$ in the class ${\cal P}\BT_{\omega}$ is the same as in step~D, modulo the use of modifications.
\\
\\
{\bf (F) $C^0$ regularity on $X\smallsetminus\Sigma_{\{\omega\}}$.}
\\
The proof will be complete if we show that $\varphi\in C^0(X\smallsetminus \Sigma_{\{\omega\}})$. This follows from the following statement.\hfill$\Box$
%In the case $n=1$ the solution $\varphi$ of the equation $\omega+i\partial\bar\partial \varphi=\Omega>0$ is always smooth.
%%%%%%%%%%%%%%%%%%%%%%%%%%%%%%%%%%%%%%%%%%%%%%%%%%%%%%%%%%%%%%%%%%%%%%%%%%%%%%%%%%%%%%%%%%%%%%%%%%%%%%%%%%%%%%%%%%%%%%%%%%%%%%%%%%%%%%%%%%%%%%%%%%%%%%%%%%%%%%%%%%%%%%%%%%%%%%%%%%%%%%%%%%%%%%%%%%%%%%%%%%
\begin{theorem}.\label{ExistSol}
Let $X$ be a compact connected K\"ahler manifold of complex dimension \hbox{$n\geq 2$}, let $\omega\geq 0$ be a big closed smooth $(1,1)$-form and let $\Omega>0$ be a smooth volume form.
Let also $f\in L\log^{n+\delta}L(X)$, $\delta>0$ such that $\int_X\omega^n=\int_Xf\,\Omega$ and $\lambda\geq 0$ be a real number. Then there exists a unique solution $\varphi\in{\cal P}BT_{\omega}$ of the degenerate complex Monge-Amp\`ere equation
\begin{eqnarray}\label{DG-MA-RHS}
(\omega+i\partial\bar\partial \varphi)^n=f\,e^{\lambda\,\varphi}\Omega\,,
\end{eqnarray}
which in the case $\lambda=0$ is normalized by $\sup_X\varphi=0$. 
The solution $\varphi$ is in the class ${\cal P}_{\omega}\cap L^{\infty}(X)\cap C^0(X\smallsetminus \Sigma_{\{\omega\}})$ and satisfies the $L^{\infty}$-estimate
$$
\|\varphi\|_{L^{\infty}(X)}\leq C(\delta, \omega, \Omega)\,I_{\omega,\delta}(f)^{\frac{n}{\delta}}+1\,.
$$
Moreover the constant $C(\delta, \omega, \Omega)>0$ stays bounded for perturbations of $\omega\geq 0$ as in Statement $(C)$ of Theorem~\ref{Kolo}.
\end{theorem}
$Proof$. We consider a regularizing family $(f_j)_j\subset C^{\infty}(X)$, $f_j>0$ of $f$ in $L\log^{n+\delta}L(X)$. (The existence of such family follows from \cite{Ra-Re} page 364 or \cite{Iw-Ma}, Theorem 4.12.2, page 79.) We can assume as usually
$\int_X \omega^n=\int_Xf_j\,\Omega$. 
By the proof of Theorem \ref{HgRegMA} and the $L^{\infty}$-estimate in corollary \ref{C0est} we deduce the existence of a  unique solution of the degenerate complex Monge-Amp\`ere equation 
\begin{eqnarray}\label{DG-MA-RHS-j}
(\omega+i\partial\bar\partial \varphi_j)^n=f_j\,e^{\lambda\,\varphi_j}\Omega\,,
\end{eqnarray}
with the properties $\varphi_j\in {\cal P}_{\omega}\cap L^{\infty}(X)\cap C^{\infty}(X\smallsetminus \Sigma_{\{\omega\}})$ and 
\begin{eqnarray}\label{C0-DG-MA-RHS-j}
\|\varphi_{j}\|_{L^{\infty}(X)}\leq C:=C(\delta, \omega, \Omega)\,I_{\omega,\delta}(f)^{\frac{n}{\delta}}+1\,.
\end{eqnarray}
(With $\sup_X\varphi_j=0$ in the case $\lambda=0$.) We deduce in particular the uniform estimate
\begin{eqnarray}\label{LlgLf-ST}
\|f_j\,e^{\lambda\,\varphi_j}\|_{L\log^{n+\delta}L(X)}\leq Ke^{\lambda C}\|f\|_{L\log^{n+\delta}L(X)}\,,
\end{eqnarray} 
for all $j$. 
(See \cite{Ra-Re} page 364 or \cite{Iw-Ma}, Theorem 4.12.2, page 79.) 
On the other hand the uniform estimate \eqref{C0-DG-MA-RHS-j} implies (see \cite{Dem2}, chapter 1) the existence of a $L^1$-convergent subsequence $(\varphi_j)_j$ (which by abuse of notation we denote in the same way).
We can apply the $L^{\infty}$-stability estimate of Theorem \ref{Kolo}~(B) to the complex Monge-Amp\`ere equation \eqref{DG-MA-RHS-j} thanks to the estimates \eqref{C0-DG-MA-RHS-j} and \eqref{LlgLf-ST}. Notice that by \eqref{C0-DG-MA-RHS-j}, the $L^{\infty}$-stability estimate of Theorem \ref{Kolo}~(B) applies even if in the case $\lambda>0$, when the solutions $\varphi_j$ are not necessarily normalized by the supremum condition.
We infer that the sequence $(\varphi_j)_j$ is a Cauchy sequence in the $L^{\infty}(X)$-norm, thus convergent to some function
$ \varphi\in  {\cal P}_{\omega}\cap L^{\infty}(X)\cap C^0(X\smallsetminus \Sigma_{\{\omega\}})$. This yields weakly convergent limits
$$
(\omega+i\partial\bar\partial \varphi)^n=\lim_{j\rightarrow +\infty}(\omega+i\partial\bar\partial \varphi_j)^n=\lim_{j\rightarrow +\infty}f_j\,e^{\lambda\,\varphi_j}\Omega=f\,\,e^{\lambda\,\varphi}\Omega\,,
$$
over $X\smallsetminus \Sigma_{\{\omega\}}$. Moreover the fact that the global Monge-Amp\`ere measure $(\omega+i\partial\bar\partial \varphi)^n$ does not carry any mass on complex analytic sets of $X$ implies that $\varphi$ is the unique (in the class ${\cal P}BT_{\omega}$) global solution of the degenerate complex Monge-Amp\`ere equation \eqref{DG-MA-RHS} with the required regularity and with $\|\varphi\|_{L^{\infty}(X)}\leq C$. (We remark that the uniqueness of the solution in the case $\lambda>0$ follows from the same argument in step (D) in the proof of the Theorem \ref{HgRegMA} .)\hfill$\Box$
\\
\\
%%%%%%%%%%%%%%%%%%%%%%%%%%%%%%%%%%%%%%%%%%%%%%%%%%%%%%%%%%%%%%%%%%%%%%%%%%%%%%%%%%%%%%%%%%%%%%%%%%%%%%%%%%%%%%%%%%%%%%%%%%%%%%%%%%%%%%%%%%%%%%%%%%%%%%%%%%%%%%%%%%%%%%%%%%%%%%%%%%%%%%%%%%%%%%%%%%%%%%%%%%
\emph{Proof of Theorem} \ref{GenTyp}. 
\\
A result of Kawamata \cite{Kaw} claims that in our case the canonical bundle is base point free, and so, for all $m\gg 0$ sufficiently big and divisible, $mK_X$ has no base points. So we can fix $m$ such that the pluricanonical map $f_m:X\rightarrow \C\proj^N$ is holomorphic. Consider also the semipositive and big K\"{a}hler form  $\omega_m:=f_m^*\omega_{FS}/m\in 2\pi c_1(K_X)$, where $\omega_{FS}$ is the Fubini-Study metric of $\C\proj^N$. 
Let $\Omega>0$ be a smooth volume form over $X$ such that $\int_X\Omega=\int_X\omega_m^n$ and $\Ric(\Omega)=-\omega_m$ (these conditions prescribe $\Omega$ in a unique way). According to Theorem \ref{HgRegMA} we can find a unique solution 
$\varphi\in {\cal P}\BT_{\omega_m}$ of the degenerate complex Monge-Amp\`ere equation
\begin{eqnarray*}
(\omega_m+i\partial\bar\partial \varphi)^n=e^{\varphi}\, \Omega\,.
\end{eqnarray*}
Moreover $\varphi\in {\cal P}_{\omega_m}\cap L^{\infty}(X)\cap C^{\infty}(X\smallsetminus \Sigma_{\{\omega_m\}})$, and so $\omega_{_{E}}:=\omega_m+i\partial\bar\partial \varphi$ is the required unique Einstein current in the class $\BT^{\log}_{\, 2\pi c_1(K_X)}$.\hfill$\Box$
\\
\\
%%%%%%%%%%%%%%%%%%%%%%%%%%%%%%%%%%%%%%%%%%%%%%%%%%%%%%%%%%%%%%%%%%%%%%%%%%%%%%%%%%%%%%%%%%%%%%%%%%%%%%%%%%%%%%%%%%%%%%%%%%%%%%%%%%%%%%%%%%%%%%%%%%%%%%%%%%%%%%%%%%%%%%%%%%%%%%%%%%%%%%%%%%%%%%%%%%%%%%%%%%
\emph{Proof of Theorem} \ref{G-GenTyp}.
%%%%%%%%%%%%%%%%%%%%%%%%%%%%%%%%%%%%%%%%%%%%%%%%%%%%%%%%%%%%%%%%%%%%%%%%%%%%%%%%%%%%%%%%%%%%%%%%%%%%%%%%%%%%%%%%%%%%%%%%%%%%%%%%%%%%%%%%%%%%%%%%%%%%%%%%%%%%%%%%%%%%%%%%%%%%%%%%%%%%%%%%%%%%%%{\bf Appendix D. Proof of the conjecture of Tian.} 
\\
The uniqueness statement in the theorem \ref{G-GenTyp} follows from the corollary \ref{cmp-KE}. In order to prove the existence of a K\"ahler-Einstein current $\omega_{_{E}}\in 2\pi c_1(K_X)$
let $m$ be a sufficiently large integer such that the base locus of $mK_X$ coincides with the stable base locus $\SB$ and let
$$
f_m:X\smallsetminus \SB\longrightarrow X_m:=\overline{f_m(X\smallsetminus \SB)}\,,
$$
be the rational map associated to the linear system $H^0(X,mK_X)$.
Let $\hat \Gamma_m$ be the desingularization of the Zariski closure of the graph $\Gamma_m\subset X\times X_m$ of $f_m$, let $\pi_m:\hat \Gamma_m\rightarrow X$ and 
$p_m:\hat \Gamma_m\rightarrow X_m$ be the natural projections. By definition of the graph there exists a Zariski dense open set $U_m\subset \hat \Gamma_m$ such that $X\smallsetminus \SB=\pi_m(U_m)$ and $p_m=f_m\circ \pi_m$ over $U_m$. 
Consider also bases
$$
(\sigma_{m,j})_{j=1}^{N_m}\subset H^0(X,mK_X)\,,
$$ 
and the induced curvature currents
$$
0\;\le \;\gamma_m\;:=\;\frac{1}{r_m}\,f_m^*\,\omega_{FS,m}\;=\;-\Ric(\Omega_m)\in 2\pi c_1(K_X)\,,
$$
where $\omega_{FS,m}$ is the Fubini-Study metric of $\C\proj^{N_m-1}$ and $\Omega_m^{-1}$ is the induced singular hermitian metric over $mK_X$. Explicitly
$$
\Omega_m\;=\;\left(\;\sum_{j=1}^{N_m}\;\left|\frac{\sigma_{m,j}}{\kappa^{m}}\right|^2\right)^{1/m}i^{n^2}\kappa\wedge \bar\kappa
\;=\;
\left(\;\sum_{j=1}^{N_m}\;\left|\sigma_{m,j}\right|^2_{{\Omega^{-1}}}\right)^{1/m}\Omega\,,
$$ 
for arbitrary $\kappa\in H^0(X,K_X)$ and $\Omega>0$ a smooth volume form. Observe now that the smooth form 
$$
0\le \theta_m:=m^{-1}p_m^*\,\omega_{FS,m}\,,
$$ 
is big. 
Moreover the Zariski dense open set 
$
V_m:= \hat \Gamma_m\smallsetminus \Sigma_{\theta_m}%\,,
$
satisfies $X\smallsetminus \Sigma=\pi_m(V_m)$.
By Theorem \ref{HgRegMA} we infer the existence of a solution 
$$
\Phi_m\in ({\cal P}_{\theta_m}\cap L^{\infty})(\hat \Gamma_m)\cap C^{\infty}(V_m)\,,
$$ 
of the degenerate complex Monge-Amp\`ere equation
\begin{eqnarray}\label{dMA1}
(\theta_m+i\partial\bar\partial \Phi_m)^n=e^{\Phi_m}\,\pi_m^*\Omega_m\,, 
\end{eqnarray}
over $\hat \Gamma_m$. The fact that $\theta_m=\pi_m^*\gamma_m$ over $U_m$ and the fibers of $\pi_m$ are connected allows to $\pi_m$-push forward the equation \eqref{dMA1}. We infer a solution $\varphi_m\in L^{\infty}(X)\cap C^{\infty}(X\smallsetminus \Sigma)$
of the degenerate complex Monge-Amp\`ere equation
\begin{eqnarray}\label{dMA2}
(\gamma_m+i\partial\bar\partial\varphi_m)^n=e^{\varphi_m}\Omega_m \,,
\end{eqnarray}
over $X\smallsetminus \SB$. We observe that \eqref{dMA2} can be rewritten in an equivalent way as
\begin{eqnarray*}%\label{dMA3}
\left(-\Ric(\Omega)+i\partial\bar\partial\psi_m\right)^n=e^{\psi_m}\Omega \,,
\end{eqnarray*}
over $X\smallsetminus \SB$, with 
$$
\psi_m:=\varphi_m+m^{-1}\log \sum_{j=1}^{N_m}\left|\sigma_{m,j}\right|^2_{{\Omega^{-1}}}\,.
$$ 
Thus $\omega_{_{E}}:=-\Ric(\Omega)+i\partial\bar\partial\psi_m$ is the required K\"ahler-Einstein current.
\hfill$\Box$
\\
\\
\emph{Proof of the conjecture of Tian} \ref{ti-conj}
\\
The hypothesis (C1) of Statement (C) in Theorem \ref{Kolo} is obviously satisfied. The hypothesis (C2b) is also satisfied since
$$
\lim_{t\rightarrow 0}\,\frac{(\pi^*\omega_Y+t\omega_X)^n}{K_t\,\omega^n_X}
=
\left(\;\int\limits_{y\in Y}\omega^m_Y(y)\cdot\kern-12pt\int\limits_{z\in \pi^{-1}(y)}\omega^{n-m}_X\right)^{-1}\frac{\pi^*\omega_Y^m\wedge \omega^{n-m}_X }{\omega^n_X}<+\infty.
$$ 
We deduce $\Osc(\psi_t)\leq C<+\infty$ for all $t\in (0,1)$ by Statements (C) and (A) of Theorem \ref{Kolo}. This solves in full generality the conjecture of Tian \ref{ti-conj}. \hfill$\Box$
\section{Appendix}
{\bf Appendix A. Computation of a complex Hessian.}
Let $\sigma\in H^0(X,E)$ be a holomorphic section of a holomorphic hermitian vector bundle $(E,h)$ and set $S_{\varepsilon}:=\log(|\sigma|^2+\varepsilon)$, for some $\varepsilon>0$. We denote by $\{\cdot,\cdot\}$ the exterior product of $E$-valued forms respect to the hermitian metric $h$. We have
$$
i\partial S_{\varepsilon}=\frac{i\{\partial_h \sigma,\sigma\}}{|\sigma|^2+\varepsilon}\,,
$$
since $\sigma$ is a holomorphic section. We compute now the complex hessian
\begin{eqnarray*}
i\partial\bar\partial S_{\varepsilon}
&=&
-\bar\partial\,\frac{i\{\partial_h \sigma,\sigma\}}{|\sigma|^2+\varepsilon}
\\
\\
&=&%\{,\}
\frac{-i\{\bar\partial \partial_h \sigma,\sigma\}+i\{\partial_h \sigma,\partial_h \sigma\}}{|\sigma|^2+\varepsilon}
\,+\,
i\{\partial_h \sigma,\sigma\}\wedge \bar\partial\left(\frac{1}{|\sigma|^2+\varepsilon}\right)
\\
\\
&=&
\frac{i\{\partial_h \sigma,\partial_h \sigma\}
\,-\,
\{i{\cal C}_{E,h}\sigma,\sigma\}}{|\sigma|^2+\varepsilon}
\,-\,
\frac{i\{\partial_h \sigma,\sigma\}\wedge\{\sigma,\partial_h \sigma\} }{(|\sigma|^2+\varepsilon)^2}
\\
\\
&=&
\underbrace{\frac{(|\sigma|^2+\varepsilon)i\{\partial_h \sigma,\partial_h \sigma\}-i\{\partial_h \sigma,\sigma\}\wedge\{\sigma,\partial_h \sigma\} }{(|\sigma|^2+\varepsilon)^2}}_{iT(S_{\varepsilon})}
\,-\,
\frac{\{i{\cal C}_{E,h}\sigma,\sigma\}}{|\sigma|^2+\varepsilon}
\end{eqnarray*}
where ${\cal C}_{E,h}\in C^\infty(X,\Lambda^{1,1}T^*_X\otimes\End(E,E))$ is the curvature tensor of $(E,h)$. We show that the $(1,1)$-form $iT(S_{\varepsilon})$ is nonnegative. In fact by using twice the Lagrange inequality
$$
i\{\partial_h \sigma,\sigma\}\wedge\{\sigma,\partial_h \sigma\} \leq |\sigma|^2\,i\{\partial_h \sigma,\partial_h \sigma\}
$$
(which is an equality in the case of line bundles), we get
\begin{eqnarray*}
iT(S_{\varepsilon}) 
\geq 
\frac{\varepsilon i\{\partial_h \sigma,\partial_h \sigma\} }{(|\sigma|^2+\varepsilon)^2}
\geq
\frac{\varepsilon i\{\partial_h \sigma,\sigma\}\wedge\{\sigma,\partial_h \sigma\} }{|\sigma|^2(|\sigma|^2+\varepsilon)^2}=\frac{\varepsilon}{|\sigma|^2}\,i\partial S_{\varepsilon}\wedge \bar\partial S_{\varepsilon}\geq 0\,.
\end{eqnarray*}
Observe that the last form is smooth. Consequently, we find the 
inequalities
\begin{eqnarray*}
i\partial\bar\partial S_{\varepsilon}
&\geq &
\frac{\varepsilon}{|\sigma|^2}\,i\partial S_{\varepsilon}\wedge \bar\partial S_{\varepsilon} 
\,-\,
\frac{\{i{\cal C}_{E,h}\sigma,\sigma\}}{|\sigma|^2+\varepsilon}
\\
\\
&\geq &
\frac{\varepsilon}{|\sigma|^2}\,i\partial S_{\varepsilon}\wedge \bar\partial S_{\varepsilon} 
\,-\,
\Vert{\cal C}_{E,h}\Vert_{h,\omega}\,\frac{|\sigma|^2}{|\sigma|^2+\varepsilon}\,\omega
\end{eqnarray*}
where $\omega$ is a positive $(1,1)$-form.
\\
\\
{\bf Appendix B. Proof of the estimate \eqref{SecndOrdB} in Lemma \ref{YauIth}.} We will apply the computations of step (B) in the proof of Theorem \ref{HgRegMA} to the non-degenerate complex Monge-Amp\`ere equation 
$$
(\omega+i\partial\bar\partial\varphi'_j)^n=e^{h+L\varphi'_j-\varphi'_{j-1}}\,\omega^n\,.
$$
In this setting, the notation of setup (A) in the proof of the Theorem~\ref{HgRegMA} reduces to $\delta=l=h=0$, $\tilde\omega_{\varepsilon}=\omega$ and $i\partial\bar\partial h\geq -K_0\,\omega$. By replacing the term $f$ with $h-\varphi'_{j-1}$ in the expansion of the term $\sum_p\tilde A_{p,\bar p}/u_{p,\bar p}$ in step (B) in the proof of Theorem \ref{HgRegMA}, we infer
$$
0\geq e^{\frac{-L\varphi'_j-h+\varphi'_{j-1}}{n-1}}\,u_{n,\bar n}^{\frac{1}{n-1}}-\frac{(\varphi'_{j-1})_{n,\bar n}}{u_{n,\bar n}}-C'_0\,,
$$
Thus
\begin{eqnarray}\label{LAPYauIther}
0\geq C'_1\,u_{n,\bar n}^{\frac{1}{n-1}}-\frac{2n+\max_X \Delta_{\omega}\varphi'_{j-1}}{4u_{n,\bar n}}-C'_0\,,
\end{eqnarray}
by the estimates 
\begin{eqnarray}\label{IneqYauItherBis}
\varphi''_{0}\leq \varphi''_{j-1}\leq \varphi''_{j}\leq \varphi'_{j}\leq \varphi'_{j-1}\leq \varphi'_{0}\,.
\end{eqnarray}
This estimate implies also that at the maximum point $x_j$ we have
$$
u_{n,\bar n}(x_j)=\Lambda^{\omega}_{\varphi'_j}=e^{k\varphi'_j(x_j)}{\cal B}_j(x_j)\geq C'_2 B_j\,,
$$
with $B_j:=\max_X {\cal B}_j>0$. Then estimate \eqref{SecndOrdB} in Lemma \ref{YauIth} follows from \eqref{LAPYauIther} and the fact that 
$$
0<2n+\Delta_{\omega}\varphi'_j\leq 2ne^{k\max_X\varphi'_j}\,B_j\leq C\,B_j\,,
$$
which is itself a consequence of \eqref{IneqYauItherBis}.\hfill $\Box$
\\
\\
{\bf Appendix C. Relation with other works.}
As explained in the introduction the present work has its foundations in the papers \cite{Yau}, \cite{Be-Te} and especially in \cite{Kol1}, \cite{Kol2}.
A few months after that the first version of the present paper appeared on the arXiv server, P.~Eyssidieux, V.~Guedj, A.~Zeriahi posted on the same server a related preprint \cite{E-G-Z2}. In this preprint the authors obtain a weaker version of Statement (C) given in our Theorem \ref{Kolo}, which is sufficient to imply  Tian's conjecture as stated in \cite{Ti-Ko}. The statement in \cite{E-G-Z2} is weaker since it requires the (somehow stronger) assumption $\Omega/\omega^n\in L^{\varepsilon}(X)$, where $\omega\ge 0$ is smooth, big and degenerate. For the same reason a weaker version of Lemma \ref{comprCap} is stated in \cite{E-G-Z1}. 

At this point one should observe that the essence of the capacity method introduced in \cite{Kol1} does not allow to produce the required $L^{\infty}$-estimate in the case of a big and non nef class. It is possible to see that in this case the constants blow-up. This blow-up phenomenon has been one of the motivations of our work, which has led us to the proof of Tian's conjecture \cite{Ti-Ko}. Moreover fixed point methods do not produce a priori the $L^{\infty}$-estimate needed to construct singular K\"{a}hler-Einstein metrics and to investigate their regularity.

We wish to point out that in a quite recent preprint \cite{Di-Zh} the authors claim (in Theorem 1.1) boundedness and continuity of the solutions of some particular type of degenerate complex Monge-Amp\`{e}re equations. No proof of this claim seems to be provided. The authors also claim a stability result which is not sufficient to imply the continuity of solutions in the degenerate case. In fact a sequence of discontinuous functions converging in $L^{\infty}$-norm does not have necessarily a continuous limit$\,$! Moreover the same claim (Theorem 1.1) has been stated in \cite{Zh1}, \cite{Zh2}, but again no proof of continuity seems to be given (see page 12 in \cite{Zh1} and page 146 in \cite{Zh2}). The arguments for the boundedness of the solutions in \cite{Zh1}, \cite{Zh2} are quite informal in the degenerate case and seem impossible to follow. 

Concerning the stability of the solutions, the continuity assumption is quite natural and often available in the applications. In fact in the applications one works with smooth solutions provided by the Aubin-Yau solution of the Calabi conjecture with respect to variable  K\"{a}hler forms of type $\omega+\varepsilon\alpha$, as in the proof of theorem \ref{HgRegMA}. This perturbation process is one of the reasons of trouble for the continuity of the solutions. Moreover the stability with respect to the data $f$ considered in \cite{Di-Zh} is not essential in this context since one has $L^1$-compactness of quasi-plurisubharmonic functions normalized by the supremum condition. In fact a particular case of the stability result, namely Theorem \ref{Kolo} B, implies the continuity of the solution of the complex Monge-Amp\`{e}re equation $(\omega+i\partial\bar\partial \varphi)^n=e^{\lambda\varphi}f\,\Omega$, whenever $\omega>0$ is a K\"{a}hler metric and $f\in L\log^{n+\varepsilon}L$. This fact has been observed also in \cite{Kol2}.

Finally we mention that a nice and simple proof of the regularization of quasi-plurisubharmonic functions in the case of zero Lelong numbers can be found in \cite{Bl-Ko}.
%%%%%%%%%%%%%%%%%%%%%%%%%%%%%%%%%%%%%%%%%%%%%%%%%%%%%%%%%%%%%%%%%%%%%%%%%%%%%%%%%%%%%%%%%%%%%%%%%%%%%%%%%%%%%%%%%%%%%%%%%%%%%%%%%%%%%%%%%%%%%%%%%%%%%%%%%%%%%%%%%%%%%%%%%%%%%%%%%%%%%%%%%%%%%%%%%%%%%%%%%%%%%%%%%%%%%%%%
%%%%%%%%%%%%%%%%%%%%%%%%%%%%%%%%%%%%%%%%%%%%%%%%%%%%%%%%%%%%%%%%%%%%%%%%%%%%%%%%%%%%%%%%%%%%%%%%%%%%%%%%%%%%%%%%%%%%%%%%%%%%%%%%%%%%%%%%%%%%%%%%%%%%%%%%%%%%%%%%%%%%%%%%%%%%%%%%%%%%%%%%%%%%%%%%%%%%%%%%%%%%%%%%%%%%%%%%%%%%%%%%%%%%%%%%%%%%%%%%%%%%%%%%%%%%%%%%%%%%%%%%%%%%%%%%%%%%%%%%%%%%%%%%%%%%%%%%%%%%%%%%%%%%%%%%%%%%%%%%%%%%%%%%%%%%%%%%%%%%%%%%%%%%%%%%%%%%%%%%%%%%%%%%%%%%%%%%%%%%%%%%%%%%%%%%%%%%%%%%%%%%%%%%%%
\\
\\
{\bf Acknowledgments.} The second named author is grateful to Professor Gang Tian for bringing this type of problems to his attention. He also expresses his gratitude to the members of Institut Fourier for providing an excellent research environment. He thanks Adrien Dubouloz, Herv\'e Pajot, Olivier Labl\'{e}e, Eric Dumas and Gabriele La Nave for useful conversations.
The authors warmly thank the referee for detailed and constructive
criticism of the exposition of the manuscript.

\vspace{1cm}
\noindent
Jean-Pierre Demailly
\\
Universit\'{e} de Grenoble I, D\'epartement de Math\'ematiques
\\
Institut Fourier, 38402 Saint-Martin d'H\`{e}res, France
\\
E-mail: \textit{demailly@fourier.ujf-grenoble.fr}
\\
\\
Nefton Pali
\\
Universit\'{e} Paris Sud, D\'epartement de Math\'ematiques 
\\
B\^{a}timent 425 F91405 Orsay, France
\\
E-mail: \textit{nefton.pali@math.u-psud.fr}
\end{document}